\documentclass[a4paper]{article}
\usepackage{geometry}
\input{assiomi.tex}

\RequirePackage{amsmath}
\RequirePackage{amssymb}
\RequirePackage{amsthm}
\usepackage{epsfig}
\usepackage{graphicx}
\usepackage[dvips]{psfrag}
\usepackage[latin1]{inputenc}
\newcommand{\SSS}{\mathbb{S}}
\newcommand{\RR}{\mathbb{R}}
\newcommand{\NN}{\mathbb{N}}

\newcommand{\CC}{\mathbb{C}}
\newcommand{\eps}{\varepsilon}

\newcommand{\action}{{\mathcal A}}
\newcommand{\cont}{{\mathcal C}}
\newcommand{\dint}{\displaystyle \int}
\newcommand{\dlim}{\displaystyle \lim}
\newcommand{\ds}{\displaystyle}
\newcommand{\dt}{\, dt}

\newcommand{\uztag}{(U0)}

\newcommand{\uutag}{(U1)}

\newcommand{\udtag}{(U2)}

\newcommand{\udhtag}{(U2)$_{\mbox{h}}$}

\newcommand{\udltag}{(U2)$_{\mbox{l}}$}

\newcommand{\uthtag}{(U3)$_{\mbox{h}}$}

\newcommand{\utltag}{(U3)$_{\mbox{l}}$}

\newcommand{\uqhtag}{(U4)$_{\mbox{h}}$}

\newcommand{\uqltag}{(U4)$_{\mbox{l}}$}

\newcommand{\uctag}{(U5)}

\newcommand{\useitag}{(U6)}

\newcommand{\usetteltag}{(U7)$_{\mbox{l}}$}

\newcommand{\usettehtag}{(U7)$_{\mbox{h}}$}

\newcommand{\minus}{{\smallsetminus}}
\newcommand{\from}{\colon\thinspace}

\theoremstyle{plain}
\newtheorem{theorem}{Theorem}
\newtheorem{lemma}{Lemma}
\newtheorem{corollary}[lemma]{Corollary}
\newtheorem{proposition}[lemma]{Proposition}
\newtheorem*{painleve}{Painlev\'e's Theorem}
\newtheorem*{vonzeipel}{Von Zeipel's Theorem}

\theoremstyle{definition}
\newtheorem{remark}[lemma]{Remark}
\newtheorem{definition}[lemma]{Definition}
\newtheorem{example}{Example}

\numberwithin{equation}{section}
\numberwithin{lemma}{section}

\title{On the singularities of generalized solutions \\to $n$--body type problems%
\footnote{This work is partially supported by Italy MIUR, national project
``Variational Methods and Nonlinear Differential Equations''.}}
\author{Vivina Barutello\footnote{Supported by Istituto Nazionale di Alta
Matematica.}\and Davide L.~Ferrario \and Susanna Terracini%
\footnote{\protect\address}
}

\newcommand{\address}{{Universit\`a di Milano Bicocca,}\\
{Dipartimento di Matematica e Applicazioni,}\\
{via  Cozzi 53, 20125 Milano.}\\
\small{e-mail {\ttfamily vivina.barutello@unimib.it,davide.ferrario@unimib.it,susanna.terracini@unimib.it}.}
}

\begin{document}
%===========================
\maketitle
%===========================
\begin{abstract}
\noindent The validity of Sundman-type asymptotic estimates for collision solutions
 is established for a wide class of dynamical systems with singular forces, including the classical $N$--body problems with Newtonian, quasi--homogeneous and logarithmic potentials. The solutions are meant in the generalized sense of Morse (locally --in space and time-- minimal  trajectories with respect to compactly supported variations) and their uniform limits. The analysis includes the extension of the Von Zeipel's Theorem and the proof of isolatedness of collisions. Furthermore, such asymptotic analysis is applied to prove the absence of collisions for  locally minimal trajectories.
\end{abstract}
%===========================
\noindent 2000 {\it Mathematics Subject Classification.} 
\section{Introduction}\label{sec:intro}
%===========================
%\cite{DPS1,DPS2,DS,diacu_book}
Many systems of interacting bodies of interest in Celestial and other areas of classical 
Mechanics have the form
\begin{equation} \label{DS_intro}
m_i\ddot x_i = \frac{\partial U}{\partial x_i}(t,x), \qquad i=1,\ldots,n
\end{equation}
where the  forces $\frac{\partial U}{\partial x_i}$ are undefined 
on a singular set $\Delta$. This is for example the set of collisions
between two or more particles in the $n$--body problem.  Such singularities play
a fundamental role in the phase portrait (see, e.g. \cite{diacu_book}) and strongly influence the global orbit
structure, as they can be held responsible, among others,
of the  presence of chaotic motions (see, e.g. \cite{Dev_survey}) and of motions
becoming unbounded in a finite time \cite{MM,xia}.

Two are the major steps in the analysis of the impact of the singularities 
in the $n$--body problem:  the first consists in performing the asymptotic
analysis along a single collision  (total or partial)  trajectory  and
goes back, in the classical case, to the works by  Sundman (\cite{sundman}), 
Wintner (\cite{wintner}) and, in
more recent years by Sperling, Pollard, Saari, Diacu and other authors
(see for instance \cite{polsaa1,polsaa2,saa72,sperling,elbialy,diacureg}).
The second step consists in blowing--up the singularity by a
suitable change of coordinates introduced by McGehee in \cite{mcgehee1} 
and replacing it by an invariant
boundary --the collision manifold-- where the flow can be extended in a
smooth manner. It turns out that, in many 
interesting applications, the flow on the collision manifold has a simple
structure: it is a gradient--like, Morse--Smale flow featuring a few stationary
points and heteroclinic connections (see, for instance, the surveys \cite{Dev_survey, moeckel88}).
The analysis of the extended flow allows us to obtain a full picture of the
behavior of solutions 
near the singularity, despite the flow fails to be fully regularizable (except
in a few cases).

The geometric approach, via the McGehee coordinates and the
collision manifold, 
can be successfully applied also to obtain asymptotic estimates
in some cases, such as  the 
collinear three--body problem (\cite{mcgehee1}), the anisotropic Kepler problem 
(\cite{Dev1,Dev2,DPS2,DS}),
the three--body problem both in the planar isosceles case (\cite{Dev3}) and  the
full perturbed three-body, as described in \cite{Dev_survey,diacuab}.
Besides the quoted cases, however, one needs
to establish the asymptotic estimates before
blowing--up the singularity, in order to prove convergence of the blow--up
family. The reason is quite 
technical and mainly rests in the fact the a singularity of the $n$--body
problems needs not  be
isolated, for the possible occurrence of  partial collisions in a 
neighborhood of the total collision. 
In the literature, this problem has been usually overcame by extending the flow
on partial collisions via some regularization technique (such as Sundman's, in \cite{Dev3}, 
or Levi--Civita's in \cite{LC}). Such a
device works well only when partial collisions are binary, which are the only
singularities to be globally removable.
Thus, the extension of the geometrical analysis to the full $n$--body problems
finds a strong theoretical obstruction:
partial collisions must be regularizable,  what is known to hold true only in few cases.
Other interesting cases in which the geometric method is not effective are that  of
quasi--homogeneous potentials
(where there is a lack of regularity for the extended flow)
and that of logarithmic potentials (for the failure of the blow--up technique).

In this paper we extend the classical asymptotic estimates near collisions in three
main directions.
\begin{enumerate}
\item 
We take into account of a very general notion of solution for the dynamical
system~\eqref{DS_intro},
which fits particularly well to solutions found by variational techniques. Our notion of
solution includes, besides all classical noncollision trajectories, all the  
\emph{locally minimal solutions} (with respect to compactly supported variations) 
which are often termed  minimal the sense of Morse. Furthermore, we include in the set of 
\emph{generalized solutions} all the limits of classical and locally minimal solutions.
\item We extend our analysis to a wide class of potentials including not only
homogeneous and quasi-homogeneous potentials, but also those with weaker
singularities of logarithmic type.
\item We allow potentials to strongly depend  on time
(we only require its time derivative 
to be controlled by the potential itself -- see assumption \ref{a:U1}). 
In this way, for instance, we 
can take into account models where masses vary in time.
\end{enumerate}

Our main results on the asymptotics near total collisions (at the origin) are Theorems 
\ref{central_conf} and~\ref{central_conf2} (for quasi-homogeneous potentials) and 
Theorems~\ref{central_conf_log} and~\ref{central_conf2_log} 
(when the potential is of logarithmic type) which extend the classical 
Sundman--Sperling asymptotic estimates (\cite{sundman,sperling}) in the directions above
(see also \cite{diacuab,DPS1}).

As a consequence of the asymptotic estimates, the
presence of a total collision prevents the occurrence of partial ones for
neighboring times.

This observation plays a central role when extending the
asymptotic estimates to the full $n$--body problem  since it allows us to reduce
from partial (even simultaneous) collisions to total ones by decomposing the
system in colliding clusters. Our results also lead to the extension of 
the concept of singularity for the dynamical system~\eqref{DS_intro} to the 
class of generalized solutions.
We shall prove an extension of the Von Zeipel's Theorem:
when the moment of inertia is bounded then every singularity of a generalized solution
admits a limiting configuration, hence all singularities are collisions.
The results on total collisions are then fully extended 
to partial ones in Theorem~\ref{thm:partial}.

A further motivation for the study of generalized solutions 
comes from the variational approach to the study of 
selected  trajectories to the $n$-body problem.
Indeed the exclusion of collisions is a major problem in the 
application of variational techniques as it results in the recent
literature, where many different arguments have been introduced to prove that the trajectories found in such a way
are collisionless (see \cite{amco93,ABT,BFT,BT2004,BC-Z,MR2032484,chenkyoto,CM,F3,FT,montgomery96,montgomery_contmath,mont_prepr,serter1,terven,andrea_thesis}).
As a  first application we shall be able to extend some of these techniques in order to 
prove that action minimizing trajectories are free of collisions for a wider class of 
interaction potentials. For example in the case of quasi--homoge\-neous potentials,
once collisions are isolated, the blow-up technique can be successfully applied 
to prove that locally minimal solutions are, in many circumstances, {\em free of collisions}.
In order to do that we can use the method 
of {\em averaged variations} introduced by Marchal and developed in \cite{marchal,Ch2,FT}. 
It has to be noticed that, when dealing with logarithmic--type potentials, the blow-up technique
is not available since  converging
blow-up sequences do not exists; we can anyway prove that the average 
over all possible variations is negative by taking advantage of the
harmonicity of the function $\log|x|$ in $\RR^2$. With this result we can then 
extend to quasi--homogeneous and logarithmic potentials all the analysis 
of the (equivariant) minimal trajectories carried in \cite{FT}.

Besides the direct method, other variational techniques --Morse and minimax theory--
have been applied to the search of periodic solutions in singular problems
(\cite{amco93,bahrab,mate95,riahi}). In the quoted papers, however 
only the case of \emph{strong force} interaction (see \cite{G}) has been
treated. Le us consider a  sequence of solutions to penalized problems
where an infinitesimal sequence of strong force terms is added to the potential:
then its limit enjoys the same conservation laws as the generalized solutions. Hence our
main results apply also to this class of trajectories.
We believe that our study can be usefully applied to
develop a Morse Theory that takes into account 
the topological contribution of collisions.
Partial results in this direction are given in \cite{BS,riahi2}, where the contribution
 of collisions to the Morse Index is computed.

The paper is organized as follows:

\tableofcontents

%===========================
\section{Singularities of locally minimal solutions}\label{sec:sing}
%===========================
%% __HERE__ 
\subsection{Locally minimal solutions}\label{sec:locminsol}
%===========================
We fix a metric on the \emph{configuration space} $\RR^k$
and we denote by $I(x)=|x|^2$ the \emph{moment of inertia} associated 
to the configuration $x \in
\RR^{k}$ and
\[
\mathcal{E}:=\{ x \in \RR^{k} : |x|^2 = 1 \}
\]
the inertia ellipsoid.
We define the radial and ``angular'' variables associated to $x \in \RR^{k}$ as
\begin{equation}
\label{rs}
r := |x| = I^{\frac12}(x) \in [0,+\infty), \quad
s := \frac{x}{|x|} \in {\mathcal E}.
\end{equation}

We consider the dynamical system 
\begin{equation} \label{DS}
\ddot x = \nabla U(t,x),
\end{equation}
on the time interval $(a,b)\subset \RR$, $-\infty \leq a < b \leq +\infty$.
Here $U$ is a positive time-dependent potential function 
$U \from   (a,b) \times \left(\RR^{k}\minus \Delta\right)  \rightarrow \RR^+$,
and it is supposed to be of class $\mathcal C^1$ on  its domain; 
by $\nabla U$ we denote  its gradient
with respect to the given metric.
\begin{remark}
In the case of $n$-body type systems as described in~\eqref{DS_intro}, given $m_1,\ldots,m_n$, $n \geq 2$ positive real numbers, 
we define the scalar product induced by the {\it mass metric} 
on the \emph{configuration space} $\RR^{nd}$ between $x=(x_1,\ldots,x_n)$ 
and $y=(y_1,\ldots,y_n)$, as 
\begin{equation} \label{mass_prod}
x \cdot y = \sum_{i=1}^{n}m_i \langle x_i,y_i \rangle,
\end{equation}
where $\langle\cdot,\cdot\rangle$ is the scalar product in $\RR^{d}$.
We denote by $|\cdot|$ the norm induced by the mass scalar
product~\eqref{mass_prod}. Then 
$\nabla U(t,x)$ denotes the gradient of the potential, in the mass metric, with respect to the spatial variable $x$,
that is:
\[
\nabla U(t,x)=M^{-1} \frac{\partial U}{\partial x}(t,x),
\]
where $\left(\frac{\partial U}{\partial x}\right)_i=\frac{\partial U}{\partial x_i}$,
$i=1,\ldots,n$, and $M=[M_{ij}]$, 
$M_{ij}=m_i\delta_{ij}{\mathbf 1}_d$ (${\mathbf 1}_d$ is the $d$-dimensional
identity matrix) for every $i,j=1,\ldots,n$.
\end{remark}

Furthermore we suppose that $\Delta$ is a singular set for $U$ of an attractive
type, in the sense that
\begin{assiomi}
\item[\uztag]\label{a:U0}
$\displaystyle\lim_{x \to \Delta} U(t,x) =
+\infty$, uniformly in $t$.
\end{assiomi}
Borrowing the terminology from the study of the singularities of the $n$--body problem, 
the set $\Delta$ will be often referred as \emph{collision set} and it is
required to be a \emph{cone}, that is\label{cone}
\[
x \in \Delta \quad \implies \quad \lambda x\in\Delta, \quad \forall \lambda\in
\RR.
\]
We observe that being a cone implies that  $0 \in \Delta$. 
When $x(t^*) \in \Delta$ for some $t^*\in(a,b)$ we will say that $x$ has
an \emph{interior collision} at $t=t^*$ and that $t^*$ is a
\emph{collision instant} for $x$. When $t^* = a$ or $t^* = b$ (when
finite) we will talk about a \emph{boundary collision}. In particular, if
$x(t^*)=0\in \Delta$, we will say that $x$ has a \emph{total collision at the origin} at
$t=t^*$.
A collision instant $t^*$ is termed
\emph{isolated} if there exists $\delta>0$ such that, for every
$t \in (t^*-\delta,t^*+\delta)\cap (a,b)$, $x(t)\notin \Delta$.

We consider the following assumptions on the potential $U$:

\begin{assiomi}
\item[\uutag] 
\label{a:U1}
\text{There exists a constant $C_1 \geq 0$ such that, for every $(t,x) \in (a,b)
\times\left( \RR^{k}\minus \Delta\right)$},
\[
\left|\frac{\partial U}{\partial t}(t,x)\right| \leq C_1 \left(U(t,x)+1\right).
\]
\item[\udtag] 
\label{a:U2}
\text{There exist constants $\tilde \alpha \in (0,2)$ and $C_2 \geq 0$ such
that}
\[
\nabla U(t,x)\cdot x + \tilde \alpha U(t,x) \geq - C_2.
\]
\end{assiomi}

% \begin{remark}
% When $U$ is homogeneous of degree $-\alpha$ in its second
% variable (for instance when $U(t,x)=1/|x|^\alpha$), condition \ud\,
% is satisfied and the equality is achieved with $\tilde \alpha=\alpha$ and $C_2=0$. 
% Furthermore, \ud\, is satisfied when $U(t,x)=\log(1/|x|)$ for every value of\enspace
% $\tilde \alpha$.
% \end{remark}

We then define the \textit{lagrangian action functional} on the interval $(a,b)$
as
\begin{equation}
\label{action}
\action(x,[a,b]) := \int_a^b K(\dot{x}) + U(t,x) \dt,
\end{equation}
where 
\begin{equation}
\label{kin}
K(\dot{x}) :=  \frac{1}{2} |\dot{x}|^2,
\end{equation}
is the \textit{kinetic energy}. We observe that $\action(\cdot,[a,b])$ is
bounded and ${\cal C}^2$ on the Hilbert space $H^1\left( (a,b),\RR^{k}\minus
\Delta \right)$. In terms of the variables $r$ and $s$ introduced
in~\eqref{rs}, 
the action functional reads as
\[
\action(rs,[a,b]):=\int_a^b\frac12\left(\dot r^2 + r^2|\dot
s|^2\right)+U(t,rs)\dt
\]
and the corresponding Euler--Lagrange equations, whenever $x \in H^1\left(
(a,b),\RR^{k}\minus \Delta \right)$, are
\begin{equation}
\label{eq:EL}
\begin{split}
-\ddot r + r|\dot s|^2 + \nabla U(t,r s) \cdot s & = 0 \\
-2r\dot r \dot s -r^2 \ddot s + r \nabla_{T} U(t,r s) & = \mu s
\end{split}
\end{equation}
where $\mu= r^2|\dot s|^2$ is the Lagrange multiplier due to the presence of the
constraint $|s|^2 = 1$ and the vector $\nabla_{T} U(t,r s)$ is the tangent
components to the ellipsoid ${\mathcal E}$ of the gradient $\nabla U(t,r s)$, that is
$\nabla_{T} U(t,rs)=\nabla U(t,r s)-\nabla U(t,r s)\cdot s$.

\begin{definition}\label{def:locally_minimal_sol}
A path $x \in H^1_{loc}\left((a,b),\RR^{k}\right)$ is a \emph{locally minimal
solution} for the dynamical system~\eqref{DS} if,
 for every $t_0 \in (a,b)$, there exists $\delta_0 >0$ such that the restriction
of $x$ to the interval $I_0=[t_0-\delta_0,t_0+\delta_0]$, is a local minimizer
for
$\action(\cdot,I_0)$ with respect to compactly supported variations (fixed-ends).
\end{definition}

\begin{remark}\label{rem:3}
We observe that a priori a locally minimal solution $x$ can have a large collision 
set, $x^{-1}(\Delta)$; this set, though of Lebesgue measure zero,  can very well 
admit many accumulation points. For this
reason the Euler--Lagrange equations~\eqref{eq:EL} and the dynamical 
system~\eqref{DS} do not hold for a locally minimal trajectory.
\end{remark}
\begin{remark}
When the potential is of class ${\cal C}^2$ outside $\Delta$ then every classical
noncollision solution in the interval $(a,b)$ is a a locally minimal solution.
\end{remark}
\begin{definition}\label{def:generalized_sol}
A path $x$ is a \emph{generalized solution} for the dynamical system~\eqref{DS} 
if there exists a sequence $x_n$ of locally minimal solutions such that
\begin{enumerate}
\item $x_n\to x$ uniformly on compact subsets of $(a,b)$;
\item for almost all  $t \in (a,b)$ the associated total energy
$h_n(t) := K(\dot x_n(t))-U(t,x_n(t))$ converges.
\end{enumerate}
\end{definition}

We say that a (classical, locally minimal, generalized) solution $x$ on the interval $(t_1,t_2)$,
has a {\em singularity} at  $t_2$ (finite) if it is not possible to extend $x$ as a (classical,
locally minimal, generalized) solution to a larger interval $(t_1,t_3)$ with $t_3>t_2$.

In the framework of classical solutions to $n$-body systems, 
the classical Painlev\'e's Theorem  (\cite{painleve,Dpainleve}) asserts
that the existence of a singularity
at a finite time $t^*$ is equivalent to the fact that the minimal 
of the mutual distances becomes infinitesimal as $t\to t^*$.
This fact reads as:

\begin{painleve}\label{teo:Painleve}
Let $\bar x$ be a classical solution for the $n$-body dynamical system
on the interval $[0,t^*)$. If $\bar x$ has a singularity at $t^*<+\infty$, then the
potential associated to the problem diverges to $+\infty$ as $t$ approaches $t^*$.
\end{painleve}

Painlev\'e's Theorem does not necessarily imply that a collision (i.e. that is a singularity
such that all mutual distances have a definite limit) occurs when
there is a singularity at a finite time; indeed this two facts are equivalent 
only if each particle approaches a definite configuration
(on this subject we refer to \cite{polsaa1,polsaa2,saa72}).
This result has been stated by Von Zeipel in 1908 (see \cite{vonzeipel} and also \cite{mcgehee2})
and definitely proved by Sperling in 1970 (see \cite{sperling}): 
in the $n$-body problem the occurrence of singularities (in finite time) 
which are not collisions is then equivalent to the existence of an unbounded motion.

\begin{vonzeipel}\label{teo:Zeipel}
If $\bar x$ is a classical solution for the $n$-body dynamical system
on the interval $[0,t^*)$ with a singularity at $t^*<+\infty$ and
$\lim_{t \to t^*}I(\bar x(t))<+\infty$, then $\bar x(t)$ has a definite limit
configuration $x^*$ as $t$ tends to $t^*$.
\end{vonzeipel}
We will come back later on the proof of this result (in Corollary~\ref{cor:2}
and in Section~\ref{sec:partial}). To our purposes, we give the following definition.

\begin{definition}\label{def:singularity}
We say that the (generalized) solution $\bar x$ for 
the dynamical system~\eqref{DS} has a {\em singularity} at $t=t^*$ if
\[
\lim_{t \to t^*}U(t,\bar x(t))=+\infty.
\]
\end{definition}

\begin{definition}\label{def:collision}
The singularity $t^*$ is a said to be a \emph{collision} 
for the locally minimal 
solution $\bar x$ if it admits a limit configuration as $t$ tends to $t^*$.
\end{definition}

%=======================================================
\subsection{Approximation of locally minimal solutions}
%=======================================================
Let $\bar x$ be a locally minimal solution on the interval $(a,b)$ and let
$I_0\subset(a,b)$ be an interval 
such that $\bar x$ is a (local) minimizer for $\action(\cdot,I_0)$ with respect
to compactly supported variations. Generally local minimizers need not to be
isolated; we illustrate below a penalization argument to select a particular
solution from the possibly large set of local minimizers.
To begin with, we define the auxiliary functional on the space $H^1\left( I_0,
\RR^{k}\right)$
\begin{equation}
\label{baraction}
\bar \action(x,I_0) := 
\int_{I_0} K(\dot{x}) + U(t,x) +\dfrac{|x-\bar x|^2}{2}\dt.
\end{equation}
When the interval $I_0$ is sufficiently
small, $\bar x$ is actually the global minimizer for the penalized functional
$\bar \action(\cdot,I_0)$
defined in~\eqref{baraction}. Of course we may assume that
\begin{equation}\label{bddA}
\bar \action(\bar x,I_0) = \action(\bar x,I_0) < + \infty,
\end{equation}
which is equivalent to require that $\bar \action(\cdot,I_0)$ takes a finite
value at least at one point.

\begin{proposition} \label{propo:globalmin}
Let $\bar x$ be a locally minimal solution on the interval $(a,b)$, let
$\delta_0>0$ and $t_0 \in (a,b)$ be such that $\bar x$ is a local minimizer for
$\action(\cdot,I_0)$ where  $I_0=[t_0-\delta_0,t_0+\delta_0]\subset(a,b)$.
Then there exists $\bar \delta = \bar \delta(\bar x)>0$ such that whenever
$\delta_0 \leq \bar \delta$, $\bar x$ is the unique global minimizer for $\bar
\action(\cdot,I_0)$.
\end{proposition}

\begin{proof}
For every $x \in H^1_{loc}\left( I_0, \RR^{k}\right)$ the 
inequality $\action(x,I_0)
\leq \bar \action(x,I_0)$ holds true,
and it is an equality only if $x=\bar x$. Since
$\bar x$ is a local minimizer for $\action(x,I_0)$, one easily infers, by a
simple convexity argument, the existence of $\eps > 0$ such that
\[\|x-\bar x\|_\infty < \eps\quad\implies\quad\action(\bar x,I_0) \leq
\action(x,I_0).\]
We conclude that, for every $x \in H^1_{loc}\left( I_0, \RR^{k}\right)$,
such that $0<\|x-\bar x\|_\infty < \eps$,  the following chain 
of inequalities holds
\[
\bar \action(\bar x,I_0) = \action(\bar x,I_0) \leq \action(x,I_0) < \bar
\action(x,I_0);
\]
hence $\bar x$ is a strict local minimizer for $\bar\action(\cdot,I_0)$,
independently on 
$\delta_0$.

In order to complete the proof we show that $\bar \action(\bar x,I_0) < \bar
\action(x,I_0)$ also for those functions 
$x \in H^1_{loc}\left( I_0, \RR^{k}\right)$ such that $\|x-\bar x\|_\infty \geq
\eps$, provided $\delta_0$ is sufficiently small. Indeed, since  the Sobolev
space $H^1_{loc}\left( I_0, \RR^{k}\right)$ is embedded in the space of
absolutely continuous functions, we can compute, by H\"older inequality,
\begin{equation} \label{eq:dis}
\begin{split}
|(x-\bar x)(t)|&\leq \int_{I_0} |\dot x(s)|ds + \int_{I_0}|\dot{\bar x}(s)|ds \\
               &\leq \sqrt{2\delta_0} \left( \sqrt{\int_{I_0} |\dot x(s)|^2ds} 
                     + \sqrt{\int_{I_0}|\dot{\bar x}(s)|^2ds}\right).
\end{split}
\end{equation}
By taking the supremum at both sides of~\eqref{eq:dis} it follows that 
\[
\dfrac{\|x-\bar x\|_\infty}{\sqrt{2\delta_0}} - \sqrt{\int_{I_0}|\dot{\bar
x}(s)|^2\, ds} \leq
\sqrt{\int_{I_0} |\dot x(s)|^2\, ds},
\]
and therefore, for every $x \in H^1\left( I_0, \RR^{k}\right)$,
\begin{equation} \label{eq:dis2}
\begin{split}
\bar \action(x,I_0) \geq \int_{I_0} |\dot x(s)|^2\, ds 
&\geq \left( \dfrac{\|x-\bar x\|_\infty}{\sqrt{2\delta_0}}
             -\sqrt{\int_{I_0}|\dot{\bar x}(s)|^2\, ds}\right)^2\\
&\geq \left( \dfrac{\eps}{\sqrt{2\delta_0}} -\sqrt{\int_{I_0}|\dot{\bar
x}(s)|^2\, ds}\right)^2
\end{split}
\end{equation}
Hence, by choosing $\delta_0$ such that 
$2\delta_0 < \eps\left(\sqrt{\int_{I_0}|\dot{\bar x}(s)|^2ds}+\sqrt{\bar
\action(\bar x,I_0)}\right)^{-2}$, it follows that 
\[
\bar \action(x,I_0) \geq \left( \dfrac{\eps}{\sqrt{2\delta_0}} - 
\sqrt{\int_{I_0}|\dot{\bar x}(s)|^2ds}\right)^2 > \bar \action(\bar x,I_0) 
\]
also for those paths $x \in H^1_{loc}\left( I_0, \RR^{k}\right)$ such that
$\|x-\bar x\|_\infty \geq \eps$.
This concludes the proof.
\end{proof}

We now wish to approximate the singular potential $U$ with a family of smooth
potentials $U_{\eps}\from (a,b)\times \RR^{k} \to \RR^+$, 
depending on a parameter
$\eps > 0$.
To this aim consider  the function
\[
\eta(s) = 
\begin{cases}
s & \text{if $s\in [0,1]$} \\
\dfrac{-s^2 +6s - 1}{4} & \text{if $s\in [1,3]$} \\
2 & \text{if $s \geq 3$;}
\end{cases}
\]
notice that
$\eta\in\cont^{1}\left(\RR^+,\RR^+\right) $ and,  
% \begin{equation*}
% \label{eq:eta}
% \eta(s)=s \, \text{ when } \, s \in \left[0,{1}\right], \qquad \eta(s)=2 \,
% \mbox{ when } s \geq 3
% \end{equation*}
for every $s \in [0,+\infty)$,
\begin{equation*}
% \label{eta_ineq}
\dot\eta(s) \, s \leq \eta(s) \qquad \text{and} \qquad \dot\eta(s)\leq 1.
\end{equation*}
Now let us define, for $\eps >0$,
\[
\eta_\eps(s) := \frac{1}{\eps}\eta(\eps s);
\]
then the following inequalities hold for every $s \in [0,+\infty)$
\begin{equation}\label{etaeps_ineq}
\dot\eta_\eps(s) \, s \leq \eta_\eps(s), \qquad \text{and} \qquad
\dot\eta_\eps(s)\leq 1.
\end{equation}
By means of the family $\eta_\eps$ we can regularize the potential $U$ in the
following way: 
\begin{equation}
\label{def:Ueps}
U_{\eps}(t,x) = \left\{
\begin{array}{lll}
\displaystyle
\eta_\eps\left(U(t,x)\right), &   & \mbox{if } x \in \RR^{k}\minus
\Delta,\\ 
{2}/{\eps},                   &   & \mbox{if } x \in \Delta.
\end{array} \right.
\end{equation}
It is worthwhile to understand that each $U_\eps(t,x)$ coincides with $U(t,x)$
whenever $U(t,x) \leq {1}/{\eps}$; in fact
\[
\eta_\eps(s)=\frac{1}{\eps}\eta(\eps s)=s
\]
whenever $\eps s \in [0,1]$, that is $s \in [0,1/\eps]$.
Next, we  consider the associated family of boundary value problems on the
interval $I_0 \subset (a,b)$
\begin{equation}
\label{DSbar}
\left\{
\begin{array}{l}
\ddot x = \nabla U_\eps(t,x) + (x-\bar x), \\
x\vert_{ \partial I_0 } = \bar x \vert_{\partial I_0} 
\end{array}\right.
\end{equation}
where, as usual, $\nabla U_\eps(t,x)$ is the gradient, in the mass metric, with respect
to the spatial variable $x$. Solutions of~\eqref{DSbar} are critical points of the
action functional
\begin{equation}
\label{baractioneps}
\bar \action_\eps(x,I_0) := \int_{I_0} K(\dot{x}) + U_\eps(t,x) +\dfrac{|x-\bar
x|^2}{2}\dt.
\end{equation}
We observe that $\bar \action_\eps(\cdot,I_0)$ is bounded and $\cont^2$ on
$H^1_{loc}\left( I_0, \RR^{k}\right)$, since $U_\eps$ is smooth on the whole
$\RR^{k}$.
We also remark that the infimum of $\bar \action_\eps(\cdot,I_0)$ is achieved,
for $\bar \action_\eps(\cdot,I_0)$ is a positive and coercive functional on 
$H^1_{loc}\left( I_0, \RR^{k}\right)$.

In the next proposition we prove that a locally minimal solution 
has the fundamental property to be the limit of a
sequence of global minimizers for the approximating functionals $\bar
\action_\eps(\cdot,I_0)$, provided the interval $I_0 \subset (a,b)$ is chosen so
small  that the restriction of the minimal solution to $I_0$ is the unique
global minimizer for $\bar\action(\cdot,I_0)$.  This result is crucial, indeed, as observed in 
Remark~\ref{rem:3}, the Euler--Lagrange equations and the dynamical system 
hold for a locally minimal solution; we will anyway be able to use the ones corresponding
to the approximating global minimizers (for the regularized problems) to prove the fundamental
properties of locally minimal (and generalized) solutions in the rest of the paper.

\begin{proposition}\label{propo:conv}
Let $\bar x$ and $I_0$ be given by Proposition~\ref{propo:globalmin}.
Let $\epsilon>0$ and $x_\eps$ be a global minimizer for $\bar
\action_{\eps}(\cdot,I_0)$. Then, up to subsequences, as $\eps \to 0$,
\begin{enumerate}
\item[(i)] $U_{\eps}(t,x_{\eps}) \to U(t,\bar x)$ almost everywhere and in $L^1$;
\item[(ii)] $x_{\eps}\to \bar x$ uniformly; 
\item[(iii)] $\dot x_{\eps}\to \dot{\bar x}$ in $L^2$;
\item[(iv)] $\dot x_{\eps}\to \dot{\bar x}$ almost everywhere;
\item[(v)] $\dfrac{\partial U_{\eps}}{\partial t}(t,x_{\eps}) \to \dfrac{\partial
U}{\partial t}(t,\bar x)$ almost everywhere and in $L^1$.
\end{enumerate}
\end{proposition}
\begin{proof}
As we have already observed, for every $\eps>0$, the potential $U_\eps$
coincides with $U$ on the sublevel $\{(t,x):U(t,x)\leq 1/\eps\}$ and, by its
definition, for every $(t,x) \in I_0\times \RR^{k}\minus\Delta$
\[
U_\eps(t,x) \leq U(t,x).
\]
Therefore 
\[
\bar \action_\eps(x,I_0) \leq \bar \action(x,I_0)
\]
for every $x \in H^1_{loc}\left(I_0,\RR^{k}\right)$.
It follows from~\eqref{bddA} that
\begin{equation}\label{bddineq}
\bar \action_\eps(x_\eps,I_0) = \inf_{x \in H^1_{loc}} \bar \action_\eps(x,I_0)
\leq \bar \action(\bar x,I_0) < +\infty,
\end{equation}
which implies the boundedness of the family $\left\{\int_{I_0}|\dot{x}_\eps|^2
+|x_\eps-\bar x|^2\right\}_\eps$.
Hence we deduce the existence of a sequence $\left(x_{\eps_n} \right)_{\eps_n}
\subset (x_\eps)_\eps$ such that $\left(\dot x_{\eps_n} \right)_{\eps_n}$
converges weakly in $L^2$ and uniformly to some limit $\tilde x$. In addition we
observe that 
\[
\lim_{\eps_n \to 0}U_{\eps_n}(t,x_{\eps_n}(t)) = U(t,\tilde x(t))
\]
for every $t \in I_0$, regardless the finiteness of $U(t,\tilde x(t))$.\\
From~\eqref{bddineq} we also deduce the boundedness of the following integrals 
\[
\int_{I_0}U_{\eps_n}(t,x_{\eps_n}) \,dt \leq \bar \action_{\eps_n}(x_{\eps_n},I_0)
< +\infty.
\]
and therefore, since the sequence $\left( U_{\eps_n}(t,x_{\eps_n})
\right)_{\eps_n}$
is positive, by applying Fatou's Lemma  one deduces that 
\[
\int_{I_0} U (t,\tilde x) \leq \liminf \int_{I_0}
U_{\eps_n}(t,x_{\eps_n})<+\infty.
\]
Hence from the weak semicontinuity of the norm in $L^2$ (the sequence
$(\dot x_{\eps_n})_{\eps_n}$ converges weakly in $L^2$ to $\tilde x$) we obtain the
inequalities
\begin{equation*}
\bar \action (\tilde x,I_0) \leq \liminf \bar \action_{\eps_n}(x_{\eps_n},I_0)
\leq \bar \action(\bar x,I_0)
\end{equation*}
which contradict Proposition~\ref{propo:globalmin}, unless $\tilde x = \bar x$
and
\begin{equation}\label{disA}
\liminf \bar \action_{\eps_n}(x_{\eps_n},I_0) = \bar \action(\bar x,I_0).
\end{equation}
Therefore we deduce the $L^2$-convergence of the sequence $\left(\dot x_{\eps_n}
\right)_{\eps_n}$
and its convergence almost everywhere to $\dot{\bar x}$, up to subsequences.
From~\eqref{disA} it follows  also that
\begin{equation}\label{conv_int}
\lim_{\eps_n \to 0}\dint_{I_0} U_{\eps_n}(t,x_{\eps_n}) =\dint_{I_0}U(t,\bar x).
\end{equation}
From the convergence almost everywhere of $(U_{\eps_n}(t,x_{\eps_n}))_{\eps_n}$
together with~\eqref{conv_int} we conclude its convergence in $L^1$ to $U(t,\bar
x)$.

We now turn to the convergence of the sequence
$(\varphi_n(t))_{\eps_n}=\left(\dfrac{\partial U_{\eps_n}}{\partial
t}(t,x_{\eps_n})\right)_{\eps_n}$. To this aim, we observe that condition \ref{a:U1}
together with~\eqref{etaeps_ineq} imply the following chain of inequalities
\[
\begin{split}
\left|\frac{\partial U_{\eps_n}}{\partial t}(t,x_{\eps_n}(t))\right| 
& = \dot\eta_{\eps_n}\left(U(t,x_{\eps_n}(t))\right)\left|\frac{\partial
U}{\partial t}(t,x_{\eps_n}(t))\right| \\
& \leq C_1 \dot\eta_{\eps_n}\left( U(t,x_{\eps_n}(t))\right) 
\left(U(t,x_{\eps_n}(t))+1\right)\\
& \leq C_1 \left(\eta_{\eps_n}\left( U(t,x_{\eps_n}(t))\right)+1\right) \\
& = C_1 \left(U_{\eps_n}(t,x_{\eps_n}(t))+1\right).
\end{split}
\]
We already know that $U_{\eps_n}(t,x_{\eps_n}(t))$ converges in $L^1$. 
This implies the finiteness almost everywhere of $(\varphi_n(t))_{\eps_n}$ and
hence its almost everywhere 
convergence is due to the uniform convergence of
$(x_{\eps_n})_{\eps_n}$.
We obtain the $L^1$ convergence of $(\varphi_n(t))_{\eps_n}$ to $\dfrac{\partial
U}{\partial t}(t,\bar x)$ from the Dominated Convergence Theorem.
\end{proof}

%==============================================
\subsection{Conservation laws}
\label{subsec:conlaws}
%==============================================
Now the sequence of solutions to the regularized problems are used to prove
the conservation of the energy for locally minimal solutions.

\begin{proposition}\label{propo:convenergy}
Let $\bar x$ and $I_0$ be given by Proposition~\ref{propo:globalmin}. 
Then the energy associated to\enspace$\bar x$
\begin{equation}
\label{eq:energy}
h\from I_0 \to \RR, \qquad h(t) := K(\dot{\bar x}(t)) - U(t,\bar x(t)) 
\end{equation}
is of class $W^{1,1}$ on $I_0$ and its weak derivative is
\[
\dot h(t)= \frac{\partial U}{\partial t}(t,\bar x).
\]
\end{proposition}
\begin{proof}
Let $(x_{\eps})_{\eps}$ be the sequence of global minimizers for the
corresponding functionals $\bar \action_{\eps}(\cdot,I_0)$ convergent to $\bar
x$ whose existence is proved in Proposition~\ref{propo:conv}.
Let $h_{\eps}$ be the energy associated to $x_{\eps}$, that is
\begin{equation}
\label{eq:energyeps}
h_{\eps} \from  I_0 \to \RR,  \qquad 
h_{\eps}(t) := K({\dot x}_{\eps}(t))-U_{\eps}(t,x_{\eps}(t))+\frac{1}{2}|\bar
x(t)- x_{\eps}(t)|^2
\end{equation}
From Proposition~\ref{propo:conv} we immediately deduce that the sequence
$(h_{\eps})_{\eps}$ converges pointwise to $h(\bar x)$; moreover 
from~\eqref{bddineq} and~\eqref{eq:energyeps} 
we have
\[
\int_{I_0} |h_{\eps}(t)|\dt \leq \bar\action_{\eps}\left(x_{\eps},I_0\right) <
\bar\action\left(\bar x,I_0\right).
\]
From the Dominated Convergence Theorem we obtain that the sequence
$(h_{\eps})_{\eps}$ converges in $L^1$ to the integrable function $h$.

We still have to prove that $h$ admits weak derivative. To this end, let us
consider a test function $\varphi \in \cont_0^\infty(I_0)$; we can  write
\[
\begin{split}
\int_{I_0}h(t)\dot \varphi(t)\dt 
& = \lim_{\eps \to 0} \int_{I_0}h_{\eps}(t)\dot \varphi(t)\dt \\
& = \lim_{\eps \to 0} -\int_{I_0}\frac{\partial U_{\eps}}{\partial
t}(t,x_{\eps}(t))\varphi(t)\dt.
\end{split}
\]
In consequence of Proposition~\ref{propo:conv}, the sequence $\left(\frac{\partial
U_{\eps}}{\partial t}(t,x_{\eps}(t))\right)_\eps$ converges to
$\frac{\partial U}{\partial t}(t,\bar x(t))$ in $L^1$; then
\begin{equation}\label{limdU}
\lim_{\eps \to 0} \int_{I_0}\frac{\partial U_{\eps}}{\partial t}(t,x_{\eps}(t))
\varphi(t)\dt
=  \int_{I_0}\frac{\partial U}{\partial t}(t,\bar x(t)) \varphi(t)\dt, \qquad 
\forall \varphi \in \cont_0^\infty(I_0)
\end{equation}
and hence
\[
\int_{I_0}h(t)\dot \varphi(t) \dt 
= - \int_{I_0}\dfrac{\partial U}{\partial t}(t,\bar x)\varphi(t) \dt, \qquad 
\forall \varphi \in \cont_0^\infty(I_0)
\]
which means that 
$\frac{\partial U}{\partial t}(t,\bar x)$ is the weak derivative of
$h(\bar x)$.
\end{proof}
The next corollary follows straightforwardly.
\begin{corollary}\label{cor:en_abs_cont}
The energy associated to a locally minimal solution on the interval $(a,b)$ is
in the Sobolev space $W^{1,1}_{loc}\left((a,b),\RR\right)$.
\end{corollary}

We now investigate the behavior of the moment of inertia of a locally minimal solution
when a singularity occurs
(see Definition~\ref{def:singularity}). The results contained 
in Proposition~\ref{propo:Iconvex} and Corollary~\ref{cor:I''>0} are the natural extension
of the classical Lagrange--Jacobi inequality to locally minimal solutions (see \cite{wintner}).

\begin{proposition}\label{propo:Iconvex}
Let $\bar x$ be a locally minimal solution and $I_0$ be given by Proposition
\ref{propo:globalmin}.
Then 
\begin{equation} \label{disconv}
\frac{1}{2}\int_{I_0}I(\bar x(t))\ddot \varphi (t)\dt \geq 
\int_{I_0} \left[2h(\bar x(t))+(2-\tilde\alpha)U(t,\bar x(t))-C_2\right]\varphi (t)\dt
\end{equation}
for every $\varphi \in \cont_0^\infty\left(I_0,\RR\right)$, $\varphi(t)\geq 0$.
\end{proposition}
\begin{proof}
Let $(x_{\eps})_{\eps}$ be the sequence of global minimizers for the
corresponding functionals $\bar \action_{\eps}(\cdot,I_0)$ convergent to $\bar
x$ whose existence is proved in Proposition~\ref{propo:conv}. 
When we compute the second derivative of the moment of inertia of $x_{\eps}$ we
obtain
\[
\begin{split}
\frac{1}{2}\ddot I(x_{\eps}(t)) & = |\dot x_{\eps}(t)|^2 + \ddot x_{\eps}(t)
\cdot x_{\eps}(t)\\
&= 2h_{\eps}(t)+2U_{\eps}(t,x_{\eps}(t))-|\bar x(t)-x_{\eps}(t)|^2 \\
&\hspace*{3.5cm} + \left[\nabla U_{\eps}(t,x_{\eps}(t))+(x_{\eps}(t)-\bar
x(t))\right]\cdot x_{\eps}(t)\\
&= 2h_{\eps}(t)+2U_{\eps}(t,x_{\eps}(t))+\bar x(t)\cdot(x_{\eps}(t)-\bar x(t)) \\ 
&\hspace*{3.5cm} + \dot\eta_{\eps}\left(U(t,x_{\eps})\right) \nabla
U(t,x_{\eps}(t)) \cdot x_{\eps}(t)
\end{split}
\]
hence, by assumption \ref{a:U2} on the potential $U$ and
inequality~\eqref{etaeps_ineq}, it follows that 
\begin{equation}\label{dis:I''eps}
\begin{split}
\frac{1}{2}\ddot I(x_{\eps}(t)) &\geq 2h_{\eps}(t)+2U_{\eps}(t,x_{\eps}(t)) +
\bar x(t)\cdot(x_{\eps}(t)-\bar x(t))\\
&\hspace*{4cm}-\dot\eta_{\eps}\left(U(t,x_{\eps})\right)\left[\alpha
U(t,x_{\eps})+C_2\right]\\
&\geq 2h_{\eps}(t) + (2-\tilde\alpha)U_{\eps}(t,x_{\eps}(t))
+\bar x(t)\cdot(x_{\eps}(t)-\bar x(t)) -C_2
\end{split}
\end{equation}
for some $\tilde \alpha\in (0,2)$ and $C_2>0$.
Therefore, since $x_{\eps} \in \cont^2(I_0)$, for every $\varphi \in
\cont_0^\infty\left(I_0,\RR\right)$, $\varphi(t)\geq 0$
\[
\frac{1}{2}\int_{I_0}I(x_{\eps}(t))\ddot \varphi (t)\dt =
\frac{1}{2}\int_{I_0}\ddot I(x_{\eps}(t))\varphi (t)\dt
\]
and, from~\eqref{dis:I''eps},
\[
\begin{split}
&\frac{1}{2}\int_{I_0}I(x_{\eps}(t))\ddot \varphi (t)\dt \\ 
& \hspace{1cm} \geq 
\int_{I_0} \left[ 2h_{\eps}(t)+(2-\tilde\alpha)U_{\eps}(t,x_{\eps}(t))+\bar
x(t)\cdot(x_{\eps}(t)-\bar x(t))-C_2 \right]\varphi (t)\dt.
\end{split}
\]
We conclude by passing to the limit as $\eps \to 0$ in~\eqref{dis:I''eps} and
using the $L^1$-convergences proved in Propositions~\ref{propo:conv} and
\ref{propo:convenergy}.
\end{proof}

The next corollaries follow directly.
\begin{corollary}[\textbf{Lagrange--Jacobi inequality}]\label{cor:I''>0}
Let $\bar x$ be given by Proposition~\ref{propo:globalmin}.
Then the following inequality holds in the distributional sense
\[
\frac{1}{2}\ddot I(\bar x(t)) \geq 
2h(t)+(2-\tilde \alpha)U(t,\bar x(t))-C_2, \qquad \forall t \in (a,b).
\]
\end{corollary}

\begin{corollary}\label{cor:14}
Let $\bar x$ be given by Proposition~\ref{propo:globalmin}.
Then its moment of inertia is convex on $I_0$ whenever $\bar x$ has a singularity
in $t_0$ and $\delta_0$ is small enough.
\end{corollary}
\begin{proof}
Whenever $\eps$ and $\delta_0$ are sufficiently small, the right hand side 
of inequality~\eqref{dis:I''eps} is
strictly positive, indeed $h_{\eps}(t)$ is bounded, $x_\eps$ converges to $\bar
x$ uniformly and $U_{\eps}(t,x_{\eps}(t))$ diverges to $+\infty$. Whenever
$\eps$ is small enough we conclude that
\[
\ddot I(x_{\eps}(t)) > 0
\]
and hence $I(x_{\eps})$ are strictly convex functions in a neighborhood of
$t_0$. Since the sequence $I(x_{\eps})$ uniformly converges to $I(\bar x)$ we
conclude that  also $I(\bar x)$ is convex on the interval\enspace$I_0$.
\end{proof}

We now investigate 
the possibility that a sequence of singularities accumulates
at the right bound of the interval $(a,b)$; in this section we will
suppose that $b<+\infty$.

\begin{lemma}
\label{le:1}
Let $\bar x$ be given in Proposition~\ref{propo:globalmin}, $h$ be its energy
defined in~\eqref{eq:energy} and fix $\tau \in (a,b)$ be such that 
\begin{equation}\label{eq:lambda}
\lambda := \frac{2-\tilde\alpha}{2}-C_1(b-\tau)
\end{equation}
is a strictly positive constant. Then there exists a constant $K>0$ such that 
\begin{equation}
\label{eq:lemma1}
\left|\int_{\tau}^{t} h(s) ds\right| \leq \left(
\frac{2-\tilde\alpha}{2}-\lambda\right) \int_{\tau} ^{t} U(s,\bar x(s))ds + K, \qquad
\forall t \in (\tau,b).
\end{equation}
\end{lemma}
\begin{proof}
Since $h$ is absolutely continuous on every interval $[\tau,t] \subset (a,b)$
(Corollary~\ref{cor:en_abs_cont}) we have
\[
|h(t)| \leq |h(\tau)| + \int_{\tau}^t |\dot{h}(\xi)| \,d\xi, \qquad \forall t \in
(\tau,b).
\]
From Proposition~\ref{propo:convenergy} and assumption \ref{a:U1} we obtain
\begin{equation}
\label{bddh}
\begin{split}
|h(t)|
& \leq |h(\tau)| + \int_{\tau}^t \left|\frac{\partial U}{\partial
\xi}(\xi,\bar x(\xi))\right| \,d\xi \\
& \leq |h(\tau)| + C_1\int_{\tau}^t \left(U(\xi,\bar x(\xi))+1\right) \,d\xi
\end{split}
\end{equation}
and integrating both sides of the inequality on the interval $[\tau,t]$
\[
\int_{\tau}^{t}|h(s)|ds \leq |h(\tau)|(t-\tau)+C_1\frac{(t-\tau)^2}{2}
+ C_1\int_{\tau}^t ds \int_{\tau}^s U(\xi,\bar x(\xi)) \,d\xi.
\]
Since $U$ is positive, the integral $\int_{\tau}^s U(\xi,\bar x(\xi)) \,d\xi$
increases in the variable $s$, hence we conclude
\[
\left|\int_{\tau}^{t} h(s) ds\right| \leq \int_{\tau}^{t}|h(s)|ds \leq K +
C_1(b-\tau)\int_{\tau}^t U(\xi,\bar x(\xi)) \,d\xi.
\]
where $K:=|h(\tau)|(t-\tau)+C_1 (t-\tau)^2/2$.
\end{proof}
\begin{lemma}
\label{le:2}
Let $\bar x$ be given in Proposition~\ref{propo:globalmin} and $\tau$ be chosen
as in Lemma~\ref{le:1}. Suppose that there exist $\delta,C>0$ such that
\[
I(\bar x(t)) \leq C, \mbox{ for every } t \in (b-\delta,b) 
\]
and 
\begin{equation}
\label{eq:liminf}
\liminf_{t\to b^-} \dot I(\bar x(t)) \leq C.
\end{equation}
Then there exists $\tau \in (a,b)$ such that
\[
\int_{\tau} ^{b} U(t,\bar x(t)) \dt < +\infty.
\]
\end{lemma}
\begin{proof}
If $b$ is not a singularity for $\bar x$,
the assertion follows from assumption~\eqref{bddA}.
Otherwise, 
it follows from \eqref{eq:liminf} that 
% since $\liminf_{t\to b^-} \dot I(\bar x(t)) \leq C$ 
there exists
an increasing sequence $(t_n)_n$ such that
\[
t_n \to b \mbox{ as } n\to+\infty\quad \mbox{and} \quad  \dot I(\bar x(t_n)) \leq C, \forall n.
\]
Now let $N$ be an integer
such that $t_N \in (b-\delta,b)$ and the constant $\lambda$
defined in~\eqref{eq:lambda}, with $\tau = t_N$, is strictly positive.
Hence, for every index $n>N$,
\[
2C \geq \dot I(\bar x(t_n)) - \dot I(\bar x(t_N))
= \int_{t_N}^{t_n} \ddot I(\bar x(t)) \dt.
\]
Corollary~\ref{cor:I''>0} implies that
\[
C \geq 2 \int_{t_N}^{t_n} h(t) \dt + (2-\tilde\alpha)\int_{t_N}^{t_n} U(t,\bar x(t))
\dt - C_2(t_n-t_N).
\]
We now apply Lemma~\ref{le:1} to deduce that 
\begin{equation}
\label{eq:end_le2}
2\lambda \int_{t_N}^{t_n} U(t,\bar x(t)) \dt \leq C_2(t_n-t_N) + C + 2K.
\end{equation}
Since $\lambda>0$ is fixed, as $n \to +\infty$ the proof is completed.
\end{proof}
\begin{corollary} \label{cor:2} 
Let $\bar x$ be a generalized solution on $(a,b)$. 
Suppose that 
\begin{equation}\label{eq:limitazioni}
\limsup_{t\to b^-}I(\bar x(t))<+\infty,\qquad\mbox{and}\qquad
\liminf_{t\to b^-} \dot I(\bar x(t)) <+\infty.
\end{equation} 
Then, if $-\infty<a<\tau<b<+\infty$ there hold
\begin{enumerate}
\item[(i)] $\displaystyle \int_{\tau}^{b} U(t,\bar x(t)) \dt < +\infty$;
\item[(ii)] $\displaystyle \left|\int_{\tau}^{b} h(t) \dt \right| < +\infty$;
\item[(iii)] $\displaystyle \int_{\tau}^{b} K(\dot{\bar x}(t)) \dt < +\infty$;
\item[(iv)] $\displaystyle \|h\|_\infty < +\infty$ on $[\tau,b)$.
\item[(v)] $\displaystyle \lim_{t\to b^-} \bar x(t)$ exists.
\end{enumerate}
\end{corollary}
\begin{proof}
We first prove the assertions in the case of locally minimal solutions.
The boundedness of the first integral follows from the assumption of local
boundedness of the action functional on a locally minimal trajectory,
assumption~\eqref{bddA}, and from Lemma~\ref{le:2}. Concerning the second one, we use
Corollary~\ref{cor:en_abs_cont} and inequality~\eqref{eq:lemma1}; {\it (iii)}
follows straightforwardly from {\it (i)}, {\it (ii)} and the definition of the energy 
$h$.
The boundedness of $\| h\|_\infty$ on $(a,b)$ follows from Corollary
\ref{cor:en_abs_cont} and inequality~\eqref{bddh}.
To deduce {\it (v)} it is sufficient to remark that, from {\it (iii)},
$\bar x$ is  H\"{o}lder-continuous on $(a,b)$.

In order to extend the proof to
generalized solutions, we first remark that all the constants and bounds
appearing in the 
proof 
above do not depend on the specific solution $\bar x$, but only on the potential
and the limits in equations~\ref{eq:limitazioni}, 
and the total energy
valued at single instant $h(\tau)$ of the interval. Hence the assertions {\it (i)}, {\it (ii)}
still hold true when passing to pointwise limit such that the energy $h(\tau)$ is bounded. 
The other assertions then follow from the fist two.
\end{proof}
\begin{remark}
In Corollary~\ref{cor:2} {\it (v)}, 
Von Zeipel's Theorem is proved, for generalized solutions, under the additional 
assumption~\eqref{eq:liminf}. The proof will be completed in Section~\ref{sec:partial}.
\end{remark}
\begin{remark}
From inequality~\eqref{bddh} we can easily understand that in Definition
\ref{def:generalized_sol} the convergence of the energy of the approximating sequence of
locally minimal solutions can be assumed only at one point.
\end{remark}

%=======================================================
\section{Asymptotic estimates at total collisions}\label{sec:asymp}
%=======================================================

The purpose of this section is to deepen the analysis of the asymptotics of
generalized solutions as they approach a total collision at the origin; to this aim, 
we introduce some further hypothesis on the potential $U$.
Though we will perform all the analysis in a left neighborhood of the collision
instant, the analysis concerning right neighborhoods is the exact analogue.

We recall that $\bar x$ has a {\emph total collision at the origin} at $t=t^*$
if $\lim_{t\to t^*}\bar x(t)=x(t^*)=0$. Since by our assumptions $0$ belongs to the singular set 
$\Delta$ of the potential, assumption~\ref{a:U0} reads that a total collision instant 
is a singularity for $\bar x$.

The results proved in \S~\ref{subsec:conlaws} have some relevant consequences
in the case of total collisions at the origin; in particular Corollary~\ref{cor:14} now reads
\begin{corollary}\label{cor:Icont}
Let $\bar x$ be given by Proposition~\ref{propo:globalmin}.
If $|\bar x(t_0)|=0$, then there exists $\delta_0>0$ such that $I(\bar x)$ is
continuous on $I_0=[t_0-\delta_0,t_0+\delta_0]$, it admits weak derivative
almost everywhere, the function $\dot I(\bar x)$ is monotone increasing and
$\dot I(\bar x) \in BV(I_0)$.
Furthermore the following inequalities hold in the distributional sense
\[
\begin{aligned}
\ddot I(\bar x(t)) > 0 & \qquad    \forall t \in \bar I_0 \\
\dot I(\bar x(t)) < 0 &  \qquad  \forall t \in (t_0-\delta_0,t_0)\\
\dot I(\bar x(t)) > 0 &  \qquad  \forall t \in (t_0,t_0+\delta_0).
\end{aligned}
\]
\end{corollary}

Furthermore, since $I(\bar x(t))\geq 0$, and 
$I(\bar x(t^*)) = 0$  if and only if $\bar x$ has a total collision at the origin at $t=t^*$, from
Corollary~\ref{cor:14} one can 
deduce that, whenever a total collision occurs at $t=t_0$, 
no other total collisions take place in the interval $I_0$.
Concerning the occurrence of total collision at the boundary of the interval 
$(a,b)$, we argue as in Lemmata~\ref{le:1} and~\ref{le:2} and we use
the convexity of the function $I$ to deduce that also boundary total collisions are isolated.
It is worthwhile noticing that this fact does not prevent, at this stage, the occurrence of 
infinitely many other singularities in a neighborhood of a total collision at the origin.
We summarize these remarks in the next theorem.
\begin{theorem}\label{thm:totcolliso}
Let $\bar x$ be a generalized solution for the dynamical system~\eqref{DS}.
Suppose that $-\infty<a<b<+\infty$ and that there exists $t_0 \in [a,b]$ such
that $|\bar x(t_0)|=0$. Then there exists $\delta>0$ such that, for every
$t\in(t_0-\delta,t_0+\delta)\cap[a,b]$, $t\neq t_0$,
we have $|\bar x(t)|\neq 0$.
\end{theorem}

In terms of the radial variable $r$ Corollary~\ref{cor:Icont} and Theorem
\ref{thm:totcolliso} state that whenever $r(t_0)=0$ for some $t_0\in (a,b)$ ($t_0$
can coincide with $a$ or $b$ when finite) then there exists $\delta >0$ such
that
\begin{equation}\label{eq:isol_coll}
\begin{split}
& r(t) \neq 0, \qquad \dot r(t) < 0, \qquad \forall t \in (t_0-\delta,t_0)\\
& r(t) \neq 0, \qquad \dot r(t) > 0, \qquad \forall t \in (t_0,t_0+\delta).
\end{split}
\end{equation}
We moreover rewrite the bounded energy function as
\begin{equation}
\label{eq:energy_rs}
h(t)=\frac{1}{2}\left( \dot r^2 + r^2|\dot s|^2 \right) - U(t,rs).
\end{equation}
Similarly, denoting by 
$(x_\eps)_\eps$ the sequence of global minimizers for $\bar
\action_\eps(\cdot,[t_0-\delta,t_0+\delta])$ converging to the locally minimal
collision solution $\bar x$ whose existence is proved in Proposition
\ref{propo:conv},  we define, for every $\eps$, 
\[
r_\eps := |x_\eps| \in \RR \qquad \text{and} \qquad s := \frac{x_\eps}{|x_\eps|}
\in {\mathcal E}
\]
and we write the energy in~\eqref{eq:energyeps} as 
\[
h_{\eps}(t)=\frac{1}{2}\left( \dot r_\eps^2 + r_\eps^2|\dot s_\eps|^2 \right) -
U_\eps(t,r_\eps s_\eps)+\frac{1}{2}|rs - r_\eps s_\eps|^2.
\]
Furthermore the approximating action functional and the corresponding
Euler--Lagrange equations in the new variables are respectively 
\[
\bar \action_\eps(r_\eps s_\eps,[t_0-\delta,t_0+\delta]) :=
\int_{t_0-\delta}^{t_0+\delta} \frac{1}{2}\left( \dot r_\eps^2 + r_\eps^2|\dot
s_\eps|^2 \right) + U_\eps(t,r_\eps s_\eps)+\frac{1}{2}|rs - r_\eps s_\eps|^2
\dt
\]
and 
\begin{equation}
\label{eq:ELeps}
\begin{split}
-\ddot r_\eps + r_\eps|\dot s_\eps|^2 + \nabla U_\eps(t,r_\eps s_\eps) \cdot
s_\eps -(rs - r_\eps s_\eps)\cdot s_\eps & = 0 \\
-2r_\eps\dot r_\eps \dot s_\eps -r_\eps^2 \ddot s_\eps + r_\eps \nabla_{T}
U_\eps(t,r_\eps s_\eps) -r_\eps(rs - r_\eps s_\eps)& = \mu_\eps s_\eps,
\end{split}
\end{equation}
where $\mu_\eps = r_\eps^2|\dot s_\eps|^2-r_\eps(rs - r_\eps s_\eps)\cdot
s_\eps$ is the Lagrange multiplier due to the presence of the constraint
$|s_\eps|^2 = 1$ and the vector $\nabla_{T} U_\eps(t,r_\eps s_\eps)$ is the
tangent components to the ellipsoid ${\mathcal E}$ of the gradient $\nabla
U_\eps(t,r_\eps s_\eps)$. A similar approximation procedure will be implicitly
done for generalized solutions.

\bigskip

To proceed with the analysis of the asymptotic behavior near total collisions at the origin we need 
some stronger conditions on the potential $U$ when the radial variable $r$ tends to 0.
These additional conditions includes quasi--homogeneous potential and logarithmic ones, 
in the following analysis, however we will treat separately the two different cases.

%----------------------------------------
\subsection{Quasi-homogeneous potentials}
%----------------------------------------
In this section we impose some stronger assumptions on the behavior of the potential
when $|x|$ is small. The following conditions are trivially satisfied by $\alpha$-homogeneous 
potentials and mimic the behavior of combination of such homogeneous potentials:
\begin{assiomi}
\item[\udhtag]
\label{a:U2h}
There exist $\alpha \in (0,2)$, $\gamma>l$ and $C_2 \geq 0$ such that
\[
\nabla U(t,x)\cdot x + \alpha U(t,x) \geq - C_2 |x|^\gamma U(t,x),
\]
whenever $|x|$ is small.
\end{assiomi}

% \begin{remark}
% We observe that when $U$ is homogeneous of degree $-\alpha$ in its second
% variable (for instance when $U(t,x)=1/|x|^\alpha$), the equality in condition
% \udh\, is achieved with $C_2 =0$; this assumption is satisfied also when we take
% into account potentials of the form $U(t,x)=U_\alpha(x)+U_\beta(x)$ where
% $U_\alpha$ is homogeneous of degree $-\alpha$, $U_\beta$ is homogeneous of
% degree $-\beta$, $0<\beta < \alpha$, and both  $U_\alpha$ and $U_\beta$ are
% positive. Indeed in this case  condition \udh\, is verified (with the strict
% inequality) with $\gamma = C_2 =0$.
% On the other hand when we consider a negative homogeneous singular perturbation
% $U_\beta$ (we think for instance to the potential $U(t,x)=1/|x|^{\alpha} -
% 1/|x|^{\beta}$, with $0<\beta < \alpha$) we have that \udh\, holds, when $|x|$
% is sufficiently small, with $C_2=\alpha-\beta$ and $0<\gamma < \alpha-\beta$.
% \end{remark}
\begin{remark}
\ref{a:U2h} implies \ref{a:U2} (for small values of $|x|$); in fact, 
by choosing $2>\tilde\alpha>\alpha>0$, one obtains
\[
\begin{split}
\nabla U(t,x)\cdot x + \tilde\alpha U(t,x) & =
\nabla U(t,x)\cdot x + \alpha U(t,x) + (\tilde\alpha-\alpha) U(t,x) \\
& \geq - C_2 |x|^\gamma U(t,x)+(\tilde\alpha-\alpha) U(t,x),
\end{split}
\]
and the last term remains bounded below as $|x|\to0$ since
$\tilde\alpha-\alpha>0$.
\end{remark}

Furthermore we suppose the existence of a function $\tilde U$ defined and of
class $\cont^1$ on
$(a,b)\times \left({\mathcal E}\minus \Delta\right)$ such that
\begin{equation}
\label{tildeU}
\inf_{(a,b)\times \left({\mathcal E}\minus \Delta\right)}\tilde U(t,s)  >0
\quad \text{and} \quad \lim_{s\to {\mathcal E}\cap\Delta}\tilde U(t,s)=+\infty \qquad
\mbox{uniformly in } t.
\end{equation}
The potential $U$ is then supposed to verify the following condition uniformly
in the variables $t$ and $s$ (on the compact subsets of $(a,b)\times
\left({\mathcal E}\minus \Delta\right)$):

\begin{assiomi}
\item[\uthtag]
\label{a:U3h}
$\dlim_{r \to 0} r^\alpha U(t,x) = \tilde U(t,s)$.
\end{assiomi}

\begin{remark}
In \ref{a:U2h} and \ref{a:U3h} the value of $\alpha$ must be the same. We shall refer to potentials satisfying such assumptions as quasi--homogeneous (cf.~\cite{diacuab}).
\end{remark}

\begin{lemma}
\label{le:integ}
Let $\bar x$ be a generalized solution, let $t_0 \in (a,b]$ be a total
collision instant and let $\delta$ be given in Theorem~\ref{thm:totcolliso}. Then, 
for every $\alpha'\in(\alpha,2)$, we have 
\[
\int_{t_0-\delta}^{t_0}-r^{\alpha'}\frac{\dot r}{r}U(t,rs)\dt < +\infty,
\]
where $\alpha \in (0,2)$ is the constant fixed in assumption \ref{a:U2h}.
\end{lemma}
\begin{proof}
We consider the function
\[
\Gamma_{\alpha'} (t) := r^{\alpha'} \left( \dfrac{1}{2}r^2|\dot s|^2 - U(t,rs)
\right), 
\quad \alpha' \in (\alpha,2);
\]
Replacing in~\eqref{eq:energy_rs} we have 
\[
\Gamma_{\alpha'} (t) = h(t)r^{\alpha'} -\dfrac{1}{2}\dot r^2 r^{\alpha'} \leq
h(t)r^{\alpha'};
\]
since $h$ is bounded (see Corollary~\ref{cor:2}, {\em(iv)}) and $r$ tends to $0$, we
conclude that the function $\Gamma_{\alpha'}$ is bounded above on the interval
$[t_0-\delta,t_0]$. 
We consider the corresponding functions (still bounded above) for the
approximating problems:
\[
\begin{split}
\Gamma_{\alpha',\eps} (t) &= r_\eps^{\alpha'} \left( \dfrac{1}{2}r_\eps^2|\dot
s_\eps|^2 - U_\eps(t,r_\eps s_\eps) +\frac{1}{2}|rs-r_\eps s_\eps|^2\right)\\
&= h_\eps(t)r_\eps^{\alpha'} -\dfrac{1}{2}\dot r_\eps^2 r_\eps^{\alpha'} \leq
h_\eps(t)r_\eps^{\alpha'},
\end{split}
\]
and we observe that the sequence $\left(\Gamma_{\alpha',\eps}\right)_\eps$
converges almost everywhere and $L^1$ to $\Gamma_{\alpha'}$, as $\eps \to 0$.
We compute the derivative of $\Gamma_{\alpha',\eps} (t)$ with respect to  time as
\begin{equation}
\begin{split}
\label{eq:dGamma1eps}
&\frac{d}{dt}\Gamma_{\alpha',\eps} (t)= 
 \frac{2+\alpha'}{2}r_\eps^{1+\alpha'}\dot r_\eps |\dot s_\eps|^2\\
&\qquad\qquad\qquad + r_\eps^{2+\alpha'} \dot s_\eps \cdot \ddot s_\eps 
+ \alpha'r_\eps^{\alpha'-1}\dot r_\eps\left[\frac12|rs-r_\eps s_\eps|^2 -
U_\eps(t,r_\eps 
  s_\eps)\right]\\
& \;\; - r_\eps^{\alpha'}\left[ \frac{\partial U_\eps}{\partial t}(t,r_\eps
s_\eps) 
+ \nabla U_\eps(t,r_\eps s_\eps)(\dot r_\eps s_\eps+r_\eps\dot s_\eps) \right] 
+ r_\eps^{\alpha'}(rs-r_\eps s_\eps)\frac{d}{dt}(rs-r_\eps s_\eps).
\end{split}
\end{equation}
Now 
we multiply the Euler--Lagrange equation~\eqref{eq:ELeps}$_2$ by $\dot
s_\eps$ to obtain 
(we recall that $\nabla_T U_\eps(t,r_\eps s_\eps)\cdot \dot s_\eps = \nabla
U_\eps(t,r_\eps s_\eps)\cdot \dot s_\eps$ since $s_\eps$ and $\dot s_\eps$ are orthogonal)
\begin{equation}
\label{eq:EL2s'}
r_\eps^2 \ddot s_\eps\cdot \dot s_\eps = - 2r_\eps\dot r_\eps |\dot s_\eps|^2 +
r_\eps \nabla U_\eps(t,r_\eps s_\eps)\cdot \dot s_\eps - r_\eps r s\cdot s_\eps.
\end{equation}
Replacing~\eqref{eq:EL2s'} in~\eqref{eq:dGamma1eps} we have
\begin{multline}
\label{eq:gamma'3eps}
\frac{d}{dt}\Gamma_{\alpha',\eps}(t) = -\frac{2-\alpha'}{2} r^{1+\alpha'} \dot
r_\eps |\dot s_\eps|^2 - \alpha' r_\eps^{\alpha'-1} \dot r_\eps U_\eps(t,r_\eps
s_\eps) 
- r_\eps^{\alpha'} \frac{\partial U_\eps}{\partial t}(t,r_\eps s_\eps)\\
- r_\eps^{\alpha'-1} \dot r_\eps \nabla U_\eps(t,r_\eps s_\eps)\cdot (r_\eps
s_\eps)
- r^{\alpha'+1}_\eps r s \cdot s_\eps  \\
+ \frac{\alpha'}{2}r^{\alpha'-1}_\eps\dot r_\eps|rs-r_\eps s_\eps|^2
+ r_\eps^{\alpha'}(rs-r_\eps s_\eps)\frac{d}{dt}(rs-r_\eps s_\eps).
\end{multline}
We now combine assumptions \ref{a:U1}, \ref{a:U2h} and~\eqref{etaeps_ineq} to obtain the
following inequalities
\[
\begin{split}
- r_\eps^{\alpha'} \frac{\partial U_\eps}{\partial t}(t,r_\eps s_\eps) 
& = - r_\eps^{\alpha'}\dot\eta_\eps(U(t,r_\eps s_\eps))\frac{\partial
U}{\partial t}(t,r_\eps s_\eps)\\
& \geq - C_1 
r_\eps^{\alpha'}\dot\eta_\eps(U(t,r_\eps s_\eps))\left(U_\eps(t,r_\eps s_\eps)+1\right)\\
& \geq - C_1r_\eps^{\alpha'}\left(\eta_\eps(U(t,r_\eps s_\eps))+1\right)\\
& \geq -C_1 r_\eps^{\alpha'} \left(U_\eps(t,r_\eps s_\eps)+1\right),
\end{split}
\]
\[
\begin{split}
- r_\eps^{\alpha'} \frac{\dot r_\eps}{r_\eps}\nabla U_\eps(t,r_\eps s_\eps)\cdot
(r_\eps s_\eps)& =
- r_\eps^{\alpha'} \frac{\dot r_\eps}{r_\eps} \dot\eta_\eps(U(t,r_\eps
s_\eps))\nabla U(t,r_\eps s_\eps)\cdot (r_\eps s_\eps) \\
& \geq - r_\eps^{\alpha'} \frac{\dot r_\eps}{r_\eps} \dot\eta_\eps(U(t,r_\eps
s_\eps))U(t,r_\eps s_\eps)\left[-\alpha-C_2r_\eps^\gamma\right] \\
& \geq r_\eps^{\alpha'} \frac{\dot r_\eps}{r_\eps} U_\eps(t,r_\eps
s_\eps)\left[\alpha+C_2r_\eps^\gamma \right].
\end{split}
\]
Finally, by replacing in~\eqref{eq:gamma'3eps}, we obtain 
\[
\frac{d}{dt}\Gamma_{\alpha',\eps}(t) \geq \Psi_{\alpha',\eps}(t)
\]
where
\[
\begin{split}
\Psi_{\alpha',\eps}(t) &= 
-\frac{2-\alpha'}{2} r_\eps^{1+\alpha'} \dot r_\eps |\dot s_\eps|^2
-(\alpha'-\alpha) r_\eps^{\alpha'} \frac{\dot r_\eps}{r_\eps} U_\eps(t,r_\eps
s_\eps)
- C_1 r_\eps^{\alpha'} \left(U_\eps(t,r_\eps s_\eps)+1\right) \\
& + C_2 r_\eps^{\alpha'+\gamma} \frac{\dot r_\eps}{r_\eps} U_\eps(t,r_\eps
s_\eps)
- r^{\alpha'+1}_\eps r s \cdot s_\eps\\
& + \frac{\alpha'}{2}r^{\alpha'-1}_\eps\dot r_\eps|rs-r_\eps s_\eps|^2
+ r_\eps^{\alpha'}(rs-r_\eps s_\eps)\frac{d}{dt}(rs-r_\eps s_\eps).
\end{split}
\]
Since $\gamma>0$ and $r_\eps \to 0$ as $t \to t_0$, for every $\eps >0$, there
exists a positive $\lambda_\eps\leq (\alpha'-\alpha)/2$ such that
\[
C_2r_\eps^{\alpha'+\gamma} U_\eps(t,r_\eps s_\eps) \leq 
\lambda_\eps r_\eps^{\alpha'} U_\eps(t,r_\eps s_\eps)
\]
whenever $\delta$ is small enough; furthermore, since $- \frac{2-\alpha'}{2}
r_\eps^{1+\alpha'} \dot r_\eps |\dot s_\eps|^2$ is positive, we have
\[
\begin{split}
\Psi_{\alpha',\eps}(t) & \geq
- (\alpha'-\alpha-\lambda_\eps) r_\eps^{\alpha'} \frac{\dot r_\eps}{r_\eps}
U_\eps(t,r_\eps s_\eps) 
- C_1 r_\eps^{\alpha'} \left(U_\eps(t,r_\eps s_\eps)+1\right) \\
&- r^{\alpha'+1}_\eps r s \cdot s_\eps
+ \frac{\alpha'}{2}r^{\alpha'-1}_\eps\dot r_\eps|rs-r_\eps s_\eps|^2
+ r_\eps^{\alpha'}(rs-r_\eps s_\eps)\frac{d}{dt}(rs-r_\eps s_\eps).
\end{split}
\]
Therefore, for every $\eps$, the function $\Psi_{\alpha',\eps}$ is larger than 
the sum of a positive term 
\[
- (\alpha'-\alpha-\lambda_\eps) r_\eps^{\alpha'}
\frac{\dot r_\eps}{r_\eps} U_\eps(t,r_\eps s_\eps),
\]
an integrable term 
\[
- C_1 r_\eps^{\alpha'} \left(U_\eps(t,r_\eps s_\eps)+1\right)
\]
and
a remainder
\[
- r^{\alpha'+1}_\eps r s \cdot s_\eps
+ \frac{\alpha'}{2}r^{\alpha'-1}_\eps\dot r_\eps|rs-r_\eps s_\eps|^2
+ r_\eps^{\alpha'}(rs-r_\eps s_\eps)\frac{d}{dt}(rs-r_\eps s_\eps)
\]
converging uniformly to $- r^{\alpha'+2}$ as $\eps$ tends to 0.

We can then conclude that, for every $\eps$
\begin{multline*}
\Gamma_{\alpha',\eps}(t_0) - \Gamma_{\alpha',\eps}(t_0-\delta) \geq
-\int_{t_0-\delta}^{t_0}(\alpha'-\alpha-\lambda_\eps) r_\eps^{\alpha'}
\frac{\dot r_\eps}{r_\eps} U_\eps(t,r_\eps s_\eps) \dt\\
- C_1 \int_{t_0-\delta}^{t_0}r_\eps^{\alpha'} \left(U_\eps(t,r_\eps s_\eps)+1\right)dt\\
+ \int_{t_0-\delta}^{t_0}\left[-r^{\alpha'+1}_\eps r s \cdot s_\eps  
+ \frac{\alpha'}{2}r^{\alpha'-1}_\eps\dot r_\eps|rs-r_\eps s_\eps| 
+ r_\eps^{\alpha'}(rs-r_\eps s_\eps)\frac{d}{dt}(rs-r_\eps s_\eps)\right]\dt.
\end{multline*}
The right hand side of the last inequality is bounded above 
because of the boundedness of 
$\Gamma_{\alpha',\eps}$.
Passing to the limit as $\eps \to 0$, from Proposition~\ref{propo:conv} and the
boundedness above of the function $\Gamma_{\alpha'}$ it follows that
\begin{multline*}
\Gamma_{\alpha'}(t_0) - \Gamma_{\alpha'}(t_0-\delta) \geq
-\int_{t_0-\delta}^{t_0}(\alpha'-\alpha-\lambda)r^{\alpha'}\frac{\dot r}{r}
U(t,rs)\dt\\
- C_1 \int_{t_0-\delta}^{t_0} r^{\alpha'} \left(U(t,rs)+1\right)\dt
-\int_{t_0-\delta}^{t_0}r^{\alpha'+2}\dt,
\end{multline*}
where $\lambda\leq (\alpha'-\alpha)/2$ is the limit, up to subsequences, of the
bounded sequence $(\lambda_\eps)_\eps$.
Now we recall that, by Lemma~\ref{le:2}, $U(t,rs)$ is integrable; this fact implies the
 integrability of the function
$r^{\alpha'} U(t,rs)-r^{\alpha'+2}$ and hence the existence of a constant $K$ such that
\[
\int_{t_0-\delta}^{t_0} - r^{\alpha'} \frac{\dot r}{r} U(t,rs) dt \leq K <
+\infty.
\]
\end{proof}
\begin{lemma}[{\bf Monotonicity Formula}]
\label{le:Gamma}
Let $\bar x$ be a generalized solution, let $t_0 \in (a,b]$ be a total
collision instant and let $\delta>0$ be the constant 
obtained in Theorem~\ref{thm:totcolliso}. Then
the function 
\begin{equation}
\label{eq:function_Gamma}
\Gamma_\alpha (t) := r^\alpha \left[ \frac12 r^2|\dot s|^2 - U(t,rs) \right]
\end{equation}
is bounded on $[t_0-\delta,t_0)$ and
\begin{equation}\label{eq:ineq:dGamma}
\begin{split}
\Gamma_\alpha (t)  \geq \Gamma_\alpha(t_0-\delta) 
&- \int_{t_0-\delta}^{t}\frac{2-\alpha}{2} r^{1+\alpha} \dot r |\dot s|^2 \,d\xi\\
&- C_1 \int_{t_0-\delta}^{t} r^\alpha \left(U(\xi,rs)+1\right)\,d\xi
 + C_2 \int_{t_0-\delta}^{t} r^{\alpha+\gamma} \frac{\dot r}{r} U(\xi,rs) \,d\xi
\end{split}
\end{equation}
where $t\in[t_0-\delta,t_0]$.
\end{lemma}
\begin{proof}
Replacing in~\eqref{eq:energy_rs} the expression of the function $\Gamma_\alpha$
we have 
\[
\Gamma_\alpha (t) = h(t)r^\alpha - \dot r^2r^\alpha \leq h(t)r^\alpha;
\]
since $h$ is bounded (see Corollary~\ref{cor:2}) and $r$ tends to $0$, we
conclude that the function $\Gamma_\alpha$ is bounded above.
Using the same approximation arguments described in Lemma~\ref{le:integ}, we
obtain~\eqref{eq:ineq:dGamma}.
From Lemma~\ref{le:integ} we deduce the integrability of the negative function 
$r^{\alpha+\gamma} \frac{\dot r}{r} U(t,rs)$.
Hence, since $-\frac{2-\alpha}{2} r^{1+\alpha} \dot r |\dot s|^2$
is positive and both $r^\alpha U(t,rs)$ and $r^{\alpha+\gamma} \frac{\dot r}{r}
U(t,rs)$ are integrable (Lemma~\ref{le:integ}), the boundedness below of
the function $\Gamma_\alpha$ follows from~\eqref{eq:ineq:dGamma}.
\end{proof}

\begin{corollary}
\label{cor:lim2}
In the same setting of Lemma~\ref{le:Gamma} we have 
$\displaystyle \int_{t_0-\delta}^{t_0}-r^{1+\alpha}\dot r |\dot s|^2 < +\infty$.
\end{corollary}
\begin{proof}
It follows from the boundedness above of the function
$\Gamma_\alpha$ 
and inequality~\eqref{eq:ineq:dGamma} since the terms 
$- C_1 \int_{t_0-\delta}^{t} r^\alpha \left(U(t,rs)+1\right) dt$ and
$ C_2 \int_{t_0-\delta}^{t} r^{\alpha+\gamma} \frac{\dot r}{r} U(t,rs)\,dt$ 
are
negative.
\end{proof}

\begin{lemma} \label{le:lim1}
Let $\varphi(t):= -\dot r(t) r^{\alpha/2}(t)$, $t \in [t_0-\delta,t_0]$. 
Then there exist two constants depending on $\alpha$, $c_{1,\alpha} \leq
c_{2,\alpha}$, such that
for all  $t \in [t_0-\delta,t_0]$
\[
c_{1,\alpha} \leq \varphi(t) \leq c_{2,\alpha}. 
\]
\end{lemma}
\begin{proof}
Since the energy function $h$ is bounded (see Corollary~\ref{cor:2}, {\em(iv)}) 
and we assume that $r$ tends to $0$ as $t$ tends to $t_0$,
the function $r^\alpha h(t)=\dfrac{1}{2}\varphi^2(t) + \Gamma_\alpha(t)$ 
is also bounded and by Lemma~\ref{le:Gamma} we can deduce that 
$\varphi(t)=\sqrt{2\left[r^\alpha h(t) - \Gamma_\alpha(t) \right]}$ is bounded
below and above by a pair of constants $0\leq c_{1,\alpha} \leq c_{2,\alpha}$ on
the interval $[t_0-\delta,t_0]$.
\end{proof}

\begin{corollary}
\label{cor:int2}
In the same setting of Lemma~\ref{le:Gamma} we have 
$\displaystyle \lim_{t \to t_0}\int_{t_0-\delta}^{t} \frac{1}{r^{\alpha/2
+1}} = +\infty$.
\end{corollary}
\begin{proof}
We can write the boundedness above of the function $\varphi$ (proved in
Lemma~\ref{le:lim1}) as 
\begin{equation}
\label{eq:in}
-\frac{\dot r}{r} \leq \frac{c_{2,\alpha}}{r^{\alpha/2 +1}}, \quad t \in
[t_0-\delta,t_0].
\end{equation}
Integrating inequality~\eqref{eq:in} on the interval $[t_0-\delta, t]$, when $t
\to t_0$, we obtain 
\[
\lim_{t \to t_0} c_{2,\alpha}\int_{t_0-\delta}^{t}\frac{d\xi}{r^{\alpha/2 +1}} 
\geq \lim_{t \to t_0} \int_{t_0-\delta}^{t} -\frac{\dot r}{r}d\xi
= \log r(t_0-\delta)-\lim_{t \to t_0} \log r(t) = +\infty
\]
since $r$ tends to 0 as $t \to t_0$.
\end{proof}

\begin{lemma}
\label{le:lim2}
The lower bound $c_{1,\alpha}$ 
of the function $\varphi$ defined in Lemma~\ref{le:lim1} can be chosen
strictly positive, that is $c_{1,\alpha} > 0$.
\end{lemma}
\begin{proof}
We start proving an estimate above of the derivative of the function
$\varphi$. With this purpose we consider the approximating sequence
$(\varphi_\eps)_\eps$ where 
\[
\varphi_\eps (t)=-\dot r_\eps(t) r_\eps^{\alpha/2}(t)
\]
and, for every $\eps>0$, we compute the first derivative of the smooth function
$\varphi_\eps$ and we use the Euler--Lagrange equation~\eqref{eq:ELeps}$_1$ 
for the approximating problem to obtain
\begin{equation*}
\begin{split}
\dot\varphi_\eps(t)
& =-\frac{\alpha}{2}r_\eps^{\alpha/2-1}\dot r_\eps^2-r_\eps^{\alpha/2}\ddot
r_\eps \\
& = -\frac{\alpha}{2} r_\eps^{\alpha/2 -1}\dot r_\eps^2 -
r_\eps^{\alpha/2+1}|\dot s_\eps|^2   -r_\eps^{\alpha/2-1} \nabla U_\eps(t,r_\eps
s_\eps)\cdot (r_\eps s_\eps)+r_\eps^{\alpha/2}(rs-r_\eps s_\eps)\cdot s_\eps.
\end{split}
\end{equation*}
Arguing as in the proof of Lemma~\ref{le:integ} we use assumptions \ref{a:U2h}
and~\eqref{etaeps_ineq} to deduce
\begin{equation*}
\begin{split}
\dot\varphi_\eps(t)
& \leq r_\eps^{\alpha/2-1}\left[-\frac{\alpha}{2}\dot r_\eps^2-r_\eps^2|\dot
s_\eps|^2  +(\alpha + C_2r_\eps^\gamma) U_\eps(t,r_\eps s_\eps)+(rs-r_\eps
s_\eps)\cdot (r_\eps s_\eps) \right] \\
& = \frac{1}{r_\eps^{\alpha/2 +1}}\Big[\frac{2-\alpha}{2}\varphi_\eps^2(t) -
2r_\eps^\alpha h_\eps(t)-(2-\alpha)r_\eps^{\alpha}U_\eps(t,r_\eps s_\eps) +
C_2r_\eps^{\alpha+\gamma} U_\eps(t,r_\eps s_\eps) \\
& \hspace{3cm} + r_\eps^\alpha(rs-r_\eps s_\eps)\cdot (rs) \Big]
\end{split}
\end{equation*}
and then, for every $ \in (t_0-\delta,t_0)$,
\begin{equation*}
\begin{split}
\varphi_\eps(t)
& \leq \varphi_\eps^0 +\int_{t_0-\delta}^{t}
\frac{1}{r_\eps^{\alpha/2 +1}}\Big[\frac{2-\alpha}{2}\varphi_\eps^2(\xi) -
2r_\eps^\alpha h_\eps(\xi)-(2-\alpha)r_\eps^{\alpha}U_\eps(\xi,r_\eps s_\eps)\\
& \hspace{3cm} + C_2r_\eps^{\alpha+\gamma} U_\eps(\xi,r_\eps s_\eps) 
               + r_\eps^\alpha(rs-r_\eps s_\eps)\cdot (rs) \Big]d\xi
\end{split}
\end{equation*}
where $\varphi_\eps^0=\varphi_\eps(t_0-\delta)$.
As $\eps\to 0$, from Proposition~\ref{propo:conv} we have 
\begin{multline*}
\varphi(t) \leq \varphi^0 \\+\int_{t_0-\delta}^{t}
\frac{1}{r^{\alpha/2 +1}}\left[\frac{2-\alpha}{2}\varphi^2(\xi) - 2r^\alpha
h_\eps(\xi)-(2-\alpha)r^{\alpha}U(\xi,rs) + C_2 r^{\alpha+\gamma} U(\xi,rs)
\right]d\xi,
\end{multline*}
where $\varphi^0=\varphi(t_0-\delta)$.
Since $\gamma>0$ there exists $\lambda \in (0,2-\alpha)$, such that
\begin{equation*}
\varphi(t) \leq \varphi^0 +\int_{t_0-\delta}^{t}
\frac{1}{r^{\alpha/2 +1}}\left[\frac{2-\alpha}{2}\varphi^2(\xi) + C(\xi)
-(2-\alpha-\lambda)r^{\alpha}U(\xi,rs) \right]d\xi,
\end{equation*}
where $C(t)$ is such that $|C(t)| \to 0$ as $t \to t_0$ and $2r^\alpha h(t) \leq
C(t)$ on $[t_0-\delta,t_0]$. 
Furthermore, the uniform convergence assumed in condition \ref{a:U3h} implies that,
denoting by
$\tilde U_0$ the minimal value assumed by $\tilde U$ on the ellipsoid
${\mathcal E}$, there exist two positive constants $k_1,k_2 >0$ such that
\begin{equation*}
\varphi(t) \leq \varphi^0 +\int_{t_0-\delta}^{t}\frac{k_1}{r^{\alpha/2
+1}}\left( \varphi^2(\xi) - k_2 \tilde U_0 \right)d\xi
\end{equation*}
whenever $\delta$ is sufficiently small.
We will conclude showing that necessarily $\varphi^2(t) \geq k_2 \tilde U_0$ and
then choosing $c_{1,\alpha}:= \sqrt{k_2 \tilde U_0} > 0$.

By the sake of contradiction we suppose the existence of $\hat t$ such that 
$\varphi^2(\hat t) < k_2 \tilde U_0$; 
then $\varphi^2 - k_2 \tilde U_0 < 0$ in a neighborhood of $\hat t$ and 
\begin{equation*}
\varphi(t) \leq \varphi(\hat t) +\int_{\hat t}^{t}\frac{k_1}{r^{\alpha/2
+1}}\left( \varphi^2(\xi) - k_2 \tilde U_0 \right)d\xi < \varphi(\hat t)
\end{equation*}
for every $t \in (\hat t,t_0)$.
We deduce the existence of a strictly positive constant $\hat k$ such that, 
for every $t \in (\hat t,t_0)$,
\begin{equation*}
\varphi(t)-\varphi(\hat t) \leq - \hat k \int_{\hat
t}^{t}\frac{d\xi}{r^{\alpha/2 +1}}
\end{equation*}
Since the right hand side tends to $-\infty$ as $t$ approaches $t_0$ (see
Corollary~\ref{cor:int2}), the last inequality contradicts the boundedness of
the function $\varphi$.
\end{proof}

\begin{corollary}
\label{cor:stime_r}
There exist two strictly positive constants $0<k_{1,\alpha} \leq k_{2,\alpha}$
such that 
\[
k_{1,\alpha}(t_0-t)^{\frac{2}{\alpha+2}} \leq r(t) \leq
k_{2,\alpha}(t_0-t)^{\frac{2}{\alpha+2}},
\]
whenever $t \in [t_0-\delta_0,t_0]$.
\end{corollary}
\begin{proof}
The statement follows from Lemmata~\ref{le:lim1} and~\ref{le:lim2} with 
$k_{i,\alpha} := \left(\frac{\alpha+2}{2}c_{i,\alpha}
\right)^{\frac{2}{\alpha+2}}$, $i=1,2.$
\end{proof}

\begin{corollary}
\label{cor:limGamma}
There exists $b>0$ such that 
\[
\lim_{t \to t_0^-} \Gamma_\alpha(t) = -b \qquad \mbox{and}\qquad
\lim_{t \to t_0^-} \dot r^2 r^\alpha = 2b.
\]
\end{corollary}
\begin{proof}
Since $\Gamma_\alpha$ is bounded and inequality~\eqref{eq:ineq:dGamma} holds,
$\Gamma_\alpha$ admits a limit when $t$ tends to $t_0$ from the right.
We call this limit $-b \in \RR$.
Since $\Gamma_\alpha (r,s) = h(t)r^\alpha - \dfrac{1}{2}\dot r^2 r^\alpha$, the
energy $h$ is bounded  and $r$ tends to $0$ as $t \to t_0$ we conclude that $\dot r^2 r^\alpha$ 
converges to $2b$ and, using Lemma~\ref{le:lim2} we
deduce that $b>0$.
\end{proof}

\begin{theorem}
\label{central_conf}
Let $\bar x$ be a generalized solution for the dynamical system~\eqref{DS},
let $t_0 \in (a,b)$ (if $b<+\infty$ $t_0$ can coincide with $b$) be a total
collision instant and let $\delta>0$ be the constant
obtained in Theorem~\ref{thm:totcolliso}. 
Let $r,s$ be the new variables defined in~(\ref{rs}); if the potential $U$
satisfies assumptions~\ref{a:U0}, \ref{a:U1}, \ref{a:U2h}, \ref{a:U3h} then the following
assertions hold
\begin{itemize}
\item[(a)] $\dlim_{t \rightarrow {t_0^-}} r^\alpha U(t,rs)=b$, where $b$ is
the strictly positive constant introduced in Corollary~\ref{cor:limGamma};
\item[(b)] there is a positive constant $K$ 
such that, as $t$ tends to $t_0$, 
\[
\begin{aligned}
r(t) & \sim [K
({t_0}-t)]^{\frac{2}{2+\alpha}}  \\ 
\dot r(t) & \sim -\frac{2K}{2+\alpha}
[K({t_0}-t)] ^{\frac{-\alpha}{2+\alpha}};
\end{aligned}
\]
\item[(c)] $\ds \lim_{t \rightarrow t_0^-} |\dot s(t)|({t_0}-t)=0$;
\item[(d)] for every real positive sequence $(\lambda_n)_n$, such that
           $\lambda_n \rightarrow 0$ as $n \rightarrow +\infty$ we have   
           \[
           \lim_{n \rightarrow +\infty} |s({t_0}-\lambda_n) - s({t_0}-\lambda_n
t)|=0, \quad \forall t > 0.
           \]
\end{itemize}
\end{theorem}

\begin{remark}
Condition 
{\it (a)} of Theorem~\ref{central_conf} together with assumptions \ref{a:U3h} on $U$ 
and~\eqref{tildeU} on $\tilde U$ imply that, if $\bar x$ is generalized solution
and $|\bar x(t_0)|=0$, then there exist $\delta>0$ such that, for every $t\in(t_0-\delta,t_0)$,
$\bar x(t)\notin\Delta$, i.e., in a (left) neighborhood 
of the total collision instant no other collision is allowed: neither total nor partial. 
As a consequence, in such a neighborhood, the generalized solution $\bar x$ satisfies the
dynamical system~\eqref{DS} and the corresponding variables $(r,s)$ verify the Euler--Lagrange
equations~\eqref{eq:EL}.
\end{remark}

\begin{proof}[Proof of Theorem \ref{central_conf}]
We begin by 
proving statement \emph{(a)}.
The boundedness of the function $\Gamma_\alpha$ together with
inequality~\eqref{eq:ineq:dGamma} imply the integrability of the function
$r^{\alpha+1}\dot r |\dot s|^2$ on the interval $[t_0-\delta,t_0]$ (see
Corollary~\ref{cor:lim2}). Furthermore, since the integral of $\dot r/r$ on the
same interval diverges to $-\infty$ we conclude that
\[
\liminf_{t \to t_0^-} r^{\alpha+2}|\dot s|^2 = 0
\]
and from~\eqref{eq:function_Gamma} together with Corollary~\ref{cor:limGamma}
\begin{equation}
\label{eq:liminfU}
\liminf_{t \to t_0^-} r^{\alpha} U(t,rs) = b.
\end{equation}
It remains to prove that also $\ds \limsup_{t \to t_0^-} r^{\alpha} U(t,rs) =
b$.
Suppose, 
for the sake of contradiction, the existence of a strictly positive
$\eps$
such that 
\begin{equation}
\label{eq:limsupU}
\limsup_{t \to t_0^-} r^{\alpha} U(t,rs) = b + 3\eps.
\end{equation}
Using assumption \ref{a:U3h} we have that~\eqref{eq:liminfU} and~\eqref{eq:limsupU}
are respectively equivalent to 
\[
\liminf_{t \to t_0^-} \tilde U(t,s) = b \quad \mbox{ and } \quad 
\limsup_{t \to t_0^-} \tilde U(t,s) = b + 3\eps
\]
and Corollary~\ref{cor:limGamma} implies the existence of $t_{\eps}$ such that
$\Gamma_\alpha(t) \geq -b-\eps/2$ whenever $t \in (t_{\eps} , t_0]$. We can then
define the set 
\[
{\cal U} := \left \{ t \in (t_{\eps} , t_0) : \tilde U(t,s(t))\geq b+\eps
\right\}.
\]
We define two non-empty subsets of the ellipsoid ${\mathcal E}$ as 
\[
A:= \left\{ s(t) : \tilde U(t,s(t)) \leq b+\eps \right\} \quad \mbox{ and }
\quad
B:= \left\{ s(t) : \tilde U(t,s(t)) \geq b+2\eps \right\};
\]
since $\eps > 0$ the quantity 
\[
d:=\mathrm{dist}(A,B)=\inf_{s_1 \in A, s_2 \in B} |s_1-s_2|
\]
is strictly positive and there exists a sequence $(t_n)_{n\geq 0} \subset
[t_0-\delta,t_0]$, such that
\[
\begin{split}
& t_n \to t_0 \mbox{ as } n \to +\infty  \\
& s(t_{2k}) \in \partial A \quad \mbox{and} \quad s(t_{2k+1}) \in \partial B
\quad \mbox{ for every } k \in \NN \\
& b+\eps \leq \tilde U(t,s(t)) \leq b+2\eps, \quad \mbox{for every } t \in
(t_{2k},t_{2k+1}) \mbox{ and } k \in \NN.
\end{split}
\]
Hence $(t_{2k},t_{2k+1}) \subset {\cal U}$, for every $k$, and
from the definition of the function $\Gamma_\alpha$ in~\eqref{eq:function_Gamma}
we have that
\begin{equation}
\label{eq:ass}
r^{\alpha+2}|\dot s|^2 \geq \eps \mbox{ in the intervals } (t_{2k},t_{2k+1}). 
\end{equation}
We now estimate the integral on $(t_{2k},t_{2k+1})$ of the integrable (on
$[t_0-\delta,t_0]$) function $r^{\alpha+1}\dot r|\dot s|^2$ using~\eqref{eq:ass}
and Corollary~\ref{cor:stime_r}
\begin{equation}
\label{eq:est1}
\begin{split}
\int_{t_{2k}}^{t_{2k+1}} -\frac{\dot r}{r}r^{\alpha+2}|\dot s|^2 dt 
& \geq \eps \int_{t_{2k}}^{t_{2k+1}} - \frac{\dot r}{r} dt 
= \eps\log \frac{r(t_{2k})}{r(t_{2k+1})} \\
&\geq \frac{2\eps}{2+\alpha}\log
\frac{c_{1,\alpha}(t_0-t_{2k})}{c_{2,\alpha}(t_0-t_{2k+1})}.
\end{split}
\end{equation}
On the other hand, using H\"older inequality, we have 
\begin{equation}
\label{dis:1}
d^2 \leq |s(t_{2k+1})-s(t_{2k})| \leq \left( \int_{t_{2k}}^{t_{2k+1}} |\dot s| dt \right)^2 \leq 
         \int_{t_{2k}}^{t_{2k+1}} -r^{\alpha+2}\frac{\dot r}{r}|\dot s|^2 dt 
         \int_{t_{2k}}^{t_{2k+1}} \frac{dt}{-r^{\alpha+1}\dot r} 
\end{equation}
and from Lemma~\ref{le:lim1} and Corollary~\ref{cor:stime_r}, we obtain
\begin{equation}
\label{dis:2}
\begin{split}
\int_{t_{2k}}^{t_{2k+1}} \frac{dt}{-r^{\alpha+1}\dot r} & =
\int_{t_{2k}}^{t_{2k+1}} \frac{1}{-r^{\alpha/2}\dot r}\frac{1}{r^{\alpha/2+1}}dt
\\ 
& \leq \frac{2}{2+\alpha}\frac{1}{c^2_{1,\alpha}} \int_{t_{2k}}^{t_{2k+1}}
\frac{dt}{t_0-t} =
\frac{2}{2+\alpha}\frac{1}{c^2_{1,\alpha}} \log \frac{t_0-t_{2k}}{t_0-t_{2k+1}}.
\end{split}
\end{equation}
Combining~\eqref{dis:1} and~\eqref{dis:2} we obtain
\begin{equation}
\label{eq:est2}
\int_{t_{2k}}^{t_{2k+1}} -r^{\alpha+2}\frac{\dot r}{r}|\dot s|^2 dt \geq 
\frac{2+\alpha}{2}d^2 c_{1,\alpha}^2 \left[ \log
\frac{t_0-t_{2k}}{t_0-t_{2k+1}}\right]^{-1}.
\end{equation}
From the estimates~\eqref{eq:est1} and~\eqref{eq:est2} we deduce
\[
\int_{t_{2k}}^{t_{2k+1}} - r^{\alpha+2}\frac{\dot r}{r}|\dot s|^2 dt \geq 
\frac{\eps}{2+\alpha}\log
\frac{c_{1,\alpha}(t_0-t_{2k})}{c_{2,\alpha}(t_0-t_{2k+1})} +
\frac{2+\alpha}{4}d^2 c_{1,\alpha}^2 \left[ \log
\frac{t_0-t_{2k}}{t_0-t_{2k+1}}\right]^{-1}.
\]
Summing on the index $k$ and recalling that the positive function $-\dot r
r^{\alpha+1}|\dot s|^2$ has a finite integral on $[t_0-\delta,t_0]$ (Corollary
\ref{cor:lim2}) we have
\begin{equation}
\label{mult:rhs}
\begin{split}
+\infty & > \int_{t_0-\delta}^{t_0} -\dot r r^{\alpha+1}|\dot s|^2 dt 
> \sum_{k \geq 0} \int_{t_{2k}}^{t_{2k+1}} -\dot r r^{\alpha+1}|\dot s|^2 dt \\
& \geq \frac{\eps}{2+\alpha}\sum_{k \geq 0}\log
\frac{c_{1,\alpha}(t_0-t_{2k})}{c_{2,\alpha}(t_0-t_{2k+1})} +
\frac{2+\alpha}{4}d^2 c_{1,\alpha}^2 \sum_{k \geq 0}\left[ \log
\frac{t_0-t_{2k}}{t_0-t_{2k+1}}\right]^{-1}.
\end{split}
\end{equation}
Since $c_{2,\alpha}/c_{1,\alpha}$ is 
bounded  (see Lemma~\ref{le:lim2}),
for the last term in~\eqref{mult:rhs} to be  finite it 
is necessary that 
\begin{equation}
\label{eq:limits}
\lim_{k \to+\infty}\frac{t_0-t_{2k}}{t_0-t_{2k+1}}=\frac{c_{2,\alpha}}{c_{1,\alpha}}
\quad \mbox{and} \quad
\lim_{k \to +\infty}\frac{t_0-t_{2k}}{t_0-t_{2k+1}}=+\infty.
\end{equation}
This is a contradiction, hence we conclude that 
\[
\limsup_{t \to t_0} r^\alpha U(t,rs)=b
\]
and, after replacing the value in~\eqref{eq:function_Gamma},
\begin{equation}
\label{eq:limspunto}
\lim_{t \to t_0} r^{\alpha+2}|\dot s|^2 =0.
\end{equation}

\noindent To prove \emph{(b)},
from Corollary~\ref{cor:limGamma} we obtain
\[
\lim_{t \to t_0^-}\frac{r(t)^{\alpha/2 + 1}}{(\alpha/2 + 1)(t_0-t)} =
\lim_{t \to t_0^-}-r(t)^{\alpha/2}\dot r(t) = \sqrt{2b};
\]
we then conclude by 
defining $\displaystyle K:=\frac{2+\alpha}{2}\sqrt{2b}$.
The second estimate follows directly.
% of Theorem \ref{central_conf}
% directly follows from \emph{(b)}.

\noindent Part \emph{(c)}
directly follows from~\eqref{eq:limspunto} and \emph{(b)}.

\noindent We conclude by
proving  statement \emph{(d)}.
If $t=1$ there is nothing to prove. Suppose $t>0$, $t\neq 1$, and
consider a sequence $(\lambda_n)_n$, $\lambda_n \rightarrow 0$;
let $N$ be such that $\lambda_n < \delta/\max(1,t)$, $\forall n \geq N.$
Whenever $t>1$, for every $n\geq N$, we have
\[
t_0-\delta < t_0 -\lambda_nt < t_0 -\lambda_n < t_0
\]
and
\[
\begin{split}
|s({t_0}-\lambda_n)-s({t_0}-\lambda_n t)| 
&\leq \int_{{t_0}-\lambda_n t}^{{t_0}-\lambda_n}|\dot s|du \\
&\leq \left( \int_{{t_0}-\lambda_n t}^{{t_0}-\lambda_n}r^{1+\alpha/2}| 
\dot s|^2du\right)^{1/2} 
\left(\int_{{t_0}-\lambda_n
t}^{{t_0}-\lambda_n}\frac{du}{r^{1+\alpha/2}}\right)^{1/2}
\end{split}
\]
It is not restrictive to suppose $t>1$: indeed, when $t \in(0,1)$, we obtain an equivalent 
estimate by permuting the integration bounds.
From Corollary~\ref{cor:lim2} and Lemmata~\ref{le:lim1} and~\ref{le:lim2} we obtain
\[
+ \infty > \int_{t_0-\delta}^{t_0} r^{1+\alpha} \dot r |\dot s|^2 du 
      \geq \int_{t_0-\delta}^{t_0} c_{1,\alpha}r^{1+\alpha/2}|\dot s|^2 du.
\]
Then, since the constant $c_{1,\alpha}$ is strictly positive, we have
\[
\lim_{n \rightarrow +\infty} \int_{{t_0}-\lambda_n
t}^{{t_0}-\lambda_n}r^{1+\alpha/2}|\dot s|^2 du = 0.
\]
Moreover, as $n$ tends to $+\infty$, the second integral
$\int_{{t_0}-\lambda_n t}^{{t_0}-\lambda_n}r^{-(1+\alpha/2)} < +\infty$; 
indeed both integration bounds tend to $t_0$ and the asymptotic estimate proved
in 
\emph{(b)} holds. 
Hence, as $\lambda_n \rightarrow 0$
\[
\begin{split}
\lim_{n\to+\infty}\int_{{t_0}-\lambda_n
t}^{{t_0}-\lambda_n}\frac{du}{r^{1+\alpha/2}}  & =
\lim_{n\to+\infty}\left[\int_{{t_0}-\lambda_n
t}^{{t_0}-\lambda_n}C\frac{du}{(t_0-u)}  +o(1)\right]\\
& = C \lim_{n\to+\infty} [\log(\lambda_n)-\log(\lambda_n t)+o(1)] = -C\log t
\end{split}
\]
that is bounded since $t$ is fixed and 
$C=\left[\dfrac{\sqrt{2b}(\alpha+2)}{2}\right]^{-(\alpha+2)/2}$.
\end{proof}

\begin{theorem}\label{central_conf2}
In the same setting of Theorem~\ref{central_conf},
assume that the potential 
$U$ verifies the further assumption
\begin{assiomi}
\item[\uqhtag]
\label{a:U4h}
$\ds \lim_{r \to 0} r^{\alpha+1}\nabla_{T} U(t,x)=\nabla_{T}\tilde U(t,s)$.
\end{assiomi}
Then 
\[
\dlim_{t \rightarrow {t_0}} \mathrm{dist} \left({\cal C}^b,s(t)\right)=
\dlim_{t \rightarrow {t_0}} \inf_{\bar s \in {\cal C}^b}|s(t) - \bar s | = 0,
\]
where ${\cal C}^b$ is the set of central configurations for $\tilde U$ at level $b$, namely
the subset of critical points of the restriction of 
$\tilde U$  to the ellipsoid ${\mathcal E}$:
\begin{equation}
\label{eq:cal_C}
{\cal C}^b := \left\{ s : \tilde U(t_0,s) = b, \nabla_{T} \tilde U(t_0,s) = 0 \right\}.
\end{equation}
\end{theorem}

\begin{remark}
When $U$ is homogeneous, as in the classical keplerian potential, then 
$\tilde U$ is simply the restriction of $U$ 
 on ${\mathcal E}$ and Theorem~\ref{central_conf2} asserts
that the angular component $s$ of the motion tends to a set of central configurations.
\end{remark}

\begin{proof}
Since in {\it (a)} of Theorem~\ref{central_conf} we have already proved that 
$\lim_{t \to t_0} \tilde U(t,s(t))=b$, it remains to show that 
\[
\lim_{t \to t_0^-} |\nabla_{T} \tilde U(t,s(t))| = 0
\]
that, using condition \ref{a:U4h}, is equivalent to 
\[
\lim_{t \to t_0^-} r^{\alpha + 1}|\nabla_{T} U(t,rs)| = 0.
\]
We now consider the Euler--Lagrange equation~\eqref{eq:EL}$_2$
multiplied by $r^{\alpha}$
\[
-2 r^{\alpha + 1}\dot r\dot s - r^{\alpha + 2}\ddot s + r^{\alpha + 1}\nabla_{T}
U(t,rs) = r^{\alpha + 2}|\dot s|^2 s;
\]
since $r^{\alpha + 1}\dot r\dot s=r^{\alpha/2 + 1}\dot r r^{\alpha/2}\dot s$ is the product
of a bounded term with an infinitesimal one (see equation~\eqref{eq:limspunto} and
Lemma~\ref{le:lim1}), while $|r^{\alpha + 2}|\dot s|^2 s| = r^{\alpha + 2}|\dot s|^2$
tends to 0 for~\eqref{eq:limspunto}, we claim that
\begin{equation}\label{eq:ddot_s}
\lim_{t \to t_0^-} r^{\alpha + 2}\ddot s = 0.
\end{equation}
We perform the time rescaling (cf.~McGehee's change of coordinates in~\ref{eq:mcgehee2}) 
\begin{equation}\label{eq:mcgehee1}
\tau = \int_{t_0-\delta}^{t} \frac{d\xi}{r^{\alpha/2+1}}
\end{equation}
which maps the interval $[t_0-\delta,t_0)$ into $[0,+\infty)$ (see Corollary~\ref{cor:int2}).
If 
the prime~$\ {}'\ $ denotes
the derivative with respect to the new variable $\tau$, 
then~\eqref{eq:ddot_s} is equivalent to
\begin{equation}\label{eq:ddot_s2}
\lim_{\tau \to +\infty} s''(\tau) = 0
\end{equation}
and the limit~\eqref{eq:limspunto} reads simply 
\begin{equation}\label{eq:ddot_s3}
\lim_{\tau \to +\infty} |s'(\tau)|^2 = 0.
\end{equation}
Suppose now,
for the sake of contradiction,
that there exists a sequence $(\tau_n)_n$ 
such that $\tau_n \to +\infty$ as $n \to +\infty$ and 
\[
\lim_{n\to+\infty} \nabla_{T} \tilde U(\tau_n,s(\tau_n)) = \lim_{n\to+\infty}  s''(\tau_n)= 
\sigma
\]
for some $\sigma \neq 0$.
Since the ellipsoid ${\mathcal E}$ is compact, up to subsequences,
$\left(s(\tau_n)\right)_n$ converges to some $\bar s$. 
Furthermore, from Theorem~\ref{central_conf} we know that
$\tilde U\left(\tau_n,s(\tau_n)\right)$ tends to the finite limit $b$ as 
$n \to +\infty$, hence
$(t_0,\bar s)$ is a regular point both for $\tilde U$ and for $\nabla_T\tilde U$.
We moreover remark that, since the limit~\eqref{eq:ddot_s3} holds, for every fixed positive constant
$h>0,$ there holds
\[
s(\tau) \to \bar s, \qquad \mbox{uniformly on } [\tau_n,\tau_n+h], \; \mbox{for every }n
\]
and also
\[
\sup_{\tau \in [\tau_n,\tau_n+h]}|\nabla_T \tilde U(\tau,s(\tau))-\sigma| \to  0, 
\qquad \mbox{as }n \to +\infty.
\]
We can then compute
\[
\begin{split}
s'(\tau_n +h)-s'(\tau_n) &= \int_{\tau_n}^{\tau_n+h} s''(\tau) d\tau \\
& = \int_{\tau_n}^{\tau_n+h} \nabla_T \tilde U(\tau,s(\tau)) d\tau + o(1) \\
& = h \sigma + o(1) \qquad \mbox{as } n\to +\infty.
\end{split}
\]
We obtain the contradiction
\[
0 = \lim_{n \to +\infty}|s'(\tau_n +h)-s'(\tau_n)| = h |\sigma| \neq 0.
\]
\end{proof}

%--------------------------------
\subsection{Logarithmic potentials}
%--------------------------------
Aim of this section is  to extend  the asymptotic estimates of 
Theorem~\ref{central_conf}
to potentials having
logarithmic singularities.
We follow the same scheme and we still work in a left neighborhood of a total 
collision instant $t_0$, $(t_0-\delta,t_0)$. The main differences concern the
 monotonicity formul\ae\/ (Lemmata~\ref{le:integ_log} and~\ref{le:Gamma_log}).

In this setting, we suppose the existence of a continuous function 
\begin{equation} \label{eq:M}
M\from 
(a,b)\to\RR \quad \mbox{ such that $\dot M(t)$ is bounded on $(t_0-\delta,t_0)$}
\end{equation}
and we replace conditions \ref{a:U2h} and \ref{a:U3h} with
\begin{assiomi}
\item[\udltag]
\label{a:U2l}
There exist $\gamma>0$ and $C_2 \geq 0$ such that
\[
\nabla U(t,x)\cdot x + M(t) \geq - C_2 |x|^\gamma U(t,x),
\]
whenever $|x|$ is small.
\item[\utltag]
\label{a:U3l}
$\ds \lim_{|x| \to 0} \left[U(t,x) + M(t)\log |x| \right] = \tilde U(t,s)$, 
uniformly in $t$,
\end{assiomi}
where $\tilde U$, as in the quasi--homogeneous case, is of class $\cont^1$
on $(a,b)\times ({\mathcal E} \minus \Delta)$ and verifies~\eqref{tildeU}.

\begin{remark}
\ref{a:U2l} implies \ref{a:U2} (for small value of $|x|$) for every $\tilde \alpha \in (0,2)$.
\end{remark}

\begin{remark}\label{rem:int_log}
From Corollary~\ref{cor:2} and assumption~\ref{a:U3l} it follows
that
the positive function $-M(t)\log |x|+\tilde U(t,s)$ is integrable
in a neighborhood of a total collision at the origin.
\end{remark}

We now prove the analogue of Lemmata~\ref{le:integ} and~\ref{le:Gamma} in the setting of
logarithmic--type potentials.
\begin{lemma}
\label{le:integ_log}
Let $\bar x$ be a generalized solution, let $t_0 \in (a,b]$ be a total
collision instant and let $\delta$ be given in Theorem~\ref{thm:totcolliso}.
Let $\gamma$ be the positive exponent appearing in \ref{a:U2h}, then 
\begin{equation}
\label{eq:integ_log}
\int_{t_0-\delta}^{t_0}-r^{\gamma}\frac{\dot r}{r}U(t,rs)dt < +\infty.
\end{equation}
\end{lemma}
\begin{proof}
We define the functions
\begin{equation}
\label{eq:function_Gamma_log}
\Gamma_{log} (r,s) := \dfrac{1}{2}r^2|\dot s|^2 - \left[ U(t,rs) + M(t)\log r
\right]
\end{equation}
and
\begin{equation}
\label{eq:function_TGamma}
\tilde \Gamma_{log} (r,s) := r^{\gamma}\Gamma_{log};
\end{equation}
since 
\[
\tilde \Gamma_{log} (r,s) = r^{\gamma}\left[ h(t)-\frac12 \dot r^2 - M(t)\log r
\right]
\leq r^{\gamma}h(t) - r^{\gamma}M(t)\log r,
\]
then $\tilde \Gamma_{log}$ is bounded above, indeed $h$ is bounded, $M$
continuous and, since $\gamma>0$, $\lim_{r \to 0} r^{\gamma}\log r=0$. 
We now proceed exactly as in the proof of Lemma~\ref{le:integ}: we omit here
the approximation argument and
we formally compute the time derivative of $\tilde \Gamma_{log}$
\[
\frac{d}{dt}\tilde \Gamma_{log} (r,s) = \gamma r^{\gamma-1}\dot r \Gamma_{log}
(r,s) + r^\gamma \frac{d}{dt} \Gamma_{log} (r,s).
\]
Using the Euler--Lagrange equation~\eqref{eq:EL}$_2$, we obtain
\[
\frac{d}{dt}\Gamma_{log} (r,s) = 
-r \dot r |\dot s|^2 -\frac{\partial U}{\partial t}(t,rs)
- \frac{\dot r}{r}\nabla U(t,rs)\cdot(rs)
- \dot M(t)\log r - M(t)\frac{\dot r}{r}.
\]
From assumptions \ref{a:U1} and \ref{a:U2l} we deduce that 
\begin{equation}
\label{eq:stima_dGammalog}
\frac{d}{dt}\Gamma_{log} (r,s) 
\geq - r \dot r |\dot s|^2 - C_1 \left(U(t,rs)+1\right)
+ C_2 r^{\gamma}\frac{\dot r}{r}U(t,rs) - \dot M(t)\log r
\end{equation}
and then 
\begin{multline}
\label{eq:stima_dTGamma}
\frac{d}{dt}\tilde \Gamma_{log} (r,s) \geq 
- \frac{2-\gamma}{2}r^{\gamma+1} \dot r |\dot s|^2 
- \gamma r^\gamma \frac{\dot r}{r}U(t,rs) 
- \gamma r^\gamma \frac{\dot r}{r}M(t)\log r\\
- C_1 r^\gamma U(t,rs) - C_1 r^\gamma + C_2 r^{2\gamma}\frac{\dot r}{r} U(t,rs)
- \dot M(t)r^\gamma \log r.
\end{multline}
The first term in~\eqref{eq:stima_dTGamma} is positive, since~\eqref{eq:isol_coll}
holds; moreover, since $r$ tends to $0$ as
$t$ approaches $t_0$, there exist $\eps \in (0,\gamma)$ and 
$\delta_0\in (0,\delta]$ such that
\begin{equation}
\label{eq:ultimaGT}
-\gamma r^\gamma \frac{\dot r}{r}U(t,rs)+ C_2 r^{2\gamma}\frac{\dot r}{r}
U(t,rs) \geq -(\gamma-\eps) r^\gamma \frac{\dot r}{r}U(t,rs) \geq 0
\end{equation}
on $(t_0-\delta_0,t_0)$.
The remaining terms in~\eqref{eq:stima_dTGamma} are integrable 
functions, indeed the last term $\dot M(t)r^\gamma \log r$ is bounded as $r$ tends to 0
(see~\eqref{eq:M}),
$r^{\gamma}U \leq U$ and $U$ is integrable and we have the following estimate
\[
-\gamma r^\gamma \frac{\dot r}{r}M(t)\log r \geq 
-\gamma r^{\gamma-1} \dot r\log r \max_{t \in [t_0-\delta,t_0]} M(t)
\]
and
\[
\int_{t_0-\delta}^{t_0}\gamma r^{\gamma-1}\dot r\log r dt =  -r_0^\gamma\log r_0 +
\int_{t_0-\delta}^{t_0} r^{\gamma-1}\dot r dt
= r_0^\gamma \left(-\log r_0+\frac{1}{\gamma}\right) < +\infty
\]
where $r_0=r(t_0-\delta)$.
Hence the right hand side of~\eqref{eq:stima_dTGamma} is the sum of an
integrable function with a
positive one; since the $\tilde \Gamma_{log} (r,s)$ is bounded above
from~\eqref{eq:ultimaGT}
we have the estimate in~\eqref{eq:integ_log}.
\end{proof}

\begin{lemma}[{\bf Monotonicity Formula}]
\label{le:Gamma_log}
The function $\Gamma_{log}$ defined in~\eqref{eq:function_Gamma_log} is bounded
on $[t_0-\delta,t_0]$.
\end{lemma}
\begin{proof}
We consider the expression of the derivative of $\Gamma_{log}$ with respect to
the time variable computed in~\eqref{eq:stima_dGammalog}.
Using Lemma~\ref{le:integ_log}, the integrability of the function $U$ and
Remark~\ref{rem:int_log} we deduce the boundedness below 
(in a left neighborhood of $t_0$) 
of the function $\Gamma_{log}$ being the right hand side 
of~\eqref{eq:stima_dGammalog} the sum of a positive function with an integrable one.

To prove the boundedness above of $\Gamma_{log}$ we cannot use the
boundedness of the energy function, indeed in this case we can just estimate 
$\Gamma_{log}(r,s) + M(t)\log r = h(t) - \frac12 \dot r ^2 $.
By the sake of contradiction suppose that $\Gamma_{log}$ diverges to $+\infty$
as 
$t$ tends to $t_0$; since  $U(t,rs) + M(t)\log r$ converges uniformly to
$\tilde U(t,s)$ as
$t$ tends to $t_0$ and $\tilde U(t,s)$ is a positive function, if $\Gamma_{log}$
diverges to $+\infty$ 
\begin{equation}
\label{assurdo}
\exists t_1 \in (t_0-\delta,t_0) \mbox{ such that }  \forall \, t \in (t_1,t_0), \quad
r^2 |\dot s|^2 > \max_{t \in [t_0-\delta,t_0]} M(t).
\end{equation}
From assumption~\eqref{assurdo} we have
\begin{equation}
\label{assurdo2}
\begin{split}
\int_{t_0-\delta}^{t_0} -\frac{\dot r}{r}\left( r^2 |\dot s|^2 - M(t) \right)dt & =
\int_{t_0-\delta}^{t_1} -\frac{\dot r}{r}\left( r^2 |\dot s|^2 - M(t) \right)dt +
\int_{t_1}^{t_0} -\frac{\dot r}{r}\left( r^2 |\dot s|^2 - M(t) \right)dt \\
& \geq \text{constant} - \lim_{t \to t_0} \log r(t) = +\infty.
\end{split}
\end{equation}
We now define the function 
\[
\Omega_{log}(r,s) := \Gamma_{log}(r,s) + M(t)\log r = h(t) -\frac12 \dot r^2
\]
that is bounded above.
When we compute its derivative with respect to the time variable we obtain
the sum of a positive function with an integrable one (we use
assumption~\eqref{assurdo}
and Lemma~\ref{le:integ_log}), indeed
\[
\begin{split}
\frac{d}{dt}\Omega_{log}(r,s) 
&= \frac{d}{dt}\Gamma_{log}(r,s) + \dot M(t)\log r + M(t)\frac{\dot r}{r}\\
&\geq -\frac{\dot r}{r} \left[ r^2 |\dot s|^2 -M(t) \right] 
- C_1 \left(U(t,rs)+1\right) + C_2 r^{\gamma}\frac{\dot r}{r}U(t,rs).
\end{split}
\]
We can then conclude the boundedness of $\Omega_{log}$ on the interval $[t_0-\delta,t_0]$
and from the estimate on its derivative we have
\[
\int_{t_0-\delta}^{t_0} -\frac{\dot r}{r}\left( r^2 |\dot s|^2 - M(t) \right)dt < +\infty
\]
that contradicts~\eqref{assurdo2}.
We conclude that the function $\Gamma_{log}$ is also bounded above.
\end{proof}

\begin{corollary}
\label{cor:lim_gamma_log}
As $t$ tends to $t_0$ the limit of the function $\Gamma_{log}$ exists finite
and 
\[
\lim_{t \to t_0^+} -\frac{\dot r^2}{2\log r}=M_0
\]
where $M_0:=M(t_0)$.
\end{corollary}
\begin{proof}
We argue as in the proof of Corollary~\ref{cor:limGamma}
to show that the function $\Gamma_{log}$ has a finite limit as 
$t$ tends to $t_0$. Since $\Gamma_{log} = h(t)- \frac12 \dot r^2 - M(t)\log r$,
we conclude dividing by $\log r$ using the boundedness of the function $h$.
\end{proof}

\begin{theorem}
\label{central_conf_log}
Let ${\bar x}$ be a generalized solution for the dynamical system~\eqref{DS}
and let $t_0\in(a,b)$ (in the case $b<+\infty$, $t_0$ can coincide with $b$) be a total
collision instant.
Let $r,s$ be the new variables defined in~(\ref{rs}); if the potential $U$ satisfies
assumptions \ref{a:U0}, \ref{a:U1}, \ref{a:U2l}, \ref{a:U3l} then the following assertions hold
\begin{itemize}
\item[(a)] $\displaystyle \lim_{t \rightarrow {t_0}^-} \left[ U(t,rs) + M(t)\log r
\right] = - \lim_{t \rightarrow {t_0}^-} \Gamma_{log}(r,s) = b$;
\item[(b)]
as $t$ tends to ${t_0}$, 
\[\begin{aligned}
r(t) & \sim
({t_0}-t)\sqrt{-2M_0\log({t_0}-t)} \\
\dot r(t) & \sim
-\sqrt{-2M_0\log({t_0}-t)}; 
\end{aligned}\]
% \item[(b)] as $t$ tends to ${t_0}$, $\displaystyle r(t) \sim
% ({t_0}-t)\sqrt{-2M_0\log({t_0}-t)}$;
% \item[(c)] as $t$ tends to ${t_0}$, $\displaystyle \dot r(t) \sim
% -\sqrt{-2M_0\log({t_0}-t)}$;
\item[(c)] $\displaystyle \lim_{t \rightarrow {t_0}^-} |\dot
s(t)|({t_0}-t)\sqrt{-2M_0\log({t_0}-t)} = 0$;
\item[(d)] for every real positive sequence $(\lambda)_n$, such that
$\lambda_n \rightarrow 0$ as $n \rightarrow +\infty$ we have
\[
\lim_{n \rightarrow +\infty} |s({t_0}-\lambda_n) - s({t_0}-\lambda_n t)|=0, 
\quad \forall t > 0.
\]
\end{itemize}
\end{theorem}
\begin{proof}
\emph{(a)} The proof is essentially the same of for 
Theorem~\ref{central_conf}.

\noindent \emph{(b)} From Corollary~\ref{cor:lim_gamma_log} we deduce that 
\[
\dot r(t) \sim -\sqrt{-2M_0\log r(t)}\quad \mbox{ as } t \mbox{ tends to }
t_0.
\]
We define $R(t) := ({t_0}-t)\sqrt{-2M_0\log({t_0}-t)}$ and we remark that, as
$t$ tends to $t_0$
\[
-\log R(t) = -\log({t_0}-t)-\log\left(\sqrt{-2M_0\log({t_0}-t)}\right) \sim
-\log({t_0}-t)
\]
and
\[
\dot R(t) = -\sqrt{-2M_0\log({t_0}-t)}
            +\frac{M_0}{\sqrt{-2M_0\log({t_0}-t)}}\sim -\sqrt{-2M_0\log
R(t)}.
\]
Our aim is then to prove that the function $r(t)$ is asymptotic to $R(t)$ 
as $t$ tends to $t_0$. We define the following functions
\[
f(\xi) := -\sqrt{-2M_0\log \xi} \quad \mbox{ and } \quad
\Phi(\xi):=\int_0^\xi \frac{d\eta}{f(\eta)}, \quad \xi \in(0,1] 
\]
and we remark that $\Phi(0)=0$ and $\Phi$ is a strictly decreasing function on
$[0,1]$. Moreover
\[
\dot r(t) \sim f(r(t)), \quad \dot R(t) \sim f(R(t)) \mbox{ as } t \mbox{ tends
to } t_0
\]
or equivalently
\[
\lim_{t \rightarrow t_0}\frac{d}{dt}\Phi(r(t)) = \lim_{t \rightarrow t_0}
\frac{d}{dt}\Phi(R(t)) = 1.
\]
Since the function $\Phi(\xi)$ decreases in $\xi$ and $r(t)$, 
$R(t)>0$ decreases in
$t$ (when we stay close to
the collision instant) we have that the functions $\Phi(r(t))$ and $\Phi(R(t))$
are negative 
on $(t_0-\delta_0,t_0)$, vanishes at $t_0$ (since $r(t_0)=R(t_0)=0$) and increase in the
variable $t$.
Furthermore fixed $\bar t < t_0$, the following property holds 
\begin{equation}
\label{eq:prop_Phi}
\begin{split}
\frac{d}{dt}\Phi(r(t)) \leq 1 \leq \frac{d}{dt}\Phi(R(t)), \, \forall t \in
(\bar t,t_0)
&\quad \Rightarrow \quad \Phi(r(t)) \geq \Phi(R(t)), \, \forall t \in (\bar
t,t_0) \\
&\quad \Rightarrow \quad r(t) \leq R(t),\, \forall t \in (\bar t,t_0).
\end{split}
\end{equation}
For every $\epsilon > 0$, we consider the functions 
\[
\begin{split}
& R^+_{\eps} (t):= (1+\eps)R(t),\\
& R^-_{\eps} (t):= (1-\eps)R(t).
\end{split}
\]
Since $\dot R(t) \sim f(R(t))$, we deduce that in a left neighborhood of $t_0$
\begin{equation}
\label{eq:log_pepenult}
\begin{split}
&\dot R^+_{\eps} (t) = (1+\eps)\dot R(t) \leq \left( 1 +
\frac{\eps}{2}\right)f(R(t))
                       \leq \left( 1 + \frac{\eps}{2}\right)f(R^+_{\eps}(t)),\\
&\dot R^-_{\eps} (t) = (1-\eps)\dot R(t) \geq \left( 1 -
\frac{\eps}{2}\right)f(R(t))
                                              \geq \left( 1 -
\frac{\eps}{2}\right)f(R^-_{\eps}(t)),
\end{split}
\end{equation}
indeed 
\[
f(R(t)) = -\sqrt{-2M_0\log(R(t))} \leq
-\sqrt{-2M_0\log(1+\eps)-2M_0\log(R(t))}=f(R^+_{\eps}(t))
\]
and similarly 
\[
f(R(t)) \geq -\sqrt{-2M_0\log(1-\eps)-2M_0\log(R(t))}=f(R^-_{\eps}(t)). 
\]

From~\eqref{eq:log_pepenult} we then obtain
\begin{equation}
\label{eq:log_penult}
\frac{d}{dt} \Phi(R^+_{\eps}(t)) \geq 1 + \frac{\eps}{2} \quad \mbox{and} \quad
\frac{d}{dt} \Phi(R^-_{\eps}(t)) \leq 1 - \frac{\eps}{2}.
\end{equation}
Moreover, since $\dot r(t) \sim f(r(t))$, again in a left neighborhood of $t_0$
we have that
\begin{equation}
\label{eq:log_ult}
\left( 1 + \frac{\eps}{2}\right)f(r(t)) \leq \dot r(t) \leq \left( 1 -
\frac{\eps}{2}\right)f(r(t))
\end{equation}
and dividing~\eqref{eq:log_ult} for the negative function $f(r(t))$ and
comparing the resulting inequalities with~\eqref{eq:log_penult} we have
\[
\frac{d}{dt} \Phi(R^-_{\eps}(t)) \leq \frac{d}{dt} \Phi(r(t)) \leq \frac{d}{dt}
\Phi(R^+_{\eps}(t)).
\]
From~(\ref{eq:prop_Phi}) we deduce that, in a neighborhood of the collision
instant $t_0$,
the following chain of inequalities holds
\[
(1-\eps) \leq \frac{r(t)}{R(t)} \leq (1+\eps).
\]
The second estimate follows directly.

\noindent \emph{(c)} From the result proved in \emph{(a)} we have that $\lim_{t
\to t_0}r|\dot s|=0$;
we conclude using \emph{(b)}.

\noindent \emph{(d)} As in the proof of Theorem~\ref{central_conf}, 
if $t=1$ there is nothing to prove. We then chose $t>0$, $t\neq 1$, 
a sequence $(\lambda_n)_n$, $\lambda_n \rightarrow 0$ and
$N$, sufficiently large, such that $\lambda_n < \delta/\max(1,t)$, $\forall n \geq N$.
We then obtain
\begin{eqnarray*}
|s({t_0}-\lambda_n)-s({t_0}-\lambda_n t)| 
& \leq & \int_{{t_0}-\lambda_n t}^{{t_0}-\lambda_n}|\dot s(u)|du \\
& \leq & \left( \int_{{t_0}-\lambda_n t}^{{t_0}-\lambda_n}-r(u)\dot r(u)|\dot
s(u)|^2du \right)^{\frac{1}{2}}
\left( \int_{{t_0}-\lambda_n t}^{{t_0}-\lambda_n}-\frac{du}{r(u)\dot r(u)}\right)^{\frac{1}{2}}.
\end{eqnarray*}
The boundedness of the $\Gamma_{log}$ and the estimate on its derivative
in~\eqref{eq:stima_dGammalog} 
imply the boundedness of the integral $\int_0^{t_0}r \dot r|\dot s|^2$ and then
\begin{equation*}
\lim_{n \rightarrow +\infty} \int_{{t_0}-\lambda_n t}^{{t_0}-\lambda_n}-r(u)\dot
r(u)|\dot s(u)|^2du =0.
\end{equation*}
Moreover, as $n$ tends to $+\infty$, from \emph{(b)} and \emph{(c)} we have 
$r(u)\dot r(u) \sim -2M_0(t_0-u)\log (t_0-u)$, hence 
\begin{equation*}
\lim_{n \rightarrow +\infty}
\int_{{t_0}-\lambda_n t}^{{t_0}-\lambda_n}\frac{du}{r(u)\dot r(u)} 
= \frac{1}{M_0}\lim_{n \rightarrow +\infty}\log \frac{\log \lambda_n t}{\log \lambda_n} = 0.
\end{equation*}
The proof is now complete.
\end{proof}

The behavior of the angular part is conserved also for logarithmic potential 
and the following result can be proved following the proof of 
Theorem~\ref{central_conf2}.
\begin{theorem}\label{central_conf2_log}
In the same setting of Theorem~\ref{central_conf_log},
assuming furthermore that the potential $U$ verifies
\begin{assiomi}
\item[\uqltag]
\label{a:U4l}
$\ds \lim_{r \to 0} r\nabla_{T} U(t,x)=\nabla_{T}\tilde U(t,s)$,
\end{assiomi}
then there holds
\[
\dlim_{t \rightarrow {t_0}} \mathrm{dist} \left({\cal C}^b,s(t)\right)=
\dlim_{t \rightarrow {t_0}} \inf_{\bar s \in {\cal C}^b}|s(t) - \bar s | = 0
\]
where ${\cal C}^b$ is the central configuration subset defined in~\eqref{eq:cal_C}.
\end{theorem}

%===========================
\section{Partial collisions}
\label{sec:partial}
%===========================
This section is devoted to the study of the singularities which are not total collision at the origin.
At first we shall prove the existence of a limiting configuration for bounded trajectories, that is
the Von Zeipel's Theorem  (stated on page \pageref{teo:Zeipel}). This fact allows the reduction from partial to
total collisions through a change of coordinates. To carry on the analysis we shall extend the clustering argument proposed by McGehee in \cite{mcgehee2} 
to prove the Von Zeipel's Theorem. To this aim we need to introduce some further assumptions on the potential $U$ and its singular 
set $\Delta$.
More precisely we suppose that
\begin{equation} \label{eq:strucDelta}
\Delta = \bigcup_{\mu\in\mathcal M} V_\mu,
\end{equation}
where the $V_\mu$'s are distinct linear subspaces of $\RR^{k}$ and $\mathcal M$ 
is a finite set;
observe that the set $\Delta$ is a cone, as required on page \pageref{cone}. We endow the 
family of the $V_\mu$'s with the inclusion partial ordering and we assume the family to be 
closed with respect to  intersection
(thus we are assuming that $\mathcal{M}$ is a semilattice
of linear subspaces of $\RR^k$: it is the intersection semilattice generated 
by the arrangement of maximal subspaces $V_\mu$'s).
With each $\xi\in\Delta$ we associate 
\[
\mu(\xi)=\min\{\mu\;:\; \xi\in V_\mu\}\qquad \text{i.e.,}
\qquad V_{\mu(\xi)}=\bigcap_{\xi\in V_\mu} V_\mu.
\]
Fixed $\mu\in {\mathcal M}$ we define the  set of collision configurations satisfying
\[
\Delta_\mu = \left\{ \xi \in \Delta : \mu(\xi)=\mu \right\}
\]
and we observe that this is an open subset of $V_\mu$ and its closure 
$\overline{\Delta_\mu}$ is 
$ V_\mu$. We also notice that the map $\xi\to \dim(V_{\mu(\xi)})$ is lower semicontinuous.

We denote by $p_\mu$ the orthogonal projection onto $V_\mu$ and we write
\[
x=p_\mu(x)+w_\mu(x),
\]
where, of course,  $w_\mu=\mathbb I-p_\mu$.

We assume that, near the collision set, the potential depends, roughly, only on the 
projection orthogonal to the collision set: more precisely we assume

\begin{assiomi} 
\item[\uctag]
\label{a:U5}
For every $\xi\in\Delta$, there is $\varepsilon>0$ such that
\[
\displaystyle U(t,x) -U(t,w_{\mu(\xi)}(x))= W(t,x) 
\in {\cal C}^1\left((a,b)\times B_\eps(\xi)\right),
\]
where $B_\eps(\xi)=\{x\;:\: |x-\xi|<\varepsilon\}$.
\end{assiomi}

%
% We moreover assume the boundedness of the potential away from the singularities
% \begin{enumerate}
% \item[\us] $\displaystyle \forall \eps>0, \;\forall C>0, \; \exists M>0 :
% \left(d(x,\Delta)>\eps \mbox{ and } I(x)\leq C 
% \implies |\nabla U(t,x)| \leq M\right)$.
% where $d(x,\Delta)$ denotes the distance of $x$ from $\Delta$.
% \end{enumerate}

\begin{theorem}\label{thm:partial}
Let $\bar x$ be a generalized solution for the dynamical system~\eqref{DS} on the 
bounded interval $(a,b)$.
Suppose that the potential $U$ satisfies assumptions \ref{a:U0}, \ref{a:U1}, \ref{a:U5}, and 
\ref{a:U2h}, \ref{a:U3h}, \ref{a:U4h} (or \ref{a:U2l}, \ref{a:U3l}, \ref{a:U4l}).\\
If $\bar x$
is bounded on the whole interval $(a,b)$ then 
\begin{enumerate}
\item[(a)] $\bar x$ has a finite number of  singularities which are collisions 
(the Von Zeipel's Theorem holds).
\item[(b)] Furthermore, if $t^* \in \bar x^{-1}(\Delta)$ is a collision instant, 
$x^*$ the limit configuration  of $\bar x$ as $t$ tends to $t^*$ and
$\mu^*=\mu(x^*) \in {\mathcal M}$, then 
 $r_{\mu^*}=|w_{\mu^*}(\bar x)|$, $s_{\mu^*}=w_{\mu^*}(\bar x)/r_{\mu^*}$ and $U_{\mu^*}=
 U(t,w_{\mu^*}(\bar x))$ satisfy the asymptotic estimates given in
Theorems~\ref{central_conf} and~\ref{central_conf2} 
(or Theorems~\ref{central_conf_log} and~\ref{central_conf2_log} 
when \ref{a:U2l}, \ref{a:U3l} and \ref{a:U4l} hold).
\end{enumerate}
\end{theorem}

\begin{proof}
Let $\bar x$ be a generalized solution with a 
singularity at $t=t^*$ (see Definition~\ref{def:singularity}) and $\Delta^*$ its 
$\omega$-limit set, that is
\[
\Delta^* = \left\{ x^* : \exists (t_n)_n \mbox{ such that } t_n \to t^*
\mbox{ and } \bar x(t_n)\to x^* \right\}.
\]
It is well known that the  $\omega$-limit of a bounded trajectory is a compact and connected set. From the Painlev\'e's Theorem (on page \pageref{teo:Painleve}) we have the inclusion
\[
\Delta^*  \subset \Delta.
\]
Von Zeipel's Theorem asserts that whenever  $\bar x$ remains
bounded as $t$ approaches $t^*$, then the $\omega$-limit set of $\bar x$ 
contains just one element, that is $\Delta^* = \{ x^* \}$.

In view of Corollary~\ref{cor:2}, where we proved the Theorem in the case 
$\liminf_{t \to t^*} \dot I(\bar x(t))$ $< +\infty$, we are left with the case
when
\begin{equation*}
\lim_{t \to t^*} \dot I(\bar x(t)) = +\infty.
\end{equation*}
From this and our  assumptions it follows that $I(\bar x(t))$ is a 
definitely increasing and bounded 
function. Hence it admits a limit
\begin{equation}
\label{hp:liminf}
\lim_{t \to t^*} I(\bar x(t)) = I^*.
\end{equation}

We perform the proof of Von Zeipel's Theorem in two steps.

\bigskip

\noindent 
{\em Step 1. We suppose that $\mu\left(\Delta^*\right) = \{\mu^*\}$ for some 
$\mu^* \in {\mathcal M}$ and we show that $\Delta^* = \{ x^* \}$}.

\medskip

\noindent 

As  $\Delta^*$ is a compact and 
connected subset of $V_{\mu^*}$, we have the following inclusions 
\[
\Delta^* \subset \Delta_{\mu^*} \subset V_{\mu^*}.
\]
% Since $\Delta^*$ is compact and  $\Delta_{\mu^*}$ is open in $V_{\mu^*}$, 
% there exists $\eps>0$ such that
% the $\eps$-neighborhood of $\left(\Delta^*\right)^\eps$ does not intersect
% $V_{\mu^*}\minus\Delta_{\mu^*}$.
We consider the orthogonal projections
\[
p(t)=p_{\mu^*}(\bar x(t)),\qquad\qquad w(t)=w_{\mu^*}(\bar x(t)).
\]
Since we have assumed that $\mu(\Delta^*)=\{\mu^*\}$, then
\begin{equation}\label{eq:bar_coor}
\lim_{t \to t^*} w(t) = 0,
\end{equation}
our aim is now to prove that
\[
\lim_{t \to t^*} p(t) = x^*.
\]

Projecting on $V_{\mu^*}$ the equations of motion, we obtain from \ref{a:U5}
\begin{equation}\label{eq:DScenter}
-\ddot p=p_{\mu^*}\left(\nabla U(t,\bar x(t))\right)=p_{\mu^*}\left(\nabla W(t,\bar x(t))\right)
\end{equation}
where $\nabla W$ is globally bounded as $t\to t^*$. Indeed, fixed $\eps>0$,
%such that $B_\eps\left( \Delta^* \right) \cap \left( V_{\mu^*}\minus 
%\Delta_{\mu^*}\right) = \emptyset$, 
there exists $\delta>0$ such that 
$\bar x(t)\in B_\eps\left( \Delta^* \right)$ whenever $t \in (t^*-\delta,t^*)$, and
from assumption \ref{a:U5} and the compactness of $\Delta^*\subset\Delta_{\mu^*}$ 
we deduce the boundedness of the right hand 
side of~\eqref{eq:DScenter}. From this fact we easily deduce the existence of a limit for $(p(t))$ as 
$t$ tends to $t^*$.
A word of caution must be entered at this point. As $\bar x$ is a generalized 
solution to~\eqref{DS}, the
equation of motions are not available, because of the possible occurence of collisions, and 
therefore they can not be projected on $V_{\mu^*}$. Nevertheless, exploiting the 
regularization method exposed in Section~\ref{sec:sing} and projecting the regularized 
equations, one can easily obtain the validity of~\eqref{eq:DScenter} after passing to the 
limit.
\bigskip

\noindent 
{\em Step 2. There always exists $\mu^* \in {\cal M}$ such that 
$\mu\left(\Delta^*\right) = \{\mu^*\}$.}

\medskip

\noindent Let $\mu^*$ be the element of $\mu\left(\Delta^*\right)$ associated with the
subspace $V_{\mu^*}$ having \emph{minimal dimension}. Since the function 
$\xi\to \dim(V_{\mu(\xi)})$ is lower semicontinuous, 
the minimality of the dimension has as a main implication that 
$\Delta_{\mu^*}\cap\Delta^*$ is compact. Hence the function $\nabla W$ appearing in \ref{a:U5} can be
though to be globally bounded in a neighborhood of $\Delta_{\mu^*}\cap\Delta^*$.
In other words, when considering the orthogonal projections
$p(t)=p_{\mu^*}(\bar x(t))$ and $w(t)=w_{\mu^*}(\bar x(t))$, 
as a major consequence of the minimality of the dimension $\mu^*$ we find the following
implication: 
\begin{equation}\label{eq:min1}
\exists M>0, \; \exists \eps>0 : |w(t)|^2 < \eps
\implies|p_{\mu^*}\left(\nabla W(t,\bar x)\right)| \leq M .
\end{equation}
We now compute the second derivative (with respect to the time $t$) 
of the function $|p(t)|^2$
\[
\dfrac{d^2}{d t^2} |p(t)|^2  = 2\ddot p(t)\cdot p(t)
+ 2 \dot p(t)\cdot \dot p(t)  \geq 
-2 p_{\mu^*}\left(\nabla W(t,\bar x(t))\right)\cdot p(t)
\]
Thus, from the projected motion equation ~\eqref{eq:DScenter} and from~\eqref{eq:min1} we infer
\begin{equation}\label{eq:min2}
\exists K>0, \; \exists \eps>0 : |w(t)|^2 < \eps
\implies  \dfrac{d^2}{d t^2}|p(t)|^2 \geq -K .
\end{equation}
We now argue by contradiction, supposing that $\mu\left(\Delta^*\right) \neq \{\mu^*\}$. 
Then
\begin{equation}\label{eq:inertia2}
0 = \liminf_{t \to t^*}|w(t)|^2 < \limsup_{t \to t^*}|w(t)|^2.
\end{equation}
Since, obviously, the total moment of inertia splits as
\[
I(\bar x(t))=|p(t)|^2+|w(t)|^2,
\]
from~\eqref{eq:inertia2} and~\eqref{hp:liminf} we deduce that
\begin{equation}\label{eq:HPassurdo}
I^* = \limsup_{t \to t^*}|p(t)|^2 > \liminf_{t \to t^*}|p(t)|^2
\end{equation}
and from~\eqref{eq:HPassurdo} together with~\eqref{eq:min2} we have
\[
\exists K>0, \; \exists \eps>0 \;:\;  |p(t)|^2 \geq I^* - \eps
\implies  \dfrac{d^2}{d t^2} |p(t)|^2\geq -K .
\]
Let $(t_n^0)_n$ and $(t_n^*)_n$ be two sequences such that, fixed $\eps>0$
\[
\begin{split}
& t_n^*<t_n^0<t_{n+1}^* \quad \forall n\\
& t_n^0 \to t^* \quad t_n^* \to t^* \mbox{ as } n \to +\infty \\
& f(t_n^*) \to I^* \mbox{ as } n \to +\infty  \mbox{ and } f'(t_n^*)=0, \quad \forall n \\
& t_n^0=\inf\{t>t^*_n : |p(t)|^2 \leq I^*-\eps\}, \quad \forall n. 
\end{split}
\]
Hence $|p(t_n^0)|^2-|p(t_n^*)|^2 = \dfrac{d}{d t^2}|p(\xi)|^2(t_n^0-t_n^*)^2/2 \geq -K(t_n^0-t_n^*)^2/2$ and then 
\[
-\eps \geq \frac{-K}{2}(t_n^0-t_n^*)^2 \quad \mbox{or} \quad
(t_n^0-t_n^*)^2 \geq \frac{2\eps}{K}
\]
in contradiction with the assumptions that both sequences $(t_n^0)_n$ and $(t_n^*)_n$
tend to the finite limit $t^*$.
This concludes the proof of the Von Zeipel's Theorem. Next we prove isolatedness of collision instants.

To this aim, let us select  $t^* \in \partial\left(\bar x^{-1}(\Delta)\right)$  a collision instant
such that the dimension of $V_{\mu(\bar x(t^*))}$ is minimal among all dimensions of collision configurations $V_{\mu(\bar x(t))}$
in $(t^*-\delta,t^*+\delta)$ for some $\delta>0$. As before, let us split the components of the trajectory $\bar x(t)= p(t)+w(t)$ on $V_{\mu^*}$ and its
orthogonal complement.

Since $\mu^*$ is minimal (see~\eqref{eq:min1}), we already know from the previous discussion that the equations of motion 
projected on the subspace $V_{\mu^*}$ (equation~\eqref{eq:DScenter}) are not singular; 
on the other hand, by \ref{a:U5}, the trajectories in the orthogonal coordinates $w$ are 
generalized solutions to a dynamical system of the form
\begin{equation}\label{eq:DSproj}
-\ddot w = \nabla U(t,w)+\nabla W(t,p(t)+w).
\end{equation}
Now, since $w(t)$ has a total collisions at the origin at $t^*$, we can apply the results of Section 
\ref{sec:asymp}. More precisely, at first we deduce 
from Theorem~\ref{thm:totcolliso} that $t^*$ is isolated in the set of collisions $\Delta_{\mu^*}$; 
furthermore from Corollary~\ref{cor:2} we deduce the boundedness of the action 
 and the energy. Finally we conclude applying 
Theorems~\ref{central_conf},~\ref{central_conf2} (or Theorems~\ref{central_conf_log}, 
\ref{central_conf2_log} when \ref{a:U2l}, \ref{a:U3l} and \ref{a:U4l} hold) to the projection $w$. 
In particular from (a) in Theorem~\ref{central_conf} (or Theorem~\ref{central_conf_log})
we obtain that every collision is isolated and hence, whenever the interval $(a,b)$ is finite,
the existence of a finite number of collisions.
\end{proof}

%==============================
\section{Absence of collisions for locally minimal path}
\label{sec:blowup}
%==============================

As a matter of fact, solutions to the Newtonian $n$--body problem which are
minimals for the action are, very likely, free of any collision. This  was
discovered  in \cite{serter1} for a class of periodic three--body problems and,
since then, widely exploited in the literature concerning the variational
approach to the periodic $n$--body problem. In general, the proof goes by the
sake of the contradiction and involves the construction of a suitable variation
that lowers the action in presence of a collision.  A recent breakthrough in
this direction is due of the neat idea, due to C. Marchal in 
\cite{marchal}, 
of averaging over a family of variations parameterized on a sphere. The method
of averaged variations for Newtonian potentials has been developed and exposed
in \cite{Ch2}, and then extended to $\alpha$--homogeneous potentials and
various constrained minimization problems in \cite{FT}. This argument can be
used in most of the known cases to prove that minimizing trajectories are
collisionless. In this section we prove the absence of collisions for locally
minimal solutions when the potentials have  quasi--homogeneous or logarithmic
singularities.

We consider separately the quasi--homogeneous and the logarithmic cases;
indeed in the first case one can exploit the blow--up technique as
developed in Section 7 of \cite{FT};
in \S~\ref{sec:blowup_hom} we will just recall the main steps of this
arguments.  On the other hand, when dealing with logarithmic potentials, the
blow--up technique is no longer available and we conclude proving directly some
averaging estimates that can be used to show the nonminimality of large classes
of colliding motions.

%==============================
\subsection{Quasi--homogeneous potentials}
\label{sec:blowup_hom}
%==============================
Let $\tilde U$ be the $\cont^1$ function defined on $(a,b) \times ({\cal E}\minus\Delta)$
introduced on page \pageref{tildeU}; we extend its definition on the whole 
$(a,b) \times (\RR^{k}\minus\Delta)$ in the following way
\[
\tilde U(t,x)= |x|^{-\alpha}\tilde U(t,x/|x|).
\]
Fixed $t^*$ (in this section we will consider a locally minimal trajectory $\bar x$ 
with a collision at $t^*$) in this section, with an abuse of notation, we denote
\begin{equation}
\label{eq:tildeUristr}
\tilde U(x)=\tilde U(t^*,x).
\end{equation}
Of course, the function $\tilde U$ is homogeneous of degree $-\alpha$
on $\RR^{k}\minus\Delta$.
\begin{theorem}\label{theo:nocoll}
In addition to \ref{a:U0}, \ref{a:U1}, \ref{a:U2h}, \ref{a:U3h}, \ref{a:U4h}, \ref{a:U5}, 
assume that, for a given $\xi\in\Delta$ 
\begin{assiomi}
\item[\useitag]
\label{a:U6}
there is a $2$--dimensional linear subspace of $V_{\mu(\xi)}^{\perp}$, say $W$,
where $\tilde U$ is rotationally invariant;
%\[\tilde U(e^{i\theta}w)=\tilde U(w) , \qquad\forall w\in W, \forall \theta\in[0,2\pi];\]
\item[\usettehtag]
\label{a:U7h}
for every $x\in\RR^k$ and $\delta \in W$ there holds
\[
\tilde U(x+\delta) \leq \tilde U\left(\left(\dfrac{\tilde U(\pi_W(x))}{\tilde U(x)}\right)^{1/\alpha}\pi_W(x)+\left(\dfrac{\tilde U(x)}{\tilde U(\pi_W(x))}\right)^{1/\alpha}\delta\right)
\]
where $\pi_W$ denotes the orthogonal projection onto $W$.
\end{assiomi}
Then generalized solutions do not have collisions at the 
configuration $\xi$ at the time $t^*$.
\end{theorem}

\begin{figure}
\caption{Potential levels, with $\lambda=\left(\dfrac{\tilde U(\pi_W(x))}{\tilde U(x)}\right)^{1/\alpha}>1$}
\begin{center}
\psfrag{x}[c]{$x$}
\psfrag{xpd}{$x+\delta$}
\psfrag{x}[c]{$x$}
\psfrag{W}[c]{$W$}
\psfrag{px}[c]{$\pi_W(x)$}
\psfrag{xpx}{$\lambda\pi_W(x)+\lambda^{-1}\delta$}
\psfrag{lampx}[c]{$\lambda\pi_W(x)$}
\psfrag{Ux}[c]{$\tilde U(x)$}
\psfrag{Ulamx}[c]{$\tilde U(\lambda \pi_W(x) + \lambda^{-1}\delta)$}
\psfrag{del}{$\delta$}
\psfrag{lamdel}{$\lambda^{-1}\delta$}
\includegraphics[height=0.3\textheight]{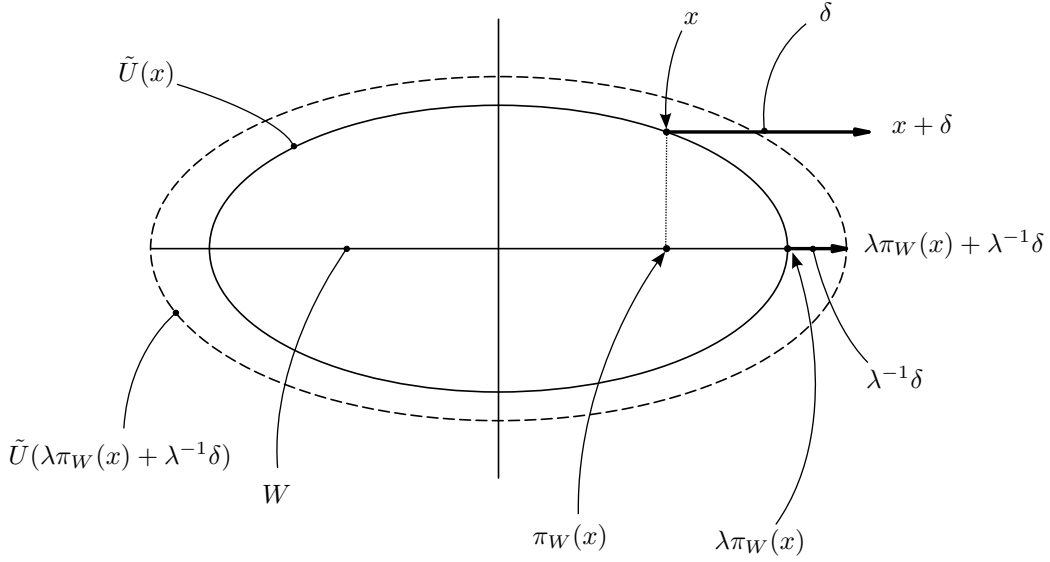}
\end{center}
\end{figure}

\begin{remark}
Some comments on assumptions \ref{a:U6} and \ref{a:U7h} are in order. Of course, as our potential 
$\tilde U$ is homogeneous of degree $-\alpha$ the function 
\[
\varphi(x)=\tilde U^{-1/\alpha}(x)
\]
is a non negative, homogeneous of degree one function, having now $\Delta$ as zero set. In most 
of our  applications $\varphi$ will be indeed a quadratic form.  
Assume that $\varphi^2$ splits in the following way:
\[
\varphi^2(x)=K|\pi_W(x)|^2+\varphi^2(\pi_{W^\perp}(x))
\]
for some positive constant $K$.
Then \ref{a:U6} and \ref{a:U7h} are satisfied. Indeed, denoting $w=\pi_W(x)$ and $z=x-w$ we have,
for every $\delta \in W$,
\[
\begin{split}
\varphi^2(x+\delta)&=K\vert w+\delta\vert^2+\varphi^2(z)\\
& = K\left|{\frac{\varphi(x)}{\varphi(w)}}w +
{\frac{\varphi(w)}{\varphi(x)}}\delta\right|^2 + 
K\frac{\varphi^2(z)}{\varphi^2(x)}|\delta|^2\\
&\geq K\left|{\frac{\varphi(x)}{\varphi(w)}}w +
{\frac{\varphi(w)}{\varphi(x)}}\delta\right|^2\\
& = \varphi^2\left({\frac{\varphi(x)}{\varphi(w)}}w +
{\frac{\varphi(w)}{\varphi(x)}}\delta\right),
\end{split}
\]
which is obviously equivalent to \ref{a:U7h}. Therefore, we have the following Proposition.
\end{remark}
\begin{proposition}\label{prop:usette}
Assume $\tilde U(x)=\mathcal Q^{-\alpha/2}(x)$ for some non negative quadratic form 
$\mathcal Q(x) = \langle Ax,x\rangle$. Then assumptions~\ref{a:U6} and~\ref{a:U7h} are satisfied whenever $W$ is 
included in an eigenspace of $A$ associated with a multiple eigenvalue.
\end{proposition}

\begin{remark}
Given two potentials satisfying \ref{a:U6} and  \ref{a:U7h}  for a common
subspace $W$, their sum enjoys the same properties. On the other hand, if they
do not admit a common subspace $W$, their sum does not satisfy \ref{a:U6} and
\ref{a:U7h}.
\end{remark}

% Let us start the proof or Theorem~\ref{theo:nocoll}.
\begin{proof}[Proof or Theorem~\ref{theo:nocoll}]
% Assume not and 
Let  $\bar x(t)$ be a generalized solution with a collision at the
time $t^*$, i.e.~ $\bar x(t^*)=\xi \in \Delta$; 
up to time-translation we assume that  the collision
instant is $t^*=0$. 
Furthermore, using the same arguments needed in the proof of the Von Zeipel's
Theorem in Section~\ref{sec:partial}, we can suppose that $\xi=0$.  We consider
the case of a boundary collision  (interior collisions can be treated in a
similar way).  Then Theorem~\ref{central_conf} 
ensures the existence of $\delta_0>0$ such that no other collision occurs in some interval
$[0,\delta_0]$.

We consider the family of rescaled generalized solutions 
\[
\bar x^{\lambda_n}(t) := \lambda_n^{-\frac{2}{2+\alpha}}\bar x(\lambda_n t),
\quad t \in [0,\delta_0/\lambda_n].
\]
where  $\lambda_n \rightarrow 0$ as $n \rightarrow +\infty$. 
From the asymptotic estimates of Theorem~\ref{central_conf2} we know that 
the angular part $(s(\lambda_n))_n$ converges, up to subsequences, to some central 
configuration $\bar s$, in particular $\bar s$ is in the $\omega$--limit of $s(t)$.

For any $\bar s$ in the $\omega$--limit of $s(t)$, 
a (right) blow-up of $\bar x$ in $t=0$ is a path defined, for $t \in [0,+\infty)$, as
\begin{equation}
\label{eq:blow-up}
\bar q(t) := \zeta t^\frac{2}{2+\alpha},\qquad\zeta=K\bar s,
\end{equation}
where the constant $K>0$ is determined by part \emph{(b)} of
Theorem~\ref{central_conf}. We note that the blow--up is a homothetic solution to the 
dynamical system associated with the homogeneous potential $\tilde U$ and that 
it has zero energy (the blow--up is parabolic).
If $s(\lambda_n) \to \bar s$ as $n \to+\infty$,
from Theorem~\ref{central_conf}, we obtain straightforwardly 
the pointwise convergence of $\bar x^{\lambda_n}$ to the blow up $\bar q$ and the
$H^1$-boundedness of $\bar x^{\lambda_n}$ implies
its uniform convergence on compact subsets of $[0,+\infty)$. Furthermore 
the convergence holds locally in the $H^1([0,+\infty))$--topology. Finally
also the 
sequence $\dot {\bar x}^{\lambda_n}$ converges uniformly on every interval $[\eps,T]$,
with arbitrary $0<\eps<T$.

The following fact has been proven in \cite{FT}, Proposition 7.9.

\begin{lemma}
Let $\bar x$ be a locally minimizing trajectory with a total collision at $t=0$ and let $\bar 
q$ be its blow--up in $t=0$. Then $\bar q$ is a locally minimizing trajectory for the 
dynamical system associated with the homogeneous potential $\tilde U$ introduced in 
\eqref{eq:tildeUristr}.
\end{lemma}

We will conclude the proof showing that $\bar q$ cannot be a locally 
minimizing trajectory for the  dynamical system associated with $\tilde U$.
Following \cite{FT}, we now introduce a class of suitable variations as follows:

\begin{definition}
\label{def:stand_var}
The standard variation associated with $\delta \in \RR^k\minus\{0\}$ is defined as 
\begin{equation*}
v^\delta(t)= \left\{
\begin{array}{ll}
\delta & \mbox{if} \,\,\, 0 \leq |t| \leq T-|\delta| \\
(T-t)\frac{\delta}{|\delta|} & \mbox{if} \,\,\, T-|\delta| \leq |t| \leq T \\
0 & \mbox{if} \,\,\, |t| \geq T, \\
\end{array}\right.
\end{equation*}
for some positive $T$.
\end{definition}

We wish to estimate the action differential corresponding to a standard variation. 
To this aim we give the next definition.
\begin{definition} The displacement potential differential associated with $\delta \in \RR^k$ 
is defined as:
$$
S(\zeta,\delta)=\int_0^{+\infty} 
\left( \tilde U(\zeta t^{2/(2+\alpha)}+\delta) -
\tilde U( \zeta t^{2/(2+\alpha)})  \right) dt
$$
where $\bar q(t) = \zeta t^{2/(2+\alpha)} $ is a blow-up of $\bar x$ in $t=0$.
\end{definition}
The quantity $S(\zeta,\delta)$ represents the \emph{potential differential} needed for 
displacing the  colliding trajectory originarily traveling along the $\zeta$--direction  
to  the point $\delta$.
It has been proven in \cite{FT} Proposition 9.2, that the function $S$ represents 
the limiting behavior, as $\delta\to0$, of the whole action differential:
\[
\Delta \mathcal{A}^{\delta} := 
\int_{-\infty}^{+\infty} \left[K(\dot{\bar q}+\dot v^\delta)+\tilde U(\bar q+v^\delta) - 
K(\dot{\bar q})-\tilde U(\bar q)\right] dt.
\]
Indeed, the fundamental estimate holds: 
\begin{lemma}\label{lemma:act_decreas}
Let $\bar q=  \zeta t^{2/(2+\alpha)}$ be a blow--up
trajectory and $v^\delta$ any standard variation. 
Then, as $\delta \rightarrow 0$
\begin{equation*}
\Delta {\mathcal A}^{\delta} 
= |\delta|^{1-\alpha/2}  S\left(\zeta, \frac{\delta}{|\delta|}\right) + O(|\delta|).
\end{equation*}
\end{lemma}

We observe that, from the homogeneity of $\tilde U$ it follows that
\begin{equation}\label{eq:homog}
S(\lambda\xi,\mu\delta)=  
\left|\lambda\right|^{-1-\alpha/2}|\mu|^{1-\alpha/2}S(\xi,\delta)
\end{equation}
(see \cite[(8.2)]{FT}) and hence, 
if $\tilde U$ is invariant under rotations, the sign of $S$ depends only on the 
angle between $\xi$ and $\delta$. To deal with the isotropic case (which is not the case 
here), the following function was introduced in \cite{FT}:
 $$\Phi_\alpha(\vartheta)=\int_0^{+\infty} 
\frac{1}{\left( t^{\frac{4}{\alpha+2}} 
-2\cos \vartheta t^{\frac{2}{\alpha+2}} 
+1 \right)^{\alpha /2}} - 
\frac{1}{ t^{\frac{2\alpha}{\alpha+2}}} dt.
$$
The value of $\Phi_\alpha(\vartheta)$ ranges from  positive to negative values, 
depending on $\vartheta$ and $\alpha$. Nevertheless,  it is always negative, 
when averaged on a circle. 
Indeed, the following inequality was obtained in \cite[Theorem 8.4]{FT}.

\begin{lemma}\label{fund_ineq} For any $\alpha \in(0,2)$ there holds
\[
\dfrac{1}{2\pi}\int_0^{2\pi}\Phi_\alpha(\vartheta)d\vartheta<0.
\]
\end{lemma}

This inequality will be a key tool in proving the following averaged estimate:
\begin{lemma}\label{lemma:ineq}
Assume \ref{a:U6} and \ref{a:U7h}, then, if $\SSS$ is the unitary circle of $W$, for any 
$\zeta\in\RR^k\minus\{0\}$ the following inequality holds
\[
% \frac{1}{|\SSS|}
\int_{\SSS} S(\zeta,\delta)d\delta<0\;.
\]
As a consequence,
\[
\forall \zeta\in\RR^k\minus\{0\}\;\exists \delta=\delta(\zeta)\in\SSS\;:\; 
S(\zeta,\delta(\zeta))<0\;.
\]
\end{lemma}
\begin{proof} 
As a first obvious application of Lemma~\ref{fund_ineq} we obtain the assertion \emph{for any 
$\zeta\in W\minus \{0\}$}. Indeed, by~\eqref{eq:homog} and \ref{a:U6} we easily obtain
\[
\zeta\in W \minus \{0\} \quad \Longrightarrow \quad
S(\zeta,\delta) = K\left|\zeta\right|^{-1-\alpha/2}\Phi_\alpha(\vartheta),
\]
where $K$ is a positive constant and $\vartheta$ denotes the angle between $\zeta$ 
and $\delta$.

Now we prove the assertion for any $\zeta \neq 0$ in the configuration space.
It follows from  the homogeneity of $\tilde U$ that
\[
\tilde U\left(\left(\dfrac{\tilde U(\pi_W(\zeta))}{\tilde 
U(\zeta)}\right)^{1/\alpha}\pi_W(\zeta)\right)=\tilde U(\zeta).
\]
Hence \ref{a:U7h} implies, for every $\delta\in\SSS$,
\[S(\zeta,\delta)\leq S\left(\left(\dfrac{\tilde U(\pi_W(\zeta))}{\tilde 
U(\zeta)}\right)^{1/\alpha}\pi_W(\zeta),\left(\dfrac{\tilde U(\zeta)}{\tilde 
U(\pi_W(\zeta))}\right)^{1/\alpha}\delta\right).
\]
Hence~\eqref{eq:homog} implies
\[S(\zeta,\delta)\leq \left(\dfrac{\tilde U(\zeta)}{\tilde U(\pi_W(\zeta))}\right)^{2/\alpha} 
S(\pi_W(\zeta),\delta),
\]
and thus 
\[
% \frac{1}{|\SSS|}
\int_{\SSS} S(\zeta,\delta)\,d\delta
\leq 
\left(\dfrac{\tilde U(\zeta)}{\tilde 
U(\pi_W(\zeta))}\right)^{2/\alpha}
% \frac{1}{|\SSS|}
\int_{\SSS} 
S(\pi_W(\zeta,\delta))\,d\delta<0\;.
\]
\end{proof}
\noindent\emph{End of the Proof of Theorem~\ref{theo:nocoll}.}  
To conclude the proof, according with Lemma~\ref{lemma:ineq} we chose 
$\delta=\delta(\zeta)\in  W\minus\{0\}$ with the property that 
$S(\zeta,\delta(\zeta)/|\delta(\zeta)|)<0$. 
As a consequence of Lemma~\ref{lemma:act_decreas}, we can 
lower the value of the action of $\bar q$ by performing the standard variation 
$v^{\delta(\zeta)}$, provided the norm of $|\delta(\zeta)|$ is sufficiently small
(in order to apply Lemma \ref{lemma:act_decreas}). 
Hence $\bar q$ can not be locally minimizing for the action.
\end{proof}

As we have already noticed, the class of potentials satisfying \ref{a:U6} and \ref{a:U7h} is not 
stable with respect to the sum of potentials. In order to deal with a class of potentials
which is closed with respect to the sum, we introduce the following variant of Theorem 
\ref{theo:nocoll}.

\begin{theorem}\label{theo:nocoll_2}
In addition to \ref{a:U0}, \ref{a:U1}, \ref{a:U2h}, \ref{a:U3h}, \ref{a:U4h}, \ref{a:U5}, 
assume that $\tilde U$ has the form
\[
\tilde U(x)=\sum_{\nu=1}^N\frac{K_\nu}{\left({\rm dist}(x,V_\nu)\right)^\alpha}
\]
where $K_\nu$ are positive constants and $V_\nu$ is a family of linear subspaces, 
with ${\rm codim}(V_\nu)\geq 2$, for every $\nu=1,\dots,N$.
Then locally minimizing trajectories do not have collisions at the time $t^*$.
\end{theorem}

\begin{proof}
Following the  arguments of the proof of Theorem~\ref{theo:nocoll}, the
assertion will be proved once we show, as in Lemma~\ref{lemma:ineq}, that, for
every index $\nu$, there holds
\begin{equation*}\label{eq:int_sfera}
% \frac{1}{|\SSS^{k-1}|}
\int_{\SSS^{k-1}} S_\nu(\zeta,\delta) \,d\delta<0,
\end{equation*}
where, of course, we denote 
\[
S_\nu(\zeta,\delta)=\int_0^{+\infty} \left(
{\rm dist}(\zeta t^{2/(2+\alpha)}+\delta, V_\nu)^{-\alpha} -
{\rm dist}(\zeta t^{2/(2+\alpha)}, V_\nu)^{-\alpha} \right) dt
\]
and $\SSS^{k-1}$ is the unit sphere of the configuration space $\RR^k$. 
This is an elementary consequence of Lemma \ref{lemma:ineq} and the  fact 
that the function $S_\nu(\zeta,\delta)$ only depends on the projection of $\zeta$ 
orthogonal to $V_\mu$ and has rotational invariance on $V_\nu^\perp$. 
Thus the integral of $S_\nu$ over the sphere is a positive multiple of its integral 
on any circle $\SSS$ orthogonal to $V_\nu$.
\end{proof}

%==============================
\subsection{Logarithmic type potentials}
\label{sec:blowup_log}
%==============================
In this section we prove the equivalent to Theorems \ref{theo:nocoll} and \ref{theo:nocoll_2}
suitable for logarithmic type potentials. Concerning the quasi-homogeneous 
case we have seen that a crucial role is played by the construction of a blow--up 
function which minimizes a limiting problem.
Before starting, let us highlight the reasons why, when dealing
with logarithmic potentials,  a blow-up limit can not exist.
Indeed, the natural scaling
should be $\bar x^{\lambda_n}(t) := \lambda_n^{-1}\bar x(\lambda_n t)$,
which does not converge, since
\[
\lim_{\lambda_n \to 0} \bar x^{\lambda_n}(t) =
\lim_{\lambda_n \to 0} \frac{r(\lambda_n t)s(\lambda_n t)}
                           {\lambda_n t\sqrt{-2M(0)\log(\lambda_n t)}} 
t\sqrt{-2M(0)\log(\lambda_n t)}
=+\infty
\]
for every $t>0$. On the other, hand, looking at~\ref{eq:blow-up},
the (right) blow--up should be, up to a change of time
scale,
\begin{equation}
\label{eq:blow-up_log}
\bar q(t) := t\bar s, \quad i \in {\bf k},
\end{equation}
where $\bar s$ is a central configuration for the system
limit of a sequence $s(\lambda_n)$ where  $(\lambda_n)_n$ is such that
$\lambda_n \rightarrow 0$.
The blow up function defined in~(\ref{eq:blow-up_log})
is the pointwise limit of the normalized sequence 
\[
\bar x^{\lambda_n}(t) := \frac{1}{\lambda_n \sqrt{-2M(0)\log\lambda_n}}\bar
x(\lambda_n t).
\]
Unfortunately the path  in~(\ref{eq:blow-up_log}) is not locally minimal for the
limiting problem, indeed since, the sequence $(\ddot{\bar x}^{\lambda_n})_n$ 
converges to $0$ as $n$ tends to $+\infty$, the blow-up in~(\ref{eq:blow-up_log})
minimizes only the kinetic part of the action functional.

We shall overcome this difficulty by proving the averaged estimate 
in a direct way from the asymptotic  estimates of Theorem~\ref{central_conf_log}
and assuming~\eqref{eq:hp_pot_log} on the potential $U$.
As we have done for the quasi--homogeneous case, we extend the function $\tilde U$,
introduced in assumption \ref{a:U3l}, to the whole $(a,b)\times \RR^k\minus \Delta$
in the natural way 
\begin{equation}
\label{eq:tildeUlog}
\tilde U(t,x) = \tilde U(t,s) - M(t)\log|x|,
\end{equation}
where $M$ has been introduced in \eqref{eq:M}.
\begin{theorem}\label{theo:nocoll_log}
In addition to \ref{a:U0}, \ref{a:U1}, \ref{a:U2l}, \ref{a:U3l}, \ref{a:U4l}, \ref{a:U5} 
assume the potential $U$ to be of the form
\begin{equation}
\label{eq:hp_pot_log}
U(t,x) = \tilde U(t,x) + W(t,x)
\end{equation}
where $\tilde U$ satisfies \eqref{eq:tildeUlog}
and $W$ is a bounded $\cont^1$ function on $(a,b) \times \RR^{k}$.
Furthermore assume that, for a given $\xi\in\Delta$ , $\tilde U$ satisfies \ref{a:U6} and
\begin{assiomi}
\item[\usetteltag]
\label{a:U7l}
for every $x\in\RR^k$ and $t \in (a,b)$ there holds
\[
\tilde U(t,x) = -\frac12 M(t)\log\left(|\pi_Wx|^2 + \psi^2(\pi_{W^\perp}x)\right)
\]
where $\pi_W$ and $\pi_{W^\perp}$ denote the orthogonal projections onto $W$ and $W^\perp$,
$\psi$ is $\cont^1$ and 
homogeneous of degree $1$.
% $\psi(0)=0$.
\end{assiomi}
Then locally minimizing trajectories do not have collisions at the configuration $\xi$ at the time $t^*$.
\end{theorem}

\begin{proof}
As in the proof of Theorem~\ref{theo:nocoll},
we consider a generalized solution $\bar x$ and we first reduce to 
the case of an isolated total collision at the origin occurring at the time $t=0$.
From Theorem~\ref{central_conf_log} we deduce the existence of $\delta_0>0$
such that no other collision occur in $[-\delta_0,\delta_0]$, hence 
we perform a local variation on the trajectory of $\bar x$
that removes the collision and makes the action decrease.

% To this aim, we first observe that,
% from Theorem~\ref{central_conf_log}, it follows that 
% $r(t)=|\bar x(t)|$ is asymptotic to $t \sqrt{-2M_0 \log t}$ as $t$ tends to $0^+$; 
% since $t \sqrt{-2M_0 \log t}>t$ whenever $t<e^{-1/2M_0}$, we can choose $\delta_0$ small enough
% and $\eps <\delta_0/2$, to have
% \begin{equation}\label{eq:>eps}
% r(t) > \eps \qquad \mbox{for every } t \in [\delta_0-\eps,\delta_0].
% \end{equation}
% Fixed $\eps >0$ verifying \eqref{eq:>eps}
% and $\delta \in \RR^d$ such that $|\delta|=\eps$, 
Consider now the 
{standard variation} $v^\delta$, defined at page \pageref{def:stand_var},
on the interval $[0,\delta_0]$ (i.e., in Definition \ref{def:stand_var} 
$T$ is replaced by $\delta_0$).
% \begin{equation}
% \label{eq:variation_mu}
% v^{\delta,\eps}:[0,\delta_0]\to\RR^{k}, \qquad
% v^{\delta,\eps}(t) := 
% \left\{
% \begin{array}{ll}
% \xi,                                        & \mbox{if } t \in [0,\delta_0-\eps]  \\
% \displaystyle(\delta_0-t)\frac{\xi}{\eps},  & \mbox{if } t \in [\delta_0-\eps,\delta_0]
% \end{array}
% \right.
% \end{equation}
% where $\xi = (\eps\delta,0,\ldots,0) \in \RR^{k}$.
%
Let $\Delta^{\delta} {\cal A}$ denote the difference 
$$
\Delta^{\delta} {\cal A} :=
{\cal A}(\bar x + v^{\delta},[0,\delta_0]) - {\cal A}(\bar x,[0,\delta_0]);
$$
generally speaking, this difference can be positive or negative, depending
on the choice of $\delta$.
Our goal is to prove that, when \emph{averaging over a suitable set of standard 
variations}, the action lowers.
Hence  $\Delta^{\delta} {\cal A}$
must be negative for at least one choise of $\delta$ and the path $\bar x$ 
can not be a local minimizer for the action.

We can write $\Delta^{\delta} {\cal A}$ as the sum of three terms
\begin{equation}
\label{eq:A_var_mu}
\Delta^{\delta} {\cal A} = \int_{0}^{\delta_0} \Delta^{\delta}{\cal K}(t) \,dt 
                   + \int_0^{\delta_0} \Delta^{\delta}{\cal U}(t) \,dt
                   + \int_{0}^{\delta_0} \Delta^{\delta}{\cal W}(t) \,dt 
\end{equation}
where $\Delta^{\delta}{\cal K}(t)$, $\Delta^{\delta}{\cal U}(t)$ and
$\Delta^{\delta}{\cal W}(t)$ are respectively the variations of the kinetic energy,
of the singular potential $\tilde U$ and of the smooth part of the potential, $W$.
More precisely  since the first derivative of the function $v^{\delta}$
vanishes everywhere on $[0,\delta_0]$, except on $[\delta_0-|\delta|,\delta_0]$, we compute
\begin{equation}
\label{eq:K_var_mu}
\Delta^{\delta}{\cal K}(t) :=
\left\{
\begin{array}{ll}
\displaystyle
0, & \mbox{if } t \in [0,\delta_0-|\delta|], \\
\displaystyle\frac{1}{2}( \left| \dot{\bar x} - \delta/|\delta|\right|^2 - |\dot{\bar x}|^2 )
   & \mbox{if } t \in [\delta_0-|\delta|,\delta_0].
\end{array}
\right.
\end{equation}
Similarly
\[
\Delta^{\delta}{\cal U}(t) := \tilde U(t,\bar x+v^{\delta}) - \tilde U(t,\bar x)
\qquad \mbox{and}
\qquad \Delta^{\delta}{\cal W}(t) := W(t,\bar x+v^{\delta}) - W(t,\bar x).
\]
We now evaluate separately the mean values 
of the tree terms of $\Delta^{\delta} {\cal A}$ over the circle $S^{|\delta|}$ of radius $|\delta|$ in $W$.
\begin{lemma}
\label{le:media_cin}
There holds
\begin{equation}
\frac{1}{2\pi|\delta|}\int_{S^{|\delta|}} \int_0^{\delta_0}
(\Delta^{\delta}{\cal K} + \Delta^{\delta} {\cal W})\,dt\,d\delta = O(|\delta|).
\end{equation}
\end{lemma}
\begin{proof}
From~\eqref{eq:K_var_mu} we obtain
\begin{equation*}
\int_0^{\delta_0} \Delta^\delta{\cal K}(t) dt =
\int_{\delta_0-|\delta|}^{\delta_0}
\frac{1}{2}(\left| \dot{\bar x} -\delta/|\delta|\right|^2 - |\dot{\bar x}|^2) dt
= \frac{1}{2}\left(|\delta|-2\int_{\delta_0-|\delta|}^{\delta_0}
\dot{\bar x}(t)\cdot\frac{\delta}{|\delta|} dt \right),
\end{equation*}
hence 
\begin{equation*}
\left|\int_0^{\delta_0} \Delta^\delta{\cal K}(t) dt \right| \leq O(|\delta|),
\end{equation*}
which does not depend on the circle $S^{|\delta|}$ where $\delta$ varies. 
Concerning the variation of the $\cont^1$ function $W$ we have
\begin{equation*}
\begin{split}
\left| \int_0^{\delta_0} \Delta^\delta{\cal W}(t) \,dt \right| & =
\left| \int_0^{\delta_0-|\delta|} \Delta^\delta{\cal W}(t) \,dt \right| +
\left| \int_{\delta_0-|\delta|}^{\delta_0} \Delta^\delta{\cal W}(t) \right| \,dt \\
& \leq W_1 |\delta| (\delta_0-|\delta|) + 2W_2 |\delta| = O(|\delta|),
\end{split}
\end{equation*}
where $W_1$ is a bound for 
$\left|\frac{\partial W}{\partial x}(t,\bar x+\lambda v^{\delta})\right|$, 
with $\lambda \in [0,1]$ and $t \in [0,\delta_0-|\delta|]$ while $W_2$ is an upper bound for
$|W(t,x)|$.
\end{proof}

In order to estimates the variation of the potential part,
$\Delta^\delta{\cal U}(t)$, we prove the next two technical lemmata.
Let us start with recalling an equivalent version of the mean value property for the
fundamental solution of the planar Laplace equation.

\begin{lemma}
\label{le:pot:2}
Fixed $z>0$, for every $y\in\RR$ such that $y\geq 2z$, we have 
\[
\frac{1}{2\pi }\int_{0}^{2\pi} \log(y+2z\cos\vartheta) \,d\vartheta =
\log\frac{y+\sqrt{y^2-4z^2}}{2}.
\]
\end{lemma}
\begin{proof}
Since $y\geq 2z$, then $\frac{y+\sqrt{y^2-4z^2}}{2}\geq z$. Let $x\in\RR^2$ be such
that 
$|x|=\frac{y+\sqrt{y^2-4z^2}}{2}$, then $y=(|x|^2+z^2)/|x|$ and for every $\delta \in S^z$, 
where $S^z$ is the circle of radius $z$, we have
\[
\begin{split}
|x+\delta|^2 &= |x|^2 + z^2 +2z|x|\cos \vartheta \\
             &= |x|\left( \frac{|x|^2 + z^2}{|x|} + 2z\cos \vartheta\right) = 
                |x|(y + 2z\cos \vartheta).
\end{split}
\]
We have, as the logarithm is the fundamental solution to the Laplace equation on the plane,
\begin{equation}
\label{le:soluz_fond}
\frac{1}{2\pi z}\int_{S^z} \log|x+\delta|^2 \,d\delta = \max\{\log|x|^2, \log z^2\} = 
\left\{
\begin{array}{ll}
\log|x|^2 ,        & \mbox{if }|x| >    z \\
\log z^2 ,  & \mbox{if }|x| \leq z.
\end{array}
\right..
\end{equation}
Consequently, when computing
\[
\begin{split}
\int_{S^z}\log|x+\delta|^2 d\delta &= \int_{S^z}\log|x| d\delta +
z\int_0^{2\pi}\log(y + 2z\cos \vartheta)d\vartheta \\
&= 2\pi z\log|x| + z\int_0^{2\pi}\log(y + 2z\cos \vartheta)d\vartheta,
\end{split}
\]
we find
\[
2\pi z \log|x|^2 = 2\pi z\log|x|+z\int_0^{2\pi}\log(y + 2z\cos \vartheta)d\vartheta.
\]
We conclude replacing $|x|=\frac{y+\sqrt{y^2-4z^2}}{2}$.
\end{proof}

Now we consider the averages of the potential with respect to a circle in $W$ (here we assume implicitly that $d\geq 3$).

\begin{lemma}
\label{le:pot:3}
Fixed $|\delta|>0$, for every circle of radius $|\delta|$, $S^{|\delta|} \subset W$, 
for every $x\in\RR^d$ and every $t \in [0,\delta_0]$, there holds
\[
\frac{1}{2\pi|\delta|}\int_{S^{|\delta|}} \left(\tilde U(x+\delta)-\tilde U(x)\right)\,d\delta  \leq 
% \left{
% \begin{array}{l}
\begin{cases}
0    \text{\ \ ( if $\ds |\pi_W x|^2+\psi^2(\pi_{W^\perp} x) > |\delta|^2$)},\\
\frac{M(t)}{2}\log(|\pi_W x|^2+\psi^2(\pi_{W^\perp} x))-\log(|\delta|^2) \text{\ \ (otherwise).}
\end{cases}
% \end{array}
% \right.
\]
\end{lemma}
\begin{proof}
We consider the orthogonal decomposition of $x$, 
$x=\pi_W x+\pi_{W^\perp} x$, and we term $u:=|\pi_W x|$ and $\eps:=\psi(\pi_{W^\perp} x)$.
Since whenever $\delta \in W$ we have
\[
|\pi_W (x+\delta)|^2+\psi^2(\pi_{W^\perp} x) = u^2  + |\delta|^2 +2u|\delta|\cos \vartheta + \eps^2\geq 0,
\]
when $\cos \vartheta = -1$ we have $\frac{u^2  + |\delta|^2 + \eps^2}{u|\delta|} \geq 2$
and, using Lemma \ref{le:pot:2} and equation \eqref{le:soluz_fond}, we compute
\[
\begin{split}
\frac{1}{2\pi |\delta|}&\int_{S^{|\delta|}} \log(|\pi_W (x+\delta)|^2+\psi^2(\pi_{W^\perp} x)) \,d\delta \\
& =\frac{1}{2\pi}\int_0^{2\pi } \log(u^2+\eps^2+|\delta|^2+2u|\delta|\cos \vartheta)
\,d\vartheta\\
& =\frac{1}{2\pi}\int_0^{2\pi} \log\left(\frac{u^2 + \eps^2 + |\delta|^2}{u|\delta|} + 2\cos
\vartheta\right) \,d\vartheta + \log (u|\delta|)\\
& =\log\left(\frac{u^2 + \eps^2 + |\delta|^2 + \sqrt{(u^2 + \eps^2 +
|\delta|^2)^2-4u^2|\delta|^2}}{2}\right)\\
& \geq\log\left(\frac{u^2 + \eps^2 + |\delta|^2 + \sqrt{(u^2 + \eps^2 + |\delta|^2)^2 - 4u^2|\delta|^2 -
4\eps^2|\delta|^2}}{2}\right)\\
& = \log \left(\frac{u^2 + \eps^2 + |\delta|^2 + |u^2 + \eps^2 - |\delta|^2|}{2}\right)\\
& = \max\left( \log(|\pi_W x|^2+\psi^2(\pi_{W^\perp} x),\log(|\delta|^2)\right)
\end{split}
\]
and the assertion easily follows.
\end{proof}

\begin{lemma}
\label{le:media_pot1}
Let $S$ be the circle of radius $|\delta|$ on $W$;
then, as $|\delta| \rightarrow 0$
\begin{equation}
\frac{1}{2\pi|\delta|}\int_{S^{|\delta|}} \int_0^{\delta_0} \Delta^{\delta} {\cal U} \, dt \, d\delta 
< -K|\delta|\sqrt{-\log|\delta|}, \qquad K>0.
\end{equation}
\end{lemma}

\begin{proof}
Let $S^{|\delta|}$ be the circle of radius $|\delta|$ on $W$,
we apply Fubini-Tonelli's Theorem and we argue as in the proof of 
Lemma~\ref{le:pot:3} to have
\[
\begin{split}
& \frac{1}{2\pi|\delta|}\int_{S^{|\delta|}} \int_{0}^{\delta_0} \Delta^{\delta}{\cal U}(t) dt \,  d\delta
= \int_{0}^{\delta_0} \frac{1}{2\pi|\delta|} 
\int_{S^{|\delta|}} \tilde U(\bar x+v^\delta) - \tilde U(\bar x) \, d\delta dt\\
&= \frac{M^*}{2} \int_{0}^{\delta_0} 
\left\{-\max\left[ \log(|\pi_W\bar x|^2 + \psi^2(\pi_{W^\perp}\bar x)), \log |v^\delta|^2\right]
+ \log(|\pi_W\bar x|^2 + \psi^2(\pi_{W^\perp}\bar x))\right\}\, dt
\end{split}
\]
where $M^* = \max_t |M(t)|$. We then straightforwardly deduce that, for every $S^{|\delta|}\subset W$
\[
\frac{1}{2\pi|\delta|}\int_{S^{|\delta|}} \int_{0}^{\delta_0} \Delta^{\delta}{\cal U}(t) 
dt \,  d\delta < 0.
\]
In order to estimate more precisely this quantity, we observe that 
\begin{equation}\label{eq:le:media_pot1}
\int_{0}^{\delta_0} \frac{1}{2\pi|v^\delta|} 
\int_{S^{|v^\delta|}} \tilde U(\bar x+v^\delta) - \tilde U(\bar x) \, d\delta dt
\leq
\int_A 
\log\frac{|\pi_W\bar x|^2 + \psi^2(\pi_{W^\perp}\bar x)}{|\delta|^2} \, dt%d\delta dt
\end{equation}
where
\[
A := \left\{ t \in [0,\delta_0-|\delta|] : |\pi_W\bar x|^2 + \psi^2(\pi_{W^\perp}\bar x) < |\delta|^2\right\}.
\]
Furthermore, there exists  a strictly positive constant $C$ such that
\[
C r^2 < |\pi_W x|^2 + \psi^2(\pi_{W^\perp} x) < C^{-1} r^2
\]
where, as usual, we denote $r^2=|\pi_W x|^2 + |\pi_{W^\perp} x|^2$ the radius of $x$. 
The left inequality follows from Theorem~\ref{central_conf_log} indeed the existence of
a finite limit of $\tilde U(t,s(t))$ prevents the projection $|\pi_W x|^2$
and the function $\psi^2(\pi_{W^\perp} x)$ to be both infinitesimal with  $r^2$.
The right inequality follows from the continuity of $\psi$.
From \eqref{eq:le:media_pot1} and the asymptotic estimates of Theorem~\ref{central_conf_log} 
we conclude that, as $|\delta| \rightarrow 0$
\[
\begin{split}
\frac{1}{2\pi|\delta|}\int_{S^{|\delta|}} \int_{0}^{\delta_0} \Delta^{\delta}{\cal U}(t) 
dt \,  d\delta 
& \leq \int_{_{t:r(t)<|\delta|/\sqrt{C}}}
\log\frac{r^2(t)}{C|\delta|^2} \, dt \\
& \sim \int_0^{|\delta|/\sqrt{C}}
2\frac{\log (r/\sqrt{C}|\delta|)}{-\sqrt{-\log r}}  dr \\
% & < \int_0^{|\delta|/\sqrt{C}}
% 2\frac{\log r}{-\sqrt{-\log r}}  dr \\
& < - 2\int_0^{|\delta|/\sqrt{C}} \sqrt{-\log r} dr
<  - K|\delta| \sqrt{-\log |\delta|}
\end{split}
\]
for some positive $K$, since $-\sqrt{-\log r}$ is an increasing function on the interval
$[0,|\delta|]$.
\end{proof}

\noindent\emph{End of the Proof of Theorem \ref{theo:nocoll}.}
Let $S^{|\delta|}$ be a circle in $W$ with radius $|\delta|$ and $\Delta^{\delta} {\cal A}$
the variation of the action functional defined in~(\ref{eq:A_var_mu}), then 
from Lemmata \ref{le:media_cin} and  \ref{le:media_pot1}
we conclude that, as $|\delta|$ tends to $0$
\[
\frac{1}{2\pi|\delta|}\int_{S^{|\delta|}} \Delta^{\delta} {\cal A} d\delta \leq O(|\delta|) 
- K|\delta| \sqrt{-\log |\delta|} < 0.
\]
\end{proof}

%===========================%
% DIM. PROPRIETA' LOGARITMO %
%===========================%
% \begin{remark}
% From Theorem \ref{thm:var<0_mu} we easily deduce that, as $\eps$ tends to $0$,
% \begin{equation*}
% \frac{1}{4\pi}\int_{S^2} \Delta^{\delta,\eps} {\cal A} \, d\delta < 0,
% \end{equation*}
% where $S^2 = \{ z  \in \RR^3 : |z|=1 \}$; indeed
% \begin{equation*}
% \frac{1}{4\pi}\int_{S^2} \Delta^{\delta,\eps} {\cal A} \, d\delta 
% = \frac{1}{4\pi} \int_0^\pi \sin \varphi \int_{S^1} \Delta^{\delta,\eps} {\cal A}
% \, d\delta \, d\varphi 
% \leq \frac{1}{4} \int_{S^1} \Delta^{\delta,\eps} {\cal A} \, d\delta < 0.
% \end{equation*}
% \end{remark}

Of course, likewise to Theorem~\ref{theo:nocoll_2}, there holds

\begin{theorem}\label{theo:nocoll_log_2}
In addition to \ref{a:U0}, \ref{a:U1}, \ref{a:U2l}, \ref{a:U3l}, \ref{a:U4l}, \ref{a:U5}, assume $\tilde U$ be of the form
\[
\tilde U(x)=-\sum_{\nu=1}^N K_\nu\log\left({\rm dist}(x,V_\nu)\right)
\]
where $K_\nu$ are positive constants and $V_\nu$ is a family of linear subspaces, with ${\rm codim}(V_\nu)\geq 2$, 
for every $\nu=1,\dots,N$.
Then locally minimizing trajectories do not have collisions at the time $t^*$.
\end{theorem}

%================================================================
\subsection{Neumann boundary conditions and $G$--equivariant minimizers}
%================================================================

As a final comment of this Section, we remark that, in our framework, the analysis allows 
to prove that minimizers to the fixed--ends (Bolza) problems are free of collisions: indeed all 
the variations of our class have compact support. However,
other type  of boundary conditions (generalized Neumann) can be treated in the same way. 
Indeed, consider a trajectory which is a (local) minimizer of the action among
all paths satisfying the boundary conditions
\[ x(0)\in X^0\qquad x(T)\in X^1,\]
where $X^0$ and $X^1$ are two given linear subspaces of  the configuration
space. Consider a (locally) minimizing  path $\bar x$: of course it has not
interior collisions.  In order to exclude boundary collisions we have to be
sure that the class of variations preserve the boundary condition; this can be
achieved by restricting to $X^i$ the points $\delta$ appearing in the standard
variations.  Hence, to complete the averaging argument, one needs
assumptions~\ref{a:U6} and \ref{a:U7h} or \ref{a:U7l}  to be fulfilled also by
the restriction of the potential to the boundary subspaces $X^i$.  This point
of view differs from that of \cite{chen}, where the boundary subspaces can not
be chosen arbitrarily in the configuration space. The argument in \cite{chen},
already introduced in \cite{FT}, does not involve any averaging on the boundary
but relies upon a suitable choice of a standard variation whose projection is
extremal. 

The  analysis of boundary conditions was a key point in the paper \cite{FT},
where symmetric periodic trajectories where constructed by reflections about
given subspaces. 
By Theorems~\ref{theo:nocoll} and~\ref{theo:nocoll_log}  one can obtain the
absence of collisions also for $G$--equivariant (local) minimizers, provided
the group $G$ satisfies the {\em Rotating Circle Property} introduced in
\cite{FT} (see Example~\ref{exe:gruppo1}). Hence, existence of $G$--equivariant
collisionless periodic solutions can be proved for the wide class of symmetry
groups described in \cite{FT,F3,BFT}, for a much larger class of interacting
potentials, 
including quasi--homogeneous and logarithmic ones. On the other hand,
Theorems~\ref{theo:nocoll_2} and~\ref{theo:nocoll_log_2} can be applied to
prove that $G$--equivariant minimals are collisionless for many relevant
symmetry groups violating the rotating circle property, such as the groups of
rotations in \cite{F2}; indeed, the idea of averaging on spheres having maximal
dimension has been borrowed from that paper (cf.~Example~\ref{exe:gruppo2}).

%The crucial assumptions for these extensions are \uc\, and \usei\enspace concerning the interaction 
%between each pair of  particles, and~\eqref{eq:strucDelta} on the structure of the singular set.
%Indeed they allow the reduction from partial collisions to total ones, in order to avoid 
%each collision from the trajectory of an equivariant minimizer.

%----------------------------------
\section{Examples and further remarks}\label{sec:comments}
%----------------------------------

We now discuss various examples of classes of potentials which fullfill our assumptions.

\begin{example}[Homogeneous isotropic potentials]\label{ex:0}
The simplest example of function satisfying  all our assumptions 
\ref{a:U0},
\ref{a:U1},
\ref{a:U2h},
\ref{a:U3h},
\ref{a:U4h},
\ref{a:U5},
\ref{a:U6} and 
\ref{a:U7h}
is the $\alpha$-homogeneous one-center problem:
\[
U_\alpha(x) =  \frac{1}{|x|^\alpha},
\]
and its associated $n$--body problem:
\[
U_\alpha(x) = \sum_{\substack{i<j\\  i,j=1}}^n \frac{m_i m_j}{|x_i-x_j|^\alpha}.
\]
Assumptions \ref{a:U0} and \ref{a:U1} are trivially satisfied since $U$ is positive, diverges to 
$+\infty$ when $x$ approaches $\Delta = \{x \in \RR^{nd}: x_i=x_j \mbox{ for some }i\neq j\}$,
and does not depend on time.
Furthermore in both \ref{a:U2} and \ref{a:U2h} the equality is achieved
with $\tilde \alpha=\alpha$ and $C_2=0$. Since $U$ is homogeneous of degree $-\alpha$,
in \ref{a:U3h} and \ref{a:U4h} the function $\tilde U$ coincides with $U$. \ref{a:U5} and \ref{a:U6} are trivially satisfied, while~\ref{a:U7h} holds by virtue of Proposition~\ref{prop:usette}
\end{example}

\begin{example}[Logarithmic potentials]\label{ex:2}
Our results apply also to logarithmic singularities of type
\[
U_{log}(x) = \sum_{\substack{i<j\\ i,j=1}}^{n} m_i m_j\log\frac{1}{|x_i-x_j|};
\]
indeed \ref{a:U2} is in this case satisfied for every value of\enspace $\tilde \alpha$
and \ref{a:U2l}, \ref{a:U3l} and \ref{a:U4l} are verified with $C_2=0$.

Dynamical systems of type~\eqref{DS} with logarithmic interactions arise in the study of vortex flows in fluid mechanics, and,
precisely, in the analysis of systems of $n$ 
almost--parallel vortex filaments, under a linearized version of the LIA 
self--interaction assumption (see \cite{KPV, KMD}).
\end{example}

\begin{example}[Anisotropic $n$--body potentials] Consider potentials having the form
\[
U(t,x) = \sum_{\substack{i<j\\ i,j=1}}^{n} U_{i,j}(t,x_i-x_j),
\]
where the interaction potentials $U_{i,j}$ have a singularity at zero, of
homogeneous or logarithmic type, but do depend on the angle. Typical examples
are the Gutzwiller potentials \cite{Gu}. Notice that the total potential
satisfies assumptions \ref{a:U0}, \ref{a:U1}, 
and \ref{a:U2h}, \ref{a:U3h}, \ref{a:U4h} (or \ref{a:U2l}, \ref{a:U3l}, \ref{a:U4l}) provided each of the $U_{i,j}$'s
do. It not difficult to see that also~\eqref{eq:strucDelta} and \ref{a:U5}
hold (in the $n$--body case), while \ref{a:U6} and \ref{a:U7h} or \ref{a:U7l}
do not. Hence we can not exclude the presence of collisions for locally
minimizing paths, though the results about isolatedness and the asymptotic
estimates are still available. More generally, we can deal with potentials of
the form
\[
U_\alpha(r s) = r^{-\alpha}\tilde U(s),
\]
where $\tilde U:{\mathcal E}\setminus\Delta \to \RR$ is positive and admits an arbitrary
 singular set on the ellipsoid ${\mathcal E}=\{I=1\}$, provided 
 \[\lim_{s\to\Delta}\tilde U(s)=+\infty\;.\]
It is worthwhile noticing that as a consequence of Theorem~\ref{central_conf}, a total collision
trajectory will not interact, definitively, with the singularities of $\tilde U$.
\end{example}

The class of potentials satisfying our assumptions is clearly stable with respect to 
the addition of arbitrary perturbations of class $\mathcal C^1$. Therefore, we are mainly 
interested in the analysis of those perturbations which are singular themselves. 

\begin{example}[$N$--body potentials with time--varying masses]\label{ex:0.5}
Although the potentials in the previous examples  do not depend on time, our assumptions
allow an effective time--dependence of the potentials.  For instance, we can choose positive 
and bounded $\cont^1$ functions $m_i(t)$, $i=1,\ldots,n$.

Obviously, the simplest example is the class of $\alpha$-homogeneous $n$-body problem
\[
U_\alpha(t,x) = \sum_{\substack{i<j\\ i,j=1}}^{n} \frac{m_i(t) m_j(t)}{|x_i-x_j|^\alpha},\qquad 0<\alpha<2.
\]
Assumptions \ref{a:U0} and \ref{a:U1} are trivially satisfied since $U$ is positive, 
diverges to  $+\infty$ when $x$ approaches 
$\Delta = \{x \in \RR^{nd}: x_i=x_j \mbox{ for some }i\neq j\}$,
and does not depend on time.
Furthermore in both \ref{a:U2} and \ref{a:U2h} the equality is achieved
with $\tilde \alpha=\alpha$ and $C_2=0$. Since $U$ is homogeneous of degree $-\alpha$,
in \ref{a:U3h} and \ref{a:U4h} the function $\tilde U$ coincides with $U$.
\end{example}

\begin{example}[Quasi--homogeneous potentials]\label{ex:1}
We can also handle homogeneous perturbations of degree $-\beta$ of the potential $U_\alpha$
\[
U(x) = U_\alpha(x) + \lambda U_\beta(x) \qquad\qquad 0<\beta<\alpha<2.
\]
Indeed, when $\lambda>0$ condition \ref{a:U2h} is verified (with the strict
inequality) with $\gamma = C_2 =0$, while, when $\lambda<0$, then \ref{a:U2} holds, when $|x|$
is sufficiently small, with $C_2=\alpha-\beta$ and $0<\gamma < \alpha-\beta$.

As pointed out in \cite{diacuab} (where the case $\beta=1$ and $\alpha>1$ was treated), quasi--homogeneous potentials generalize classical potentials such as Newton, Coulomb, Birkhoff, Manev and many others. Therefore, the range of physical applications
of quasi--homogeneous potentials spans from celestial mechanics and atomic physics to chemistry and crystallography. It is worthwhile noticing that the collision problem for quasi--homogeneous potentials exhibit an interesting and peculiar lack of regularity. Indeed, a classical framework for the study of collisions is given by the McGehee coordinates \cite{mcgehee1}
(here and below we assume, for simplicity of notations, all the masses be equal to one):
\begin{equation*}
\begin{split}
r&=|x|=I^{1/2}\\
s&=\dfrac{x}{r}\\
v&=r^{\alpha/2}(y\cdot s)\\
u&=r^{\alpha/2}(y-y\cdot s s).
\end{split}
\end{equation*}
After a reparametrization of the time--variable (see~\ref{eq:mcgehee1}):
\begin{equation}\label{eq:mcgehee2}
d\tau=r^{-1-\alpha/2}dt,
\end{equation}
the equation of motions become (here $^\prime$ denotes differentiation with respect to the new time variable $\tau$):
\begin{equation*}
\begin{split}
r^\prime&=rv\\
v^\prime&=\dfrac{\alpha}{2}v^2+|u|^2-r^{\alpha-\beta}\lambda U_\beta(s)-\alpha U_\alpha(s)\\
s^\prime&=u\\
u^\prime&=\left(\dfrac{\alpha}{2}-1\right)vu-|u|^2s+r^{\alpha-\beta}\lambda\left(U_\beta(s)s-\nabla U_\beta(s)\right)+\alpha U_\alpha(s)s+\nabla U_\alpha(s)\;.\\
\end{split}
\end{equation*}
The field depends on $r$ in a non smooth manner, unless $\alpha-\beta\geq 1$ (this last condition was indeed assumed in \cite{diacuab}). Hence the flow \emph{can not be continuously extended} to the total collision manifold $C=\{(r,s,v,u)\;:\; r=0,\,\frac{1}{2}(|u|^2+v^2)-2U_\alpha=0\}$. Another peculiar feature of this system is that the monotonicity of the variable $v$ can not be ensured close to the collision manifold. As a consequence, the usual analysis of collision and near collision motions can not be extended to this case.
\end{example}

\begin{example}[$N$--body potential reduced by a symmetry group satisfying the rotating circle property] 
The paper \cite{terven} deals with \label{exe:gruppo1}
minimal trajectories to the spatial $2N$--body problem under the \emph{hip--hop symmetry}, 
where the configuration is constrained at all time to form a regular antiprism.   
This problem has three degrees of freedom and the reduced potential of a 
configuration generated by the  point of coordinates $(u,\zeta)\in\CC\times\RR\simeq\RR^3$
decomposes as 
$$
U(u,\zeta)=\frac{K(N)}{|u|^\alpha}+U_0(u,\zeta),
$$
where 
\begin{eqnarray}
K(N)&=&\sum\limits_{k=1}^{N-1}\frac{1}{\sin^\alpha(\frac{k\pi}{N})}, \nonumber
\\
U_0(u,\zeta)&=& 
\sum\limits_{k=1}^{N} 
\frac{1}{\left(\sin^2\left(\frac{(2k-1)\pi}{2N}\right)|u|^2+\zeta^2\right)^{\frac{\alpha}{2}}}, \nonumber
\end{eqnarray}
The first term  comes from the interaction among points of the same $N$--agon and is 
singular at simultaneous partial
collisions on the $\zeta$--axis. The second term, $U_0(u,\zeta)$, comes from the interaction 
between the the upper and lower $N$--agons and is singular only at the origin.
One easily verifies that all the assumptions are satisfied, including, again by
Proposition~\ref{prop:usette}, \ref{a:U6} and \ref{a:U7h}. In general, one
easily verifies that, for a given symmetry group $G$ of the $N$--body problem,
it is equivalent to say that $\ker\tau$ has the rotating circle property and
that the reduced potential verifies \ref{a:U6} and \ref{a:U7h}.
\end{example}

\begin{example}[$N$--body potential reduced by a symmetry group not satisfying the rotating 
circle property] Consider the \label{exe:gruppo2} symmetry groups generated by rotations 
introduced in \cite{F2}: the configuration is, at all time,
an orbit of a group $Y$ of rotations about given lines in the $3$--dimensional space. 
When $Y$ is a finite group, the reduced potential takes the form required in 
Theorem~\ref{theo:nocoll_2} and minimizers can be shown to be free of collision. 
\end{example}

%  \bibliographystyle{plain}
%  \bibliography{bibliog}
% \end{document}

\end{document}